\documentclass[12pt]{amsart}
\usepackage{amsbsy,amsfonts}
\usepackage{latexsym,eucal,amsmath,amsthm,amsxtra}
\usepackage{amssymb}
\newcommand{\db}{\mathbb }

\newcommand{\dbR}{{\db R}}

\newcommand{\dbZ}{{\db Z}}
\newcommand{\dbN}{{\db N}}

\newcommand{\F}{\mbox{{\em F}}}

\newcommand{\om}{\omega}
\newcommand{\m}{\mu}

\newcommand{\x}{\xi}

\newcommand{\half}{\frac{1}{2}}
\def\gi{\frac{|\m|}{|\x|}}
\def\gia{\frac{|\m_{1}|}{|\x_{1}|}}
\def\gib{\frac{|\m_{2}|}{|\x_{2}|}}

\newcommand{\lam}{\lambda}

\def\R{{\db R}}

\newtheorem{theorem}{Theorem}
\theoremstyle{definition}
\newtheorem{definition}{Definition}
\theoremstyle{remark}
\newtheorem{remark}{Remark}
\theoremstyle{proposition}
\newtheorem{proposition}{Proposition}
\theoremstyle{lemma}
\newtheorem{lemma}{Lemma}
\theoremstyle{corollary}
\newtheorem{corollary}{Corollary}

\numberwithin{equation}{section}
\numberwithin{lemma}{section}
\numberwithin{remark}{section}

\begin{document}
\date{}
\title{Low Regularity Solutions for the Kadomtsev-Petviashvili I
Equation}
\author{J. Colliander}
\thanks{J.E.C. was supported in part by  N.S.F. Grant DMS 0100595 and N.S.E.R.C. Grant RGPIN 250233-03.}
\address{\small University of Toronto}

\author{C. Kenig}
\thanks{C.K. was supported in part by N.S.F Grant DMS 9500725}
\address{\small University of Chicago}

\author{G. Staffilani}
\thanks{G.S. was supported in part by N.S.F. Grant DMS 0100375
and a grant by the Sloan Foundation.}
\address{\small Massachusetts Institute of Technology}

\begin{abstract}
In this paper we obtain local in time existence and (suitable)
uniqueness and
continuous dependence for  the KP-I equation for
small
data in the intersection of the energy space and
a natural weighted
$L^{2}$ space.
\end{abstract}

\maketitle
\section{Introduction}
We consider the  KP initial value problem (IVP)
\begin{equation}
\left\{ \begin{array}{l}
\partial_t u + \partial_x^3 u + \gamma\partial_x^{-1}
\partial_y^2u
+\beta\partial_x (u^2) = 0, \\
 u(x,0) = u_0(x) \hspace{1.5cm}(x,y) \in \dbR^2,
\,  t \in \dbR,
\end{array}\right.
\label{ivp}\end{equation}
where $u=u(t,x,y)$ is a scalar unknown function,
$\beta\ne 0$ and $\gamma\ne 0$ are real constant.
If $\gamma<0$ the IVP (\ref{ivp}) is called KP-I
and if $\gamma>0$ it
takes the name KP-II. These equations model \cite{KP} the propagation along the $x$-axis of nonlinear
dispersive long waves on the surface of a fluid with a slow variation along the $y$-axis. They also
arise as universal models in wave propagation and may be considered as two dimensional generalizations
of the Korteweg-de Vries equation.

The first result regarding well-posedness for a KP type equation 
is due to Ukai \cite{U}. He uses a standard energy
method that does not recognize the type I or II of the equation. His
result provides local well-posedness for initial data and their
antiderivatives  in $ H^{s}, \, s\geq 3$.
Faminskii \cite{F}
observed a better smoothing effect in the KP-II evolution and used
this to prove well-posedness results.
Bourgain performed a Fourier analysis \cite{B}
of the term $\partial_x (u^2)$
in the KP-II equation in which the derivative is recovered in a nonlinear
way. The result obtained gave local well-posedness
of KP-II for initial data in $L^2 $. Since the $L^2$ norm is conserved
during the KP-II evolution, the $L^2$ local result may be iterated to
prove global well-posedness.
Takaoka \cite{Ta1} and Takaoka and Tzvetkov \cite{TT}
improved Bourgain's result by proving local well-posedness in an
anisotropic Sobolev space{\footnote{Here $H^{s_1, s_2}_{x,y} = \{ f
    \in \mathcal{S}': \int \int | \langle \xi \rangle^{s_1} \langle
    \mu \rangle^{s_2} \widehat{f} (\xi , \mu )|^2 d\xi d \mu < \infty \}$
    where $\langle \cdot \rangle = (1 + | \cdot |^2 )^\half$ and
    $\sphat$ denotes the Fourier transform from the spatial variables
    $(x,y)$ to their dual variables $(\xi, \mu)$.}} $H^{-\frac{1}{3}+\epsilon,0}_{x,y}$.
For the KP-I equation the situation is more delicate. There are
several results on local and global existence of solutions, but not a
satisfactory well-posedness theory for data with no more than
two derivatives in $L^2$. Fokas and Sung \cite{FS}, and Zhou \cite{Z},
obtained  global
existence  for small data via inverse scattering
techniques.
Schwarz \cite{Sw} proved existence of weak global
periodic solutions with small $L^{2}$ data. The smallness condition was
subsequently removed \cite{C}. Tom \cite{T}
proved existence of global weak solutions
for initial data in $H^{1}$ together with their antiderivative.
For well-posedness results, we recall the work
of Saut \cite{S}, Isaza, Mej\'ia and Stallbohom~\cite{IMS} and
finally
the work of I\'orio and Nunes~\cite{IN}. The last two authors
use the quasi linear theory of Kato, together with parabolic
regularization, to prove local well-posedness with data and their
antiderivatives in $H^{s}, s>2$. The limitation $s>2$  is
needed in order to insure that $\partial_{x}u\in L^{\infty}$, an
essential assumption for the proof. Molinet, Saut and Tzvetkov
\cite{MST1}
also proved that if one is willing to assume
more regularity for the initial data (at least three derivatives in the
$x$ variable and two in the $y$ variable need to be in $L^{2}$),
then global well-posedness holds.
 Recently we \cite{CKS} were able
to obtain
well-posedness for small data in a weighted Sobolev space with
essentially $H^{2}$ regularity, we will return to this result later.

We recall a few known facts associated with the KP equations.
If one defines the Fourier
transform for a function $f(x,y)$  as
$$\hat{f}(\xi,\mu)=\int_{\dbR^{2}}f(x,y)e^{-i(x\xi+y\mu)}dx dy,$$
then it is easy to see that the dispersive function
associated to this equation is
\begin{equation}
 \om(\x,\m)=\xi^{3} -\gamma \frac{\mu^{2}}{\xi}.
\label{disp}\end{equation}
The analysis of the KP initial value problem depends crucially on
the $sign$ of $\gamma$. We describe {\it{three differences}} due to the choice
of sign: the strength of the smoothing effect, the bilinear dispersive
identity and (non)positivity of the top order terms in the energy.

A first example of the relevance of the sign of $\gamma$
comes from  the
following observation. If we compute the gradient of $\om$, we have
that for KP-I ($\gamma=-1$, for example)
\begin{equation}
   |\nabla \om(\xi,\mu)|=|(3\x^{2}-\frac{\mu^{2}}{\x^{2}},2\frac{\m}{\x})|\geq
    |\xi|,
    \label{grad1}\end{equation}
and for KP-II ($\gamma=1$, for example)
\begin{equation}
   |\nabla \om(\xi,\mu)|=|(3\x^{2}+\frac{\mu^{2}}{\x^{2}},-2\frac{\m}{\x})|\geq
    |\xi|^{2}.
    \label{grad2}\end{equation}
Then, following the argument of Kenig, Ponce and Vega in \cite{KPV1},
we can claim that thanks to \eqref{grad2} KP-II recovers a full
derivative smoothness
along  the $x$ direction,  while by  \eqref{grad1}  KP-I recovers
only $\half$ derivative
smoothness along that same direction. Because the nonlinear term in
(\ref{ivp}) presents a derivative along the $x$ direction, this
 explains, at least
formally, why well-posedness questions for the KP-I IVP are much more
difficult to answer than for the KP-II problem.

The ``sign problem'' illustrated above appears also if one approaches
well-posedness questions using the method presented by Bourgain
in \cite{B}. This method
is based on the strength of various denominators which are controlled
using the bilinear dispersive identity
\begin{eqnarray}
 \nonumber & &\om(\xi_{1}+\x_{2},\m_{1}+\m_{2})-\om(\x_{1},\m_{1})
 -\om(\x_{2},\m_{2})\\
 &=&\frac{\x_{1}\x_{2}}{(\x_{1}+\x_{2})}
 \left(3(\x_{1}+\x_{2})^{2}+\gamma\left(\frac{\m_{1}}{\x_{1}}-
 \frac{\m_{2}}{\x_{2}}\right)^{2}\right).
\label{den}\end{eqnarray}
Clearly if $\gamma<0$ ( KP-I) this quantity could be zero, while if
$\gamma>0$ (KP-II)
$$|\om(\xi_{1}+\x_{2},\m_{1}+\m_{2})-\om(\x_{1},\m_{1})
 -\om(\x_{2},\m_{2})|\geq C|\x_{1}||\x_{2}||\x_{1}+\x_{2}|.$$
This is enough to  control the derivative in the nonlinear
term and to obtain well-posedness results for very rough data
(see also Takaoka \cite{Ta1}, Takaoka-Tzvetkov \cite{TT}).

The IVP (\ref{ivp}) has two conserved integrals, the $L^{2}$-norm and
the Hamiltonian:
\begin{eqnarray}
 \label{l2}\|u\|_{L^{2}}&=&\|u_{0}\|_{L^{2}}\\
  \label{hamiltonian} H(u)&=&\int_{\dbR^{2}}((\partial_{x}u)^{2}-\gamma
    (\partial_{x}^{-1}\partial_{y}u)^{2}-\frac{\beta}{3}u^{3})dx dy
    =H(u_{0}).
\end{eqnarray}
This time, for KP-I, the sign is favorable. In fact one can prove
\footnote{See for example \cite{C}.}
that a combination of  \eqref{l2} with  \eqref{hamiltonian} when
$\gamma=-1$,
gives
\begin{equation}
    \|u(t)\|_{L^{2}}+\|\partial_{x}u(t)\|_{L^{2}}+
    \|\partial_{x}^{-1}\partial_{y}u(t)\|_{L^{2}}\leq C
    \left(\|u_{0}\|_{L^{2}}+\|\partial_{x}u_{0}\|_{L^{2}}+
    \|\partial_{x}^{-1}\partial_{y}u_{0}\|_{L^{2}}\right),
\label{en}\end{equation}
for any sufficiently smooth solution $u$, uniformly in time.
The Sobolev space defined
by  \eqref{en} is naturally called the {\em energy space}.
It is the natural space on which the Hamiltonian is defined, and thus
it would be desirable to obtain a local well-posedness theory for KP-I
in this space. (As we mentioned before, Tom \cite{T} proved the
existence of global weak solutions, for data in the energy space,
using \eqref{en} and compactness arguments, but the uniqueness of
these weak solutions remains an open problem). Moreover, if one could
also prove that the time $T$ of existence in this (desired) local
existence theorem depends only on the norms involved in \eqref{en},
then a simple iteration argument, combined with \eqref{en}, would
yield global in time solutions for data in the energy space, and hence
the Hamiltonian would be defined globally in time, for the natural
space of initial data, providing a satisfactory ``low-regularity''
space in which KP-I is globally well-posed, and in which the
Hamiltonian is naturally defined. We next remark that this desired
dependence on $T$ above is validated by scaling considerations. In
fact, if we fix $|\beta| = 1, \gamma = -1$ and $u(x,y,t)$ is a
solution of \eqref{ivp}, then $u_\lambda (x, y,t) = \lam^2 u(\lam x ,
\lam^2 y , \lam^3 t)$ is also a solution of \eqref{ivp}, with initial
data $u_{\lam, 0} (x, y) = \lam^2 u_0 (\lam x , \lam^2 y).$ Note that
${{\| u_{\lam, 0} \|}_{L^2}} = \lam^\half {{\| u_0\|}_{L^2}}$, ${{\|
    \partial_x u_{\lam, 0} \|}_{L^2}} = \lam^{\frac{3}{2}} {{\|
    \partial_x u_0 \|}_{L^2}}, {{\| \partial_x^{-1} \partial_y u_{\lam,
      0} \|}_{L^2}} = \lam^{\frac{3}{2}} {{\| \partial_x^{-1}
    \partial_y u_{ 0} \|}_{L^2}}$, so that \eqref{ivp} is
``sub-critical'' in the energy space and thus one expects the time of
existence in a local well-posedness theorem, as the one discussed
before, to depend only on the norm of the initial data in the energy
space. (See also Remark \ref{rescaling} for further discussion of the
notion of ``criticality''). Note also that if one is only interested
in global existence of solutions of KP-I, with fairly regular initial
data, the recent work \cite{MST1} provides a very satisfactory global
existence theory, by combining the local well-posedness results of
Iorio and Nunes \cite{IN}, mentioned before, with higher order
conservation laws for KP-I (suitably regularized by the use of
Strichartz inequalities).

In the attempt to establish a local well-posedness theory for KP-I in
the energy space, one is confronted by the following difficulty, which
we have not been able to overcome: so far, in the many studies of
local well-posedness for nonlinear dispersive equations, the only
successful approach to the issue of ``low regularity'' data has been
through the use of fixed point theorems based on Picard
iteration. However, the recent counterexamples of Molinet, Saut and
Tzvetkov \cite{MST}, \cite{MST1}, show that, for KP-I, we cannot prove
local well-posedness in any type of anisotropic $L^2$-based Sobolev
space $H^{s_1, s_2}_{x,y},$ or in the energy space, by using Picard
fixed point methods for the integral equation formulation of the KP-I
initial value problem. In light of this, to study the local
well-posedness theory in the energy space, one must abandon Picard
iteration, and proceed in a new way. Since weak solutions have been
constructed in \cite{T}, as we mentioned before, the key issue is
uniqueness and one needs to establish this without relying on the
classical Gronwall inequality, which seems to require too much
regularity on the data. Possibly, recent works of Molinet and Ribaud
on dissipative generalizations of KP \cite{MRKP} and KdV \cite{MRKdV}
may prove useful in this direction, but we have not been able to
establish the required uniqueness.

Given this unsatisfactory state of affairs, an alternative is to use
spaces other than $H^{s_1, s_2}_{x,y}$, or the energy space, but with
similar regularity properties, and for which Picard iteration might
still work. For example, in our recent work \cite{CKS}, we addressed
the well-posedness question for KP-I, by restricting the space of
initial data, which we took to consist (essentially) of functions,
which together with two derivatives, belonged to the weighted
space $L^2((1+|y|)^\alpha dx dy), \alpha > \half$, and have small
norm. Our proof relied on the so called ``oscillatory integrals''
method, which combines local smoothing effects and maximal function estimates.

Our goal in the present paper is to refine the local well-posedness
result mentioned above, to reduce the number of the derivatives needed
on the initial data to bring it to a space which is close to the energy space we
discussed before. We use versions of the spaces and methods introduced
by Bourgain \cite{B}, extended to the context of weighted Besov
spaces. The weights are used to exploit the fact that the region where
\eqref{den} is small, is a region of small measure. 
The estimates
we present are sharp, in a sense that will be made clear later,
and are obtained in Besov-type spaces involving derivatives of order
$1-\epsilon$ and the weight $|y|$. We are able to remove any assumption
on the initial data concerning small frequency, but due to the fact
that in this case weighted spaces do not rescale well (see Remark
\ref{rescaling}),
our well-posedness result again holds only for small data. From now on
we will restrict ourselves to $|\beta| = 1$.

Let's now define the energy space $E$ and the weighted space $P$ as
\begin{equation}
\label{EPdefined}
E=\{f : f\in L^{2}, \partial_{x}f \in L^{2},
    \partial_{y}\partial_{x}^{-1}f \in L^{2}\}, \, \mbox{ and }
    P=\{f : yf \in L^{2}\}.
\end{equation}

% OLD REMARK 1.1 REPLACED BY BELOW ON 10 FEB 2003
%%     \begin{remark}\label{natep}
%% We consider the space $E\cap P$  natural in the context of KP-I.
%% It was proved by Saut \cite{S} that if the initial data $u_{0}$ is
%% in $E\cap P$, then the solution $u$ for \eqref{ivp} in a finite
%% interval $[-T,T]$ is such that $u\in L^{\infty}([-T,T],P)$.  (This
%% uses \eqref{en}). This result combined with \eqref{en},
%% shows that  in a fixed interval of time $[-T,T]$,
%% if the initial data $u_{0}$ is in $E\cap P$, then the solution $u$
%% enjoys the a priori bound
%% %
%% $$\|u\|_{L^{\infty}_{[-T,T]}(E\cap P)}\leq C(T, \|u_{0}\|_{E\cap P}).$$
%% \end{remark}
%% %
\begin{remark}
  \label{natep}
We consider the space $E \cap P$ natural in the context of KP-I. It
was proved by Saut \cite{S}, that for smooth solutions $u$ of KP-I,
whose initial data $u_0$ is in $E\cap P$, then for any fixed time
interval $[-T,T$, $u$ enjoys the a priori bound
\begin{equation*}
  {{\| u \|}_{L^\infty_{[-T, T]} (E \cap P)}} \leq C( T, {{\| u_0
  \|}_{E\cap P}} ).
\end{equation*}
\end{remark}

Let us   denote now  by $B_{\rho}$ the ball in $E\cap P$,
centered at zero, and
radius $\rho$. To state  the main theorem we will also
need the spaces $(E\cap P)_{1-\epsilon}$, and $Z_{1-\epsilon}$.
The first space will be defined in
\eqref{epsilon}, but for now, all the reader needs to know about it
is that it roughly has $\epsilon$ fewer derivatives than the
space $E\cap
P$. The space $Z_{1-\epsilon}$ is introduced in \eqref{zsb}. It is a
Bourgain type space (following the spaces introduced in \cite{B}), 
in which the contraction mapping theorem is
applied.
%THEOREM RESTATED BELOW 10 FEB 2003
%% \begin{theorem}\label{energytheorem}
%% Assume that $\gamma=-1$ in \eqref{ivp} and fix an interval of time
%% $[0,T]$.
%% Then there exists $\delta>0$  such
%% that for any $u_{0}\in E\cap P$, and $\|u_{0}\|_{E\cap P}\leq \delta$,
%% there exists a unique solution $u$ for  the IVP \eqref{ivp}, such
%% that $u\in L^{\infty}([0,T], E\cap P)\cap C([0,T], (E\cap
%% P)_{1-\epsilon})\cap Z_{1-\epsilon}, \, \, \|u\|_{Z_{1-\epsilon}}\leq
%% C\delta$, for any small positive $\epsilon$. Moreover, the map
%% that associates the initial data in $E\cap P$ to the solution $u$ is
%% smooth from the ball $B_{\delta}$  into the space $C([0,T], (E\cap
%% P)_{1-\epsilon})$, for any positive and small $\epsilon$.
%%     \end{theorem}

\begin{theorem}
  \label{energytheorem}
Assume that $\gamma = -1$ in \eqref{ivp} and fix an interval of time
$[0,T]$, and a small $\epsilon > 0$. Then, there exists $\delta > 0 $,
$\delta = \delta (\epsilon, T)$ such that for any $u_0 \in E\cap P$
with ${{\| u_0 \|}_{E\cap P}} < \delta$, there exists a unique
solution 
\begin{equation*}
  u \in L^\infty ([0,T]; E\cap P) \cap C([0,T]; (E\cap
  P)_{1-\epsilon}) \cap Z_{1-\epsilon}
\end{equation*}
with ${{\| u \|}_{Z_{1-\epsilon}}} \leq C \delta$. Moreover, the map
that associates the initial data in $E \cap P$ to the solution $u$ is
smooth from the ball $B_\delta$ into the space $C([0,T]; (E\cap P)_{1-\epsilon})$.
\end{theorem}

This theorem  is a consequence of a well-posedness result
involving the Besov type spaces of initial data mentioned earlier,
(see Theorem \ref{main0} below). We start
by giving  a precise definition for  these  spaces.
\begin{definition}
Let $\theta_{0}(s)=\chi_{[-1,1]}(s), \, \, \theta_{m}(s)=
\chi_{[2^{m-1},2^{m}]}(|s|), \, \, m\in \dbN$. For $(\x,\m)\in
 \dbR^{2}$
let $\chi_{1}(\x,\m)=\chi_{\{|\x|\geq \half \frac{|\mu|}{|\x|}\}}$, and
$\chi_{2}(\x,\m)=\chi_{\{|\x|< \half{\frac{|\mu|}{|\xi|}}\}}$. We define
the space $B^{2,1}_{s}$  of functions on $\dbR^{2}$ as the closure of
the Schwartz functions, for which the norm below is finite, with respect to
\begin{eqnarray}
\label{bs}    \|f\|_{B^{2,1}_{s}}&=&\sum_{m\geq 0}\|(1+|\x|+\gi)^{s}
    \theta_{m}(\x)\chi_{1}(\x,\m)\hat{f}(\x,\m)\|_{L^{2}}\\
\nonumber    &+&\sum_{n\geq 0}\|(1+|\x|+\gi)^{s}
    \theta_{n}(\m)\chi_{2}(\x,\m)\hat{f}(\x,\m)\|_{L^{2}}.
    \end{eqnarray}
We also define a ``weighted Besov space'', $P^{2,1}_{r}$ using the
norm
\begin{eqnarray}
\label{pr}    \|f\|_{P^{2,1}_{r}}&=&\sum_{m\geq 0}\|(1+|\x|+\gi)^{r}
    \theta_{m}(\x)\chi_{1}(\x,\m)\partial_{\m}\hat{f}(\x,\m)\|_{L^{2}}\\
\nonumber    &+&\sum_{n\geq 0}\|(1+|\x|+\gi)^{r}
    \theta_{n}(\m)\chi_{2}(\x,\m)\partial_{\m}\hat{f}(\x,\m)\|_{L^{2}}.
    \end{eqnarray}
    \end{definition}
\begin{remark}\label{diffregions}
Going back to the discussion of the smoothing effect involving
\eqref{grad1} and \eqref{grad2}, one can see that
the splitting into the two regions
$R_{g}=\{|\x|\geq \half \frac{|\mu|}{|\x|}\}$ and
$R_{b}=\{|\x|<\half \frac{|\mu|}{|\x|}\}$ is quite natural.
In fact in  the ``good'' region $R_{g}$ it's easy to check that
$|\nabla \om(\xi,\mu)|\gtrsim |\x|^{2}$, hence here one
should expect a gain of a full derivative. On the other hand
in  the ``bad'' region $R_{b}$ one has
$|\nabla \om(\xi,\mu)|\gtrsim |\x|$, and the gain
should be only of half  derivative, (see also Proposition \ref{regions}).
    \end{remark}
\begin{remark}\label{embedding}
    Because $l^{1}\subset l^{2}$, it follows that $B^{2,1}_{1}\cap
P^{2,1}_{0} \subseteq E\cap P$. Moreover, if $\epsilon>0$, then we
also have  $E\cap P\subseteq
B^{2,1}_{1-\epsilon}\cap P^{2,1}_{-\epsilon}$. To see  this, first
assume that $f\in E$. Then
\begin{eqnarray*}
&&    \sum_{m\geq 0}\|(1+|\x|+\gi)^{1-\epsilon}
    \theta_{m}(\x)\chi_{1}(\x,\m)\hat{f}(\x,\m)\|_{L^{2}}\\
    &\lesssim&
    \sum_{m\geq 0}\|2^{-m\epsilon}(1+|\x|)^{\epsilon}(1+|\x|+
\gi)^{1-\epsilon}
    \theta_{m}(\xi)\chi_{1}\hat{f}\|_{L^{2}}\\
    &\lesssim&\sum_{m\geq 0}2^{-m\epsilon}\|(1+|\x|+\gi)
    \theta_{m}(\xi)\chi_{1}\hat{f}\|_{L^{2}},
\end{eqnarray*}
and Cauchy Schwarz concludes this part. On the other hand
\begin{eqnarray*}
&&    \sum_{n\geq 0}\|(1+|\x|+\gi)^{1-\epsilon}
    \theta_{n}(\m)\chi_{2}(\x,\m)\hat{f}(\x,\m)\|_{L^{2}}\\
    &\lesssim& \sum_{n\geq
    0}\left\|\frac{(|\x|\gi)^{\half \epsilon}}
    {(|\x|\gi)^{\half \epsilon}}(1+|\x|+\gi)^{1-\epsilon}
    \theta_{n}(\m)\chi_{2}(\x,\m)\hat{f}(\x,\m)\right\|_{L^{2}}\\
 &\lesssim& \sum_{n\geq
    0}2^{-\frac{n}{2}\epsilon}\|
    (|\x|^{\epsilon}+(\gi)^{\epsilon})(1+|\x|+\gi)^{1-\epsilon}
    \theta_{n}(\m)\chi_{2}(\x,\m)\hat{f}(\x,\m)\|_{L^{2}},
\end{eqnarray*}
and again Cauchy-Schwarz takes care of this term. Now assume that $f\in
P$. Then
\begin{eqnarray*}
&& \sum_{m\geq 0}\|(1+|\x|+\gi)^{-\epsilon}
    \theta_{m}(\x)\chi_{1}(\x,\m)\partial_{\m}\hat{f}(\x,\m)\|_{L^{2}}\\
    &\lesssim& \sum_{m\geq 0}2^{-m\epsilon}\|
  \theta_{m}(\x)\chi_{1}(\x,\m)\partial_{\m}\hat{f}(\x,\m)\|_{L^{2}},
  \end{eqnarray*}
and we use Cauchy-Schwarz. Finally, because
$(1+|\x|+\gi)\geq 1+(|\x|\gi)^{\half}$, we obtain
\begin{eqnarray*}
&&\sum_{n\geq 0}\|(1+|\x|+\gi)^{-\epsilon}
\theta_{n}(\m)\chi_{2}(\x,\m)\partial_{\m}\hat{f}(\x,\m)\|_{L^{2}}\\
 &\lesssim& \sum_{n\geq 0}2^{-{\frac{n}{2}}\epsilon}\|
  \theta_{n}(\m)\chi_{2}(\x,\m)\partial_{\m}\hat{f}(\x,\m)\|_{L^{2}},
    \end{eqnarray*}
and also this term is estimated. We are now ready to define the space
$(E\cap P)_{1-\epsilon}$, introduced in the statement of Theorem
\ref{energytheorem}, by setting
\begin{equation}\label{epsilon}
    (E\cap P)_{1-\epsilon}=B^{2,1}_{1-\epsilon}\cap P^{2,1}_{-\epsilon},
    \end{equation}
for any $\epsilon \in \dbR$.
\end{remark}
\begin{remark}
If one could prove well-posedness with
initial data
in $B^{2,1}_{1-\epsilon}\cap P_{-\epsilon}^{2,1}$, on $[-T,T]$, with
$T$ depending only on the norm of the initial data in this space, for
some $\epsilon>0$, then for data $u_{0}\in E\cap P$ we would obtain,
in light of Remarks \ref{embedding} and \ref{natep}, a unique global
solution in $C([-T,T], B^{2,1}_{1-\epsilon}\cap P_{-\epsilon}^{2,1})\cap
L^{\infty}([-T,T], E\cap P)$, for each $T$, which would depend
continuously on the initial data, in the
$B^{2,1}_{1-\epsilon}\cap P_{-\epsilon}^{2,1}$ topology. However,
as we will explain in  Remark \ref{rescaling}, we show the required
local well-posedness only for small data in
$B^{2,1}_{1-\epsilon}\cap P_{-\epsilon}^{2,1}$, and our estimates
barely miss giving the global result\footnote{Unfortunately though, our
estimates are sharp as is shown in Proposition \ref{caunter}.}.
\end{remark}
\begin{remark}\label{cver12}
    If in the definition of $B^{2,1}_{s}$, the constant $\half$ appearing
    in $\chi_{1}$ and $\chi_{2}$ is
    replaced by $C$, we obtain the same space, with comparable norms.
    This holds also for $P^{2,1}_{r}$.
Assume that $f\in B^{2,1}_{s}$ and $0<C<\half$. We need to show that
\begin{equation}
    \label{C12}\sum_{m\geq 0}\|(1+|\x|+\gi)^{s}
    \theta_{m}(\x)\chi_{\{C\gi<|\x|<\half\gi\}}(\x,\m)
    \hat{f}(\x,\m)\|_{L^{2}}\lesssim \|f\|_{B^{2,1}_{s}}.
   \end{equation}
But if $C\gi<|\x|<\half \gi$, then $C|\m|<|\x|^{2}<\half |\m|$. If $m=0$,
then $|\x|^{2}\leq 1$ and $|\m|\leq C^{-1}$. Then
\begin{eqnarray*}
&&\|(1+|\x|+\gi)^{s}
    \theta_{0}(\x)\chi_{\{C\gi<|\x|<\half \gi\}}(\x,\m)
    \hat{f}(\x,\m)\|_{L^{2}}\\
    &\lesssim&\|(1+|\x|+\gi)^{s}
    \chi_{\{|\m|\leq C^{-1}\}}\chi_{\{|\x|<\half \gi\}}(\x,\m)
    \hat{f}(\x,\m)\|_{L^{2}}\\
    &\lesssim&C\sum_{n\geq 0}\|(1+|\x|+\gi)^{s}
    \theta_{n}(\m)\chi_{2}(\x,\m)\hat{f}(\x,\m)\|_{L^{2}}.
    \end{eqnarray*}
If $m\geq 1,  \, 2^{2m-1}\leq |\m|\leq 2^{2m}/C$. Let $n_{0}$ be the
smallest integer such that $2^{n_{0}}\geq 2^{2m}/C$. Then $n_{0}\leq
2m+C_{0}, C_{0}=C_{0}(C)$. Thus,
\begin{eqnarray*}
&&\sum_{m\geq 1}\|(1+|\x|+\gi)^{s}
    \theta_{m}(\x)\chi_{\{C\gi<|\x|<\half \gi\}}(\x,\m)
    \hat{f}(\x,\m)\|_{L^{2}}\\
    &\lesssim&\sum_{m\geq 1}\sum_{2m-1\leq n_{0}\leq 2m+C_{0}}
    \|(1+|\x|+\gi)^{s}
    \theta_{n_{0}}(\m) \chi_{2}(\x,\m)
    \hat{f}(\x,\m)\|_{L^{2}}\\
    &\lesssim&\sum_{n_{0}\geq 1}\|(1+|\x|+\gi)^{s}
    \theta_{n_{0}}(\m)\chi_{2}(\x,\m)\hat{f}(\x,\m)\|_{L^{2}},
    \end{eqnarray*}
    and \eqref{C12} follows. We also need to show that
\begin{equation}
    \label{C121}\sum_{n\geq 0}\|(1+|\x|+\gi)^{s}
    \theta_{n}(\m)\chi_{\{C\gi<|\x|<\half \gi\}}(\x,\m)
    \hat{f}(\x,\m)\|_{L^{2}}\lesssim \|f\|_{B^{2,1}_{s}}
   \end{equation}
and the argument is similar since
$$2^{n}\leq |\m|\leq 2^{n+1} \Longrightarrow C^{\half }2^{\frac{n}{2}}\leq
|\x|\leq (\half )^{\half}2^{n/2+1}.$$
The case $C>1/2$ is proved in the same way, reversing the role of $C$
and $1/2$. A similar proof can be given for the space $P_{r}^{2,1}$.
This remark will be used implicitly in our proofs.
\end{remark}
We are now ready to introduce the Banach spaces in which we
 will perform
a fixed point argument to obtain the solution for \eqref{ivp}. Below,
 we use $\widehat{f}$ to denote the Fourier transform of a function of
 $(x,y,t)$, defined in a similar fashion as for functions of
 $(x,y)$. We hope that this will not cause confusion to the reader.
\begin{definition}
Let $\chi_{0}(s)=\chi_{\{|s|<1\}}(s), \chi_{j}(s)=
\chi_{\{2^{j-1}\leq |s|<2^{j}\}}(s), \om(\x,\m)=\x^{3}+\mu^{2}/\x$ and
$w(\x,\m)=(1+|\x|+\gi)$. We
define the space $X_{s,b}$ through the following norm: 
%{\bf{[GAFA    HELP?: make label not touch the display?]}}
%
\begin{eqnarray}
    \label{xs}\|f\|_{X_{s,b}}&=&\sum_{j, m\geq 0}
    2^{jb}\left(\int_{\dbR^{3}}
    \chi_{j}(\tau-\om(\xi,\m))\chi_{1}(\xi,\m)\theta_{m}(\xi)
    w^{2s}|\hat{f}|^{2}(\x,\m,\tau)d\x d\m d\tau\right)^{\half}\\
    \nonumber&+&\sum_{j, m \geq 0}
    2^{jb}\left(\int_{\dbR^{3}}
    \chi_{j}(\tau-\om(\xi,\m))\chi_{2}(\xi,\m)\theta_{n}(\m)
    w^{2s}|\hat{f}|^{2}(\x,\m,\tau)d\x d\m d\tau\right)^{\half}.
    \end{eqnarray}
We also define the space
\begin{equation}
    \label{yrs}Y_{s,r,b}=\{f : tf \in X_{s,b}, \mbox{ and } yf \in
    X_{r,b}\},
    \end{equation}
    and the spaces
\begin{equation}
    \label{zsb}Z_{s,b}=X_{s,b}\cap Y_{s,s-1,b}, \, \,
    Z_{1-\epsilon}=Z_{1-\epsilon,\half}.
    \end{equation}
    \label{spaces}\end{definition}
\begin{remark} A statement similar to Remark \ref{cver12} holds for
these spaces.
\end{remark}
We are now ready to state the well-posedness result for initial data
in Besov spaces introduced above.

\begin{theorem}
Assume that $\gamma=-1$ in \eqref{ivp}. For any $\epsilon_{0}<\frac{1}{16}$ and
for any interval of time $[0,T]$,
there exists $\delta=\delta(\epsilon_{0},T)>0$
such that  for any $u_{0}\in B_{1-\epsilon_{0}}^{2,1}\cap
P_{-\epsilon_{0}}^{2,1}$ and
$\|u_{0}\|_{B_{1-\epsilon_{0}}^{2,1}\cap P_{-\epsilon_{0}}^{2,1}}\leq
\delta$, there exists a
unique solution $u$ for
\eqref{ivp} in $Z_{1-\epsilon_0}$. Moreover $u \in
B^{2,1}_{1-\epsilon_0} \cap P^{2,1}_{-\epsilon_0}$ 
and smoothness with respect to initial data holds in the appropriate
topology.
\label{main0}\end{theorem}
From now on we assume that $\gamma=-1$.
In the rest of the paper we  often use the notation $A \lesssim B$ if there
 exists $C>0$ such that $A\leq CB$, and $A\thicksim B$ if $A \lesssim
 B$ and $B \lesssim A$ with possibly different $C's.$

The paper is organized in the following way. In Section 2, we
introduce
some estimates for the solution of the linear KP-I initial value
problem. In Section 3, we present two bilinear
estimates (see Theorems \ref{xbil} and \ref{ybil} below)
that are  the heart of the matter for the proof of Theorem
\ref{main0}. The section concludes with a counterexample showing the
optimality of our analysis.
We finish with Section 4, in which we briefly present the proofs of Theorems
\ref{energytheorem} and \ref{main0}. Section 4 also contains a
scaling argument which reveals that the optimal analysis in Section 3
is ``endpoint critical''.

\vskip 7pt
\noindent {\bf{Acknowledgment.}} The authors are very grateful to the
referees for their extremely careful reading of the manuscript, and
their many suggestions, that have greatly clarified our original
version of the paper.

%%%%%%%%%%%%%%%%%%%%%%%%%%%%%%%%%%%%%%%%%%%%%%%%%%
\section{The Linear Estimates}

Consider the linear IVP
\begin{equation}
\left\{ \begin{array}{l}
\partial_t u + \partial_x^3 u -\partial_x^{-1}\partial_y^2u= 0, \\
 u(x,0) = u_0(x) \hspace{1.5cm}(x,y) \in \dbR^2, \,  t \in \dbR,
\end{array}\right.
\label{livp}\end{equation}
and let $S(t)u_{0}$ be the solution. By taking the Fourier transform of
the first equation in \eqref{livp} and solving the ODE one can easily
see that
$$S(t)u_{0}(x,y)=\int_{\dbR^{2}}e^{i(t\om(\x,\m)+x\x+y\m)}
\widehat{u_{0}}(\x,\m) d\x d\m.$$

We now show that the space $X_{s,\half}\cap Y_{s,s-1,\half}$ is well behaved
with respect to the group operator $S(t)$.
\begin{proposition} ({{linear homogeneous estimates}})
 Assume $\psi \in C^{\infty}_{0}(|t|\leq 1), \psi =1$ on $ |t|<\half$. Then
\begin{eqnarray}
    \label{xsgroup}\|\psi(t) S(t)u_{0}\|_{X_{s,\half}}
    &\lesssim&
    \|u_{0}\|_{B_{s}^{2,1}},\\
\label{ysgroup}\|\psi(t) S(t)u_{0}\|_{Y_{s,s-1,\half}}&
\lesssim&
\|u_{0}\|_{P^{2,1}_{s-1}}+\|u_{0}\|_{B^{2,1}_{s}}.
\end{eqnarray}
\label{linear}\end{proposition}
\begin{proof}
The proof follows the same arguments used in  \cite{KPV3}. We observe
that
\begin{equation}(\psi(t) S(t)u_{0})\sphat (\x,\m,\tau)=\hat{\psi}
(\tau-\om(\x,\m))\widehat{u_{0}}(\x,\m).
\label{fpsi}\end{equation}
Here $\sphat$ denotes the Fourier transform of a function of 3
variables on the left side of \eqref{fpsi} and also to denote the
transform of functions of 1 and 2 variables on the right side. 
Then to prove \eqref{xsgroup} we need to estimate the two integral
expressions:
\begin{eqnarray}
\label{xsgroup1}& &\sum_{j\geq 0}2^{\frac{j}{2}}\sum_{m\geq 0}
\left(\int_{\dbR^{3}}w(\x,\m)^{2s}\chi_{1}\theta_{m}(\x)
\chi_{j}(\tau-\om)|\hat{\psi}
(\tau-\om)|^{2}|\widehat{u_{0}}|^{2}d\x d\m
d\tau\right)^{\half}\\
\label{xsgroup2}& &\sum_{j\geq 0}2^{\frac{j}{2}}\sum_{n\geq 0}
\left(\int_{\dbR^{3}}w(\x,\m)^{2s}\chi_{2}\theta_{n}(\m)
\chi_{j}(\tau-\om)|\hat{\psi}
(\tau-\om)|^{2}|\widehat{u_{0}}|^{2}d\x d\m
d\tau\right)^{\half},
\end{eqnarray}
where $w(\x,\m)=(1+|\x|+\gi)$. We observe that for $j=0$
\begin{equation}
\int_{\dbR}|\hat{\psi}(\lambda)|^{2}\chi_{j}(\lambda)
d\lambda\lesssim \|\hat{\psi}\|_{L^{\infty}}^{2}
\label{lambdaj0}\end{equation}
    and for $j\geq 1$
\begin{equation}
\int_{\dbR}|\hat{\psi}(\lambda)|^{2}\chi_{j}(\lambda)
d\lambda\lesssim 2^{j}\frac{1}{(1+2^{j})^{2N}}
\|(1+|s|)^{N}\hat{\psi}(s)\|_{L^{\infty}}^{2}
\label{lambdaj1}\end{equation}
for any $N\in \dbN$. When we insert \eqref{lambdaj0} and
\eqref{lambdaj1} in \eqref{xsgroup1} we obtain the bound
\begin{equation}
  \|u_{0}\|_{B^{2,1}_{s}}\|\hat{\psi}\|_{L^{\infty}}+\sum_{j\geq 1}
  \frac{
2^{j}}{(1+ 2^{j})^{N}}\|u_{0}\|_{B^{2,1}_{s}}\|
(1+|s|)^{N}\hat{\psi}(s)\|_{L^{\infty}}.
\label{pro21}\end{equation}
It is easy to see that for $N>1$, $\sum_{j\geq 1}\frac{
2^{j}}{(1+2^{j})^{N}}\leq C$, hence \eqref{xsgroup} is proved
for \eqref{xsgroup1}. A similar argument can be used to estimate
\eqref{xsgroup2}.

To estimate \eqref{ysgroup} we first observe that
$$t\psi(t) S(t)u_{0}=\widetilde{\psi}
(t) S(t)u_{0},$$
where $\widetilde{\psi}(t)=t\psi(t)$. Hence by \eqref{xsgroup}
$$\|t\psi(t)
S(t)u_{0}\|_{X_{s,b}}\lesssim \|u_{0}\|_{B^{2,1}_{s}}.$$
We then turn to $y\psi(t)
S(t)u_{0}$. Using the fact that
$(yh(y))\sphat=-i\partial_{\m}\hat{h}(\m)$ and \eqref{fpsi},
it is easy to see that
\begin{equation}\F(y\psi(t)
S(t)u_{0})(\x,\m,\tau)=-2\widehat{\widetilde\psi}(\tau-\om)
\frac{\m}{\x} \widehat{u_{0}}(\x,\m)-
\hat{\psi}(\tau-\om)\widehat{yu_{0}}(\x,\m).
\label{ytranformation}\end{equation}
Then we can use  \eqref{xsgroup} to conclude that
$$\|y\psi(t)
S(t)u_{0}\|_{X_{s-1,b}}\lesssim \|u_{0}\|_{B^{2,1}_{s}}+
\|yu_{0}\|_{B^{2,1}_{s-1}}.$$
\end{proof}

There is  an inhomogeneous version of Proposition \ref{linear}.
\begin{proposition} (linear inhomogeneous estimates)
Assume $\psi \in C^{\infty}_{0}(|t|\leq 1), \psi =1$ on $ |t|<\half$. Then,
\begin{eqnarray}
    \label{xsnonhom}& &\|\psi(t) \int_{0}^{t}S(t-t')h(t')\,dt'
    \|_{X_{s,\half}}
    \lesssim
    \|h\|_{X_{s,-\half}},\\
\label{ysnonhom}& &\|\psi(t)\int_{0}^{t}S(t-t')h(t')
dt'\|_{Y_{s,s-1,\half}}\lesssim (\|h\|_{X_{s,-\half}}+
\|th\|_{X_{s,-\half}}+\|yh\|_{X_{s-1,-\half}}).
\end{eqnarray}
\label{linearhom}\end{proposition}
Also in this case the proof follows closely the arguments used in
\cite{KPV3}. We start with the following lemma.
\begin{lemma} (stability under time cutoff)
    There exists $C>0$ such that
    for any $s\in \dbR$,
$$
\|\psi(t)f\|_{X_{s,\half}}\leq C
\|f\|_{X_{s,\half}}.
$$
\label{xsepsilon}\end{lemma}
To prove the lemma we need the auxiliary space $\tilde{X}_{s,b}$
defined as the closure of the functions in $S$ 
for which the norm below is finite, with respect to the norm
$$
  \|f\|_{\tilde{X}_{s,b}}=
  \left(\sum_{j\geq 0}
    2^{j2b}\int_{\dbR^{3}}
    \chi_{j}(\tau-\om(\x,\m))
   (1+|\x|+\gi)^{2s}|\hat{f}|^{2}(\x,\m,\tau)d\x d\m
    d\tau\right)^{\half}.$$
Notice that while the space $X_{s,b}$ is defined using an $l^{1}$
summation with respect to $j$, the space $\tilde{X}_{s,b}$ is defined
using an $l^{2}$ summation. We have the following lemma.
\begin{lemma}
For any $b$ such that $0<b<\half$ and any $s\in \dbR$ we have
\begin{equation}
    \label{tilex1}\|\psi(t)f\|_{\tilde{X}_{s,b}}\lesssim
    \|f\|_{\tilde{X}_{s,b}}.
    \end{equation}
For any $b$ such that $\half<b<1$ and any $s\in \dbR$ we have
\begin{equation}
    \label{tilex2}\|\psi(t)f\|_{\tilde{X}_{s,b}}\lesssim
    \|f\|_{\tilde{X}_{s,b}}.
    \end{equation}
\label{tilde}\end{lemma}
\begin{proof}
We start by proving  \eqref{tilex1}. Note that
\begin{eqnarray*}& &\left(\sum_{j\geq 0}
    2^{j2b}\int_{\dbR^{3}}
    \chi_{j}(\tau-\om(\x,\m))w^{2s}|\hat{f}|^{2}(\x,\m,\tau)d\x d\m
    d\tau\right)^{\half}\\
    &\simeq &\left(
    \int_{\dbR^{3}}(1+|\tau-\om(\x,\m)|)^{2b}
    w^{2s}|\widehat{f}|^{2}(\x,\m,\tau)d\x d\m
    d\tau\right)^{\half}
\end{eqnarray*}
where $w(\x,\m)=(1+|\x|+\gi)$.
Because
$$(\psi(\cdot)f)\sphat (\x,\m,\tau)=\hat{f}*_{\tau}
(\hat\psi(\cdot)),$$
by following the arguments in \cite{KPV3}, it is  easy to see
that the proof reduces to showing that for any $a\in \dbR$, (with
$(D^b f)\sphat (\tau) = |\tau|^b \widehat{f} (\tau)$)
\begin{equation}
\int_{\dbR}|D^{b}(e^{iat}f(t)\psi(t))|^{2}dt\lesssim
\int_{\dbR}|\hat{f}(l)|^{2}(1+|l-a|)^{2b}dl.
\label{reduction}    \end{equation}
We use fractional derivatives (see appendix in
\cite{KPV2}) to obtain
$$\|D^{b}(e^{iat}f(t)\psi(t))-
D^{b}(e^{iat}f(t))\psi(t)-
e^{iat}f(t)D^{b}(\psi(t))\|_{L^{2}}\leq
\|D^{b}(e^{iat}f(t))\|_{L^{2}}\|\psi\|_{L^{\infty}}.$$
It follows that
$$\|D^{b}(e^{iat}f(t)\psi(t))\|_{L^{2}}\leq
\|D^{b}(e^{iat}f(t))\|_{L^{2}}\|\psi\|_{L^{\infty}}+
\|e^{iat}f(t)\|_{L^{2r}}
\|D^{b}(\psi(t))\|_{L^{2r'}},$$
for $\frac{1}{r}+\frac{1}{r'}=1$.
Because $b<\half$, if $\frac{1}{r} = 1-2b$, by the Sobolev
embedding theorem,
we can continue with
$$\|D^{b}(e^{iat}f(t)\psi(t))\|_{L^{2}}
\leq \left(\int_{\dbR}
|\hat{f}(l)|^{2}(1+|l-a|)^{2b}dl\right)^{\half}
(\|\psi\|_{L^{\infty}}+\|D^{b}(\psi(t))\|_{L^{2r'}}).$$
To finish observe that
$$\|D^{b}(\psi(t))\|_{L^{2r'}}<\infty,$$
 because $D^{b}(\psi)\in L^{2}\cap L^{\infty}$.

The proof of \eqref{tilex2} follows by combining  the above arguments
with those in Lemma 3.2 of \cite{KPV3}.
\end{proof}

Note that \eqref{tilex1} also follows from Lemma 2.2 in
\cite{GTV:ZS}. 

\noindent
\begin{proof}[Proof of Lemma \ref{xsepsilon}]
We recall that by real interpolation (see Theorem 5.6.1 in \cite{BL}),
if $A$ is a Banach space and
$$l^{q}_{b}(A)=\{(f_{j}) : f_{j}\in A, \, \,
\left(\sum_{j\geq 0}(2^{jb}\|f_{j}\|_{A})^{q}\right)^{\frac{1}{q}}<\infty\},$$
then for any $q_{i}\in [1,\infty]$,
$$(l^{q_{0}}_{b_{0}}(A), l^{q_{1}}_{b_{1}}(A))_{\theta,q}=
l^{q}_{b_{\theta}}(A), (b_0 \neq b_1)$$
where $\theta \in [0,1]$, $1\leq q\leq \infty$ and $b_{\theta}=\theta
b_{0}+(1-\theta)b_{1}$. We then apply this fact  to the
spaces $\tilde{X}_{s,b_{i}}$ with  $b_{0}<\half, \, b_{1}>\half, \, q=1$,
we use Lemma \ref{tilde} and then we sum with respect to $m$ and $n$
to obtain Lemma \ref{xsepsilon}.
\end{proof}
\begin{proof}[Proof of Proposition \ref{linearhom}]
We follow the proof of Lemma 3.3 in \cite{KPV3}. We write
$$\psi(t)\int_{0}^{t}S(t-t')h(t')dt'=I+II,$$
where
\begin{eqnarray*}
I&=&
\psi(t)\int_{-\infty}^{\infty}\int_{\dbR^{2}}
e^{i(x\x+y\m)}\hat{h}(\x,\m,\tau)\psi(\tau-\om)
\frac{e^{it\tau}-e^{it\om}}{\tau-\om(\x,\m)}d\x d\m d\tau\\
II&=&\psi(t)\int_{-\infty}^{\infty}\int_{\dbR^{2}}
e^{i(x\x+y\m)}\hat{h}(\x,\m,\tau)[1-\psi(\tau-\om)]
\frac{e^{it\tau}-e^{it\om}}{\tau-\om(\x,\m)}d\x d\m d\tau.
    \end{eqnarray*}
By Taylor expansion we can rewrite $I$ as
\begin{equation}\label{I}
    I=\sum_{k=1}^{\infty}\frac{i^{k}}{k!}t^{k}\psi(t)
\int_{\dbR^{2}}e^{i(x\x+y\m+t\om)}\left(\int_{-\infty}^{\infty}
\hat{h}(\x,\m,\tau)(\tau-\om(\x,\m))^{k-1}\psi(\tau-\om)d\tau\right)
d\x d\m.
\end{equation}
For $k\geq 1$ let
$$t^{k}\psi(t)=\psi_{k}(t),$$
and note that for any $k\geq 1$ and for any $s\in \dbR$,
$$|\widehat{\psi_{k}}(s)|\leq C,$$
and for $|s|>1$,
\begin{equation}
|\widehat{\psi_{k}}(s)|\leq C\frac{(1+k)^{2}}{(1+|s|)^{2}}.
\label{psik}\end{equation}
From \eqref{I} it is easy to see that
$$I=\sum_{k=1}^{\infty}\frac{i^{k}}{k!}
\psi_{k}(t)S(t)h_{k}(x,y),$$
where
$$\widehat{h_{k}}(\x,\m)=\int_{-\infty}^{\infty}
\hat{h}(\x,\m,\tau)(\tau-\om(\x,\m))^{k-1}\psi(\tau-\om)d\tau.$$
Then by Proposition \ref{linear},
in particular \eqref{pro21}, and \eqref{psik}, we obtain
$$\|I\|_{X_{s,\half}}\lesssim \sum_{k\geq 1}\frac{(1+k)^{2}}{k!}
\|h_{k}\|_{B^{2,1}_{s}}.$$
On the other hand it is easy to check that
$$\|h_{k}\|_{B^{2,1}_{s}}\lesssim \|h\|_{X_{s,-\half}},$$
which inserted above gives \eqref{xsnonhom}. We now pass to $II$.
We  write $II=II_{1}+II_{2}$, where
\begin{eqnarray*}
    II_{1}&=&\psi(t)\int_{-\infty}^{\infty}\int_{\dbR^{2}}
e^{i(x\x+y\m)}\hat{h}(\x,\m,\tau)[1-\psi(\tau-\om)]
\frac{e^{it\tau}}{\tau-\om(\x,\m)}d\x d\m d\tau,\\
II_{2}&=&-\psi(t)
\int_{\dbR^{2}}
e^{i(x\x+y\m)}\int_{-\infty}^{\infty}\hat{h}(\x,\m,\tau)
[1-\psi(\tau-\om)]
\frac{e^{it\om}}{\tau-\om(\x,\m)}d\tau d\x d\m.
\end{eqnarray*}
By Lemma \ref{xsepsilon}
$$\|II_{1}\|_{X_{s,\half}}\lesssim
\|h\|_{X_{s,-\half}}.$$
On the other hand, by Proposition \ref{linear} we have
$$\|II_{2}\|_{X_{s,\half}}\lesssim \|\tilde{h}\|_{B^{2,1}_{s}},$$
where
$$\widehat{\tilde{h}} (\xi, \mu) =\int_{-\infty}^{\infty}
[1-\psi(\tau-\om)]
\frac{\hat{h}(\x,\m,\tau)}{\tau-\om(\x,\m)}d\tau.$$
To finish the proof of \eqref{xsnonhom} one just needs to
observe that
$$\|\tilde{h}\|_{B^{2,1}_{s}}\lesssim \|h\|_{X_{s,-\half}}.$$
To prove \eqref{ysnonhom} we first observe that by
\eqref{xsnonhom}
$$\|t\psi(t)\int_{0}^{t}S(t-t')h(t')dt'
\|_{X_{s,\half}}\lesssim
\|h\|_{X_{s,-\half}}.$$
We now use \eqref{ytranformation} to write
\begin{eqnarray*}
& &\|y\psi(t)\int_{0}^{t}S(t-t')h(t')dt'
\|_{X_{s-1,\half}}\lesssim
\|\tilde{\psi}(t)\int_{0}^{t}S(t-t')
\bar{h}(t')dt'\|_{X_{s-1,\half}}\\
&+&\|\psi(t)\int_{0}^{t}S(t-t')t'\bar{h}(t')dt'
\|_{X_{s-1,\half}}+\|\psi(t)\int_{0}^{t}S(t-t')yh(t')dt'
\|_{X_{s-1,\half}},
\end{eqnarray*}
where $\tilde{\psi}(t)=t\psi(t)$ and $
\hat{\bar{h}}(\x,\m)=2\mu/\x\hat{h}(\x,\m)$. We use again
\eqref{xsnonhom} and we continue with
$$\|y\psi(t)\int_{0}^{t}S(t-t')h(t')dt'
\|_{X_{s-1,\half}}\lesssim\|h\|_{X_{s,-\half}}
+\|th\|_{X_{s,-\half}}
+\|yh\|_{X_{s-1,-\half}}.$$
This concludes the proof of Proposition \ref{linearhom}.

    \end{proof}

In this second part of the section we prove some a priori estimates
enjoyed by the solution $S(t)u_{0}$ of the linear problem \eqref{livp}.
The first estimate we present is
of Strichartz type and is due to Ben-Artzi and Saut \cite{BS}:
\begin{proposition} (linear homogeneous estimate)
Assume $u_{0} \in L^{2}$, then
\begin{equation}
    \|S(t)u_{0}(x,y)\|_{L^{4}([0,T],L^{4}(\dbR^{2}))}\lesssim
\|u_{0}\|_{L^{2}}.
\label{strichartz}\end{equation}
\end{proposition}
We would like to use \eqref{strichartz} to obtain an $L^{4}$ estimates
for any generic function $f$, not necessarily a linear solution. This
can be done by foliating the space $\dbR^{3}$ using dyadic level sets
$\Lambda_{j}=\{(\x,\m,\tau) / |\tau-\om|\sim 2^{j}\}.$
\begin{proposition}
Let $  \chi_{j}(\x,\m,\tau)=\chi_{j}(\tau-\om(\x,\m))$  as in
Definition \ref{spaces}. Then, for $\epsilon > 0$, and with $\spcheck$ denoting the inverse Fourier
transform,
 \begin{eqnarray}
     \label{st1}\|(\chi_{j}|\hat{f}|)\spcheck \|_{L^{4}}&\lesssim&
     2^{\frac{j}{2}}\|\hat{f}\chi_{j}\|_{L^{2}}\\
    \label{st2}  \|f\|_{L^{4}}&\lesssim&\left(\sum_{j\geq 0}
2^{(1+\epsilon)j}
     \|\hat{f}\chi_{j}\|_{L^{2}}^{2}\right)^{\half}.
\end{eqnarray}
\end{proposition}
\begin{proof}
Set $(\chi_{j}|\hat{f}|)\spcheck (x,y,t)=g_{j}(x,y,t)$. Then
$$
g_{j}(x,y,t)=\int_{\dbR^{3}}e^{i(x\xi+y\m+t\tau)}|\hat{f}|\chi_{j}
(\x,\m,\tau)
\,d\x d\m d\tau.
$$
We now observe that  $g_{j}$ can be written as an integral of linear solutions for
\eqref{livp} with appropriate initial data. More precisely,
by a simple change of variables one can write
\begin{eqnarray*}
g_{j}(x,y,t)&=&\int_{\dbR^{3}}e^{i(x\xi+y\m+t(\lambda-\om(\x,\m))}
|\hat{f}|(\x,\m,\lambda-\om(\x,\m))\chi_{j}(\lambda)\,d\x
d\m d\lambda\\
&=&\int_{\dbR}e^{it\lambda}\chi_{j}(\lambda)
\left[\int_{\dbR^{2}}e^{i(x\xi+y\m-t\om(\x,\m))}
|\hat{f}|(\x,\m,\lambda-\om(\x,\m))\,d\x d\m\right] d\lambda\\
&=&\int_{\dbR}e^{it\lambda}\chi_{j}(\lambda)S(-t)g_{\lambda}(x,y)
\,d\lambda,
    \end{eqnarray*}
where $\widehat{g_{\lambda}}(\x,\m)=|\hat{f}|(\x,\m,\lambda-
\om(\x,\m))$.
Then, by
\eqref{strichartz}
$$
\|g_{j}\|_{L^{4}}\lesssim \int_{\dbR}\chi_{j}(\lambda)
\|S(t)g_{\lambda}(x,y)\|_{L^{4}}\,d\lambda
\lesssim \int_{\dbR}\chi_{j}(\lambda)
\|\hat{f}(\x,\m,\lambda-\om(\x,\m))\|_{L^{2}_{\x,\m}}\,d\lambda,
$$
and after Cauchy-Schwarz in $\lambda$ \eqref{st1} is proved. To prove
\eqref{st2}
we first foliate the function $f$ over the dyadic levels
$|\tau-\om(\x,\m)|\sim
2^{j}$, that is we write $f(x,y,t)=\sum_{j\geq 0}f_{j}$, where
$$
f_{j}(x,y,t)=\int_{\dbR^{3}}e^{i(x\xi+y\m+t\tau)}\hat{f}
\chi_{j}(\x,\m,\tau)
\,d\x d\m d\tau.
$$
Then if we proceed as above and we use Minkowski's inequality, we obtain
$$
\|f\|_{L^{4}}\lesssim \sum_{j\geq 0}\int_{\dbR}\chi_{j}(\lambda)
\|S(t)f_{\lambda}(x,y)\|_{L^{4}}\,d\lambda
\lesssim \sum_{j\geq 0}\int_{\dbR}\chi_{j}(\lambda)
\|\hat{f}(\x,\m,\lambda-\om(\x,\m))\|_{L^{2}_{\x,\m}}\,d\lambda
$$
where $\widehat{f_{\lambda}}(\x,\m)=|\hat{f}|(\x,\m,
\lambda-\om(\x,\m))$. Using Cauchy-Schwarz, first in $\lam$, and then
in $j$, the proof is concluded.
\end{proof}
We also need a smoothing effect estimate and a matching  maximal function
estimate. We start by defining the operators
$P_{+}, P_{-}$ and $P_{0}$ such that
\begin{eqnarray*}
\widehat{P_{+}(f)}(\x,\m)&=&
\hat{f}(\x,\m)\chi_{\{|\x|>>\frac{|\m|}{|\x|}\}}(\x,\m),\\
\widehat{P_{-}(f)}(\x,\m)&=&
\hat{f}(\x,\m)\chi_{\{|\x|<<\frac{|\m|}{|\x|}\}}(\x,\m),\\
\widehat{P_{0}(f)}(\x,\m)&=&
\hat{f}(\x,\m)\chi_{\{|\x|\sim\frac{|\m|}{|\x|}\}}(\x,\m).
\end{eqnarray*}
We also recall the operator $D^{s}, \, s\geq 0$
defined through the Fourier transform as
$\widehat{D^{s}f}(\x)=|\x|^{s}\hat{f}$.  We then have the following
proposition.
\begin{proposition}\label{regions} (smoothing effect estimates)
For any $u_{0}\in L^{2}(\dbR^{2})$,
\begin{eqnarray}
 \label{smoothing+-}
 \|\partial_{x}S(t)P_{\pm}u_{0}\|_{L^{\infty}_{x}L^{2}_{y}L^{2}_{t}}
 &\leq& C\|u_{0}\|_{L^{2}},\\
 \label{smoothing0} \|D_{x}^{\half}S(t)P_{0}u_{0}
 \|_{L^{\infty}_{y}L^{2}_{x}L^{2}_{t}}
 &\leq& C\|u_{0}\|_{L^{2}}.
 \end{eqnarray}
\end{proposition}
\begin{proof}
    The proof follows the argument presented by Kenig,
Ponce and Vega in the proof of
the one-dimensional  KdV smoothing effect in \cite{KPV2}.
To prove \eqref{smoothing+-} we first define the regions of integration
\begin{eqnarray*}
A_{+}&=&\{|\x|>>\frac{|\m|}{|\x|}\}\\
A_{-}&=&\{|\x|<<\frac{|\m|}{|\x|}\}.
\end{eqnarray*}
Then we write
$$
\partial_xS(t)P_{\pm}u_0(x,y)= i\int_{A_{\pm}}
e^{i(x \xi + y \mu )}\xi e^{it\om(\xi,\m)}\widehat{u_0}
(\xi,\m)\,d\xi \,d\m,
$$
where $\om(\x,\m)= \xi^3+\frac{\m^2}{\xi}$.  We make the change of variables
$(\zeta, \m) = (\xi^3+\m^2/\xi,
\m)$, and it is to check that if $J(\x,\m)$ represents the jacobian,
then in $A_{\pm}$,
$|J(\xi,\m)|\gtrsim |\xi|^2$ holds. Now assume that
$\x=\theta(\zeta,\m)$, then
the term above equals
\[
\int_{\R^2} e^{i\m y+it\zeta}[e^{i\theta(\zeta,\m)x}
\theta(\zeta,\m)
\widehat{u_0}(\theta(\zeta,\m),\m)
\chi_{\tilde{A}_{\pm}}|J|^{-1}]\,d\zeta d\m,\]
where $\tilde{A}_{\pm}$ is the transformation of $A_{\pm}$ under the
given change of variables. Then by Plancherel's theorem
\begin{eqnarray*}
\|\partial_xS(t)P_{\pm}u_0(x)\|_{L^2_{y,t}}&=&
\|e^{i\theta(\zeta,\m)x}
\theta(\zeta,\m)\widehat{u_0}(\theta(\zeta,\m),\m)
\chi_{\tilde{A}_{\pm}}|J|^{-1}
\|_{L^2_{\zeta,\m}}\\
&=&\left(\int_{\tilde{A}_{\pm}}
|\theta(\zeta,\m)|^2|\widehat{u_0}|^2|J|^{-2}
\,d\zeta d\m
\right)^{\half}\\
&\lesssim&\left(\int|\xi|^2|\widehat{u_0}|^2|\xi|^{-2} d\xi d\mu \right)^{\half}
=\|u_0\|_{L^2_{x,y}}.
\end{eqnarray*}
To prove  \eqref{smoothing0} we use a similar argument. We set
$A_{0}=\{|\x|\sim\frac{|\m|}{|\x|}\}$, and we write
$$
D_x^{1/2}S(t)P_{0}u_0(x,y)=
\int_{A_{0}}
e^{i(\xi,\m)(x,y)}|\xi|^{\half} e^{it\om(\xi,\m)}
\widehat{u_0}(\xi,\m)\,d\xi d\m.
$$
We make the change of variables
$(\xi,\rho) = (\xi,\xi^3+\m^2/\xi)$, and we observe that this time the
estimate for the jacobian is $|J(\xi,\m)|\gtrsim |\xi|$.
We set $\m=\gamma(\xi,\rho)$ and we continue the chain of
inequalities above with
\[
\int e^{i\xi x+it\rho}[e^{i\gamma(\xi,\rho)y}|\xi|^{\half}
\widehat{u_0}(\xi,\gamma(\xi,\rho))\chi_{\tilde{A}_{0}}
|J|^{-1}]\,d\rho d\xi,\]
where $\tilde{A}_{0}$ is the transformation of $A_{0}$ under the above
change of variables. Then by Plancherel's theorem
\begin{eqnarray*}
\|D_x^{\half}S(t)P_{0}u_0(y)\|_{L^2_{x,t}}&=&
\|e^{i\gamma(\xi,\rho)y}
|\xi|^{\half}\widehat{u_0}(\xi,\gamma(\xi,\rho))
\chi_{\tilde{A}_{0}\}}|J|^{-1}\|_{L^2_{\xi,\rho}}\\
&=&\left(\int_{\R^2} |\xi||\widehat{u_0}|^2J^{-2}\,d\xi d\rho
\right)^{\half}\\
&\lesssim&\left(\int_{\R^2} |\xi||\widehat{u_0}|^2|\xi|^{-1} d\xi d\mu \right)^{\half}
=\|u_0\|_{L^2_{x,y}}.
\end{eqnarray*}
\end{proof}
Using the argument of foliation with $\chi_{j}$ introduced to prove
\eqref{st2},  one obtains
\begin{corollary}
 Let $  \chi_{j}(\x,\m,\tau)=\chi_{j}(\tau-\om(\x,\m))$  be as in
Definition \ref{spaces}. Then
\begin{eqnarray}
\label{smoothj0}\|(|\x|^{\half}\chi_{\{|\x|\sim \frac{|\m|}{|\x|}\}}
\hat f(\xi,\m,\tau)\chi_{j}
(\tau-\om(\x,\m)))\spcheck \|_{L^{\infty}_{y}L^{2}_{x}L^{2}_{t}}&\lesssim&
2^{\frac{j}{2}}\|\hat f\chi_{j}\|_{L^{2}_{\x,\m,\tau}}\\
\label{smoothj1}\|(|\x|\chi_{\{|\x|>>\frac{|\m|}{|\x|}\}}
\hat f(\xi,\m,\tau)\chi_{j}
(\tau-\om(\x,\m)))\spcheck \|_{L^{\infty}_{x}L^{2}_{y}L^{2}_{t}}&\lesssim&
2^{\frac{j}{2}}\|\hat f\chi_{j}\|_{L^{2}_{\x,\m,\tau}}\\
\label{smoothj2}\|(|\x|\chi_{\{|\x|<<\frac{|\m|}{|\x|}\}}
\hat f(\xi,\m,\tau)\chi_{j}
(\tau-\om(\x,\m)))\spcheck \|_{L^{\infty}_{x}L^{2}_{y}L^{2}_{t}}&\lesssim&
2^{\frac{j}{2}}\|\hat f\chi_{j}\|_{L^{2}_{\x,\m,\tau}}.
\end{eqnarray}
\end{corollary}
On the other hand, using interpolation with the trivial $L^{2}$ norm
estimate, we also obtain
\begin{corollary}
    Let $  \chi_{j}(\x,\m,\tau)=\chi_{j}(\tau-\om(\x,\m))$  be as in
Definition \ref{spaces}. Then
\begin{eqnarray}
 \label{smoothj5}\|(|\x|^{\frac{1}{4}}\chi_{\{|\x|\sim
\frac{|\m|}{|\x|}\}}
 \hat f(\xi,\m,\tau)\chi_{j}
(\tau-\om(\x,\m)))\spcheck \|_{L^{4}_{y}L^{2}_{x}L^{2}_{t}}&\lesssim&
2^{\frac{j}{4}}\|\hat f\chi_{j}\|_{L^{2}_{\x,\m,\tau}}\\
\label{smoothj3}\|(|\x|^{1/2}\chi_{\{|\x|>>\frac{|\m|}{|\x|}\}}
\hat f(\xi,\m,\tau)\chi_{j}
(\tau-\om(\x,\m)))\spcheck \|_{L^{4}_{x}L^{2}_{y}L^{2}_{t}}&\lesssim&
2^{j/4}\|\hat f\chi_{j}\|_{L^{2}_{\x,\m,\tau}}\\
\label{smoothj4}\|(|\x|^{1/2}\chi_{\{|\x|<<\frac{|\m|}{|\x|}\}}
\hat f(\xi,\m,\tau)\chi_{j}
(\tau-\om(\x,\m)))\spcheck\|_{L^{4}_{x}L^{2}_{y}L^{2}_{t}}&\lesssim&
2^{j/4}\|\hat f\chi_{j}\|_{L^{2}_{\x,\m,\tau}}
\end{eqnarray}
\end{corollary}
We finally introduce a maximal function estimate.
\begin{proposition} (maximal function estimate)
Let $T_{m}$ be the operator such that $\widehat{T_{m}(f)}(\x,\m,\tau)=
m(\m,\tau)\hat{f}(\x,\m,\tau)$. Then
\begin{equation}
    \|T_{m}(f)\|_{L^{2}_{x}L^{\infty}_{y,t}}\leq
    C\|m\|_{L^{2}_{\m,\tau}}\|\hat{f}\|_{L^{2}_{\x,\m,\tau}}.
\label{maxfun}\end{equation}
Similarly, if $\widehat{T_{w}(f)}(\x,\m,\tau)=
w(\x,\tau)\hat{f}(\x,\m,\tau)$. Then
\begin{equation}
    \|T_{w}(f)\|_{L^{2}_{y}L^{\infty}_{x,t}}\leq
    C\|w\|_{L^{2}_{\x,\tau}}\|\hat{f}\|_{L^{2}_{\x,\m,\tau}}.
\label{maxfun1}\end{equation}
\end{proposition}
\begin{proof}
We only prove \eqref{maxfun}. We first write
$$T_{m}f(x,y,t)=\int_{\dbR^{2}}\check{m}(y-y',t-t')f(x,y',t')
dy'dt'.$$
Then
$$|T_{m}f(x,y,t)|\lesssim
\|\check{m}\|_{L^{2}}\|f(x,\cdot,\cdot)\|_{L^{2}}.$$
To end the proof one just has to take the $L^{2}$ norm in the $x$
variable.
\end{proof}
It is also useful to observe that
interpolating  \eqref{maxfun} and \eqref{maxfun1} with the
trivial $L^{2}$ estimates, we obtain
\begin{eqnarray}
 \|T_{m}(f)\|_{L^{2}_{x}L^{4}_{y,t}}&\leq&
    C\|m\|_{L^{4}_{\m,\tau}}\|\hat{f}\|_{L^{2}_{\x,\m,\tau}},
\label{maxfunint}\\
\|T_{w}(f)\|_{L^{2}_{y}L^{4}_{x,t}}&\leq&
    C\|w\|_{L^{4}_{\x,\tau}}\|\hat{f}\|_{L^{2}_{\x,\m,\tau}}.
\label{maxfunint1}
\end{eqnarray}

We end this section with a simple weighted Sobolev inequality that will
be useful later.
\begin{lemma} (weighted Sobolev)
Assume $w(\x,\m)\gtrsim 1$ for any $(\x,\m)\in \dbR^{2}$. Then for
any $\epsilon_{0}\geq 0$, for any $p>2$ and $\theta=(p-2)/2p$,
\begin{equation}
    \|f\|_{L^{p}_{\m}}\lesssim
    \|w(\xi, \cdot)^{\epsilon_{0}}f\|_{L^{2}_{\m}}^{(1-\theta)}
    \|w(\xi, \cdot)^{-\epsilon_{0}}f'\|_{L^{2}_{\m}}^{\theta}.
\label{fp}\end{equation}
\label{sobolev}\end{lemma}
\begin{proof}
We write
\begin{eqnarray*} \half f^{2}(\m)&=&\half\int_{-\infty}^{\m}
    \frac{d}{d\m'}(f^{2}(\m'))d\m'\\
    &=&
\int_{-\infty}^{\m}f(\m')\frac{d}{d\m'}(f(\m'))d\m'=
\int_{-\infty}^{\m}w(\x,\m')^{\epsilon_{0}}f(\m'
)w(\x,\m')^{-\epsilon_{0}}
\frac{d}{d\m'}(f(\m'))d\m'\\
&\leq&\|w(\x,\cdot)^{\epsilon_{0}}f\|_{L^{2}_{\m}}
\|w(\x,\cdot)^{-\epsilon_{0}}f'\|_{L^{2}_{\m}}.
\end{eqnarray*}
From here it follows that
\begin{equation}
    \|f\|_{L^{\infty}}\leq
    \|w(\x,\cdot)^{\epsilon_{0}}
    f\|_{L^{2}_{\m}}^{\half}
\|w(\x,\cdot)^{-\epsilon_{0}}f'\|_{L^{2}_{\m}}^{\half}.
\label{finfty}    \end{equation}
If $p>2$
$$
\int_\R |f|^{p}d\m=\int_\R |f|^{p-2}|f|^{2}d\m
\leq\|f\|_{L^{\infty}}^{p-2}\|f\|_{L^{2}}^{2},
$$
and by \eqref{finfty} it follows that
$$\|f\|_{L^{p}_{\m}}\leq \|f\|_{L^{2}_{\m}}^{2/p}
\|w(\x,\cdot)^{\epsilon_{0}}
    f\|_{L^{2}_{\m}}^{(p-2)/2p}
\|w(\x,\cdot)^{-\epsilon_{0}}f'\|_{L^{2}_{\m}}^{(p-2)/2p}$$
and because $w(\x,\m)\gtrsim 1$
$$\|f\|_{L^{p}_{\m}}\leq \|w(\x,\cdot)^{\epsilon_{0}}
    f\|_{L^{2}_{\m}}^{(p+2)/2p}
    \|w(\x,\cdot)^{-\epsilon_{0}}f'\|_{L^{2}_{\m}}^{(p-2)/2p},$$
    and the lemma is proved.
\end{proof}

%%%%%%%%%%%%%%%%%%%%%%%%%%%%%%%%%%%%%%%%%%%%%%%%%%%%%%%%%%%%%%%%%%%

\section{The Bilinear Estimates}

As announced at the end of Section 1, the core of the well-posedness
result we present in this paper is contained in the following two
theorems.
\begin{theorem}
Assume $0 <\epsilon_{0}<\frac{1}{8}$. Then for any 
$\frac{1}{4}<\epsilon<1$, we have
\begin{eqnarray}
\label{Xbilinear}
\|\partial_{x}(uv)\|_{X_{1-\epsilon_{0},-\half}}&\leq& C
\|u\|_{X_{1-\epsilon_{0},\half}}(\|v\|_{X_{1-\epsilon_{0},\half}}+
\|v\|_{X_{1-\epsilon_{0},\half}}^{1-\epsilon }
\|v\|_{Y_{1-\epsilon_{0},-\epsilon_{0}, \half}}^{\epsilon })\\
&+&C
\|v\|_{X_{1-\epsilon_{0},\half}}(\|u\|_{X_{1-\epsilon_{0},\half}}+
\|u\|_{X_{1-\epsilon_{0},\half}}^{1-\epsilon }
\|u\|_{Y_{1-\epsilon_{0},-\epsilon_{0}, \half}}^{\epsilon }).
\end{eqnarray}
\label{xbil}\end{theorem}
The companion of the above  bilinear estimate is
\begin{theorem}
Assume $0 <\epsilon_{0}<\frac{1}{16}$. Then for any $\frac{1}{4}
<\epsilon<1$, we have
\begin{eqnarray}
\label{Ybilinear}
\|\partial_{x}(uv)\|_{Y_{1-\epsilon_{0},-\epsilon_{0},-\half}}&\leq&
C
\|u\|_{Y_{1-\epsilon_{0},-\epsilon_{0}, \half}}
(\|v\|_{X_{1-\epsilon_{0},\half}}+
\|v\|_{X_{1-\epsilon_{0},\half}}^{1-\epsilon }
\|v\|_{Y_{1-\epsilon_{0},-\epsilon_{0}, \half}}^{\epsilon })\\
&+&C
\|v\|_{Y_{1-\epsilon_{0},-\epsilon_{0}, \half}}
(\|u\|_{X_{1-\epsilon_{0},\half}}+
\|u\|_{X_{1-\epsilon_{0},\half}}^{1- \epsilon }
\|u\|_{Y_{1-\epsilon_{0},-\epsilon_{0}, \half}}^{\epsilon }).
\end{eqnarray}
\label{ybil}\end{theorem}
\begin{remark}
    We have not attempted to find the optimal value for
$\epsilon_{0}$
    for which our argument can be carried out.
    \end{remark}
\begin{remark}
An estimate for the bilinear expression $D_{x}^{1/2}(uv)$ in
spaces not involving weights already appeared in \cite{MST}.
    \end{remark}

To  give an idea of how the proof will be conducted we write the
left hand side of the bilinear inequality in Theorem  \ref{xbil}
using duality.

We have to estimate, for $g_j \geq 0$, 
\begin{eqnarray}
& &\label{duality-1}\sum_{j\geq 0}2^{-j/2}\sum_{m\geq 0}
\int_{A*} g_{j}(\x,\m,\tau)\theta_{m}(\x)
    \chi_{1}(\x,\m)\chi_{j}(\tau-\om(\x,\m))\\
\nonumber
& &\times |\x|\max(1,|\x|)^{1-\epsilon_{0}}
|\hat{u}|(\x_{1},\m_{1},\tau_{1})
    |\hat{v}|(\x_{2},\m_{2},\tau_{2})d\x_{1}d\m_{1}
    d\tau_1d\x_2d\m_{2}d\tau_{2}
\end{eqnarray}
and
\begin{eqnarray}
\label{duality0}& &\sum_{j\geq 0}2^{-j/2}\sum_{n\geq 0}
\int_{A*}
g_{j}(\x,\m,\tau)\theta_{n}(\m)
    \chi_{2}(\x,\m)\chi_{j}(\tau-\om(\x,\m))\\
\nonumber
& &\times |\x|\max(1, \gi)^{1-\epsilon_{0}}
|\hat{u}|(\x_{1},\m_{1},\tau_{1})
    |\hat{v}|(\x_{2},\m_{2},\tau_{2})d\x_{1}d\m_{1}
    d\tau_{1}d\x_{2}d\m_{2}d\tau_{2},
\end{eqnarray}
where $ A*$ is the set $\{\xi_{1}+\x_{2}=\x,
\m_{1}+\m_{2}=\mu, \tau_{1}+\tau_{2}=\tau \}$ and
$\|g_{j}\theta_{m}\chi_{1}\chi_{j}
\|_{L^{2}_{\x,\m,\tau}}\leq 1 $
and $\|g_{j}\theta_{n}\chi_{2}\chi_{j}
\|_{L^{2}_{\x,\m,\tau}}\leq 1 $.
It is clear that by symmetry one can always assume that
$|\x_{1}|\geq |\x_{2}|$.
Based on Remark \ref{diffregions} one can easily understand that many
different cases need to be considered in view of the fact that there
will be a combination of interactions between
``good'' regions of type $R_{g}$, bad
regions of type $R_{b}$, regions with relatively small or large
frequencies.
The whole analysis is complicated further by the fact that the spaces we
use are
anisotropic.  We start by subdividing $A*$ into six domains  of integration
\begin{eqnarray}
    \label{a1}&& A_{1}= A*\cap \{|\x_{1}|\geq |\x_{2}|, |\x_1|\leq
    1\},\\
    \label{a2}&& A_{2}=A*\cap \{|\x_{1}|\geq |\x_{2}|, |\x_1|>1,
    |\x_{2}|\sim |\x_{1}|\},\\
    \label{a3}&& A_{3}= A*\cap \{|\x_1|>1,
    1<|\x_{2}|\leq 10^{-10} |\x_{1}|, |\m_{2}|/|\x_{2}|\gtrsim
    \max{\{|\x_{1}|, |\m_{1}|/|\x_{1}|\}}\},\\
\label{a4}&& A_{4}= A*\cap B\cap \{|\x_1|>1, |\x_{2}|\leq 1,
    |\x_{2}|\leq 10^{-10} |\x_{1}|, |\m_{2}|/|\x_{2}|\gtrsim
    \max{\{|\x_{1}|, |\m_{1}|/|\x_{1}|\}}\},\\
    \label{a5}&& A_{5}= A*\cap B^{c}\cap \{|\x_1|>1, |\x_{2}|\leq 1,
    |\x_{2}|\leq 10^{-10}|\x_{1}|, |\m_{2}|/|\x_{2}|\gtrsim
    \max{\{|\x_{1}|, |\m_{1}|/|\x_{1}|\}}\},\\
    \label{a6}&& A_{6}= A*\cap \{|\x_1|>1,
    |\x_{2}|\leq 10^{-10} |\x_{1}|, |\m_{2}|/|\x_{2}|\leq 10^{-10}
    \max{\{ |\x_{1}|, |\m_{1}|/|\x_{1}|\}}\},
    \end{eqnarray}
where    $B=\{|\x_{1}|\geq \frac{1}{100}|\m_{1}|/|\x_{1}|\}$. We  also use
the auxiliary region
$$\tilde{A}_{5}(\epsilon_{0})=
      \{|\m_{1}/\x_{1}-\m_{2}/\x_{2}|^{2}<3/2
          |\x_{1}+\x_{2}|^{2},  |\xi_{1}|^{1+\alpha\epsilon_{0}}
          \leq |\m_{1}/\x_{1}|\},$$
where $\alpha$ will depend on $\epsilon_{0}$.
The most delicate part of our
estimate  occurs in  the region $A_{6}$ and it is only here
that we need the weighted spaces and the Besov type norms.
We start with a lemma.
\begin{lemma}
If $\epsilon_{0}<\frac{1}{8}$, then
\begin{eqnarray}
 \label{lemma1}& &\sum_{j_{1},j_{2}\geq 0}\sum_{j\geq
 0}2^{-j/2}\int_{A}g_{j}(\x_{1}+\x_{2},\m_{1}+\m_{2},
\tau_{1}+\tau_{2})
 \chi_{j}(\tau_{1}+\tau_{2}-\om(\x_{1}+\x_{2},
\m_{1}+\m_{2}))\\
 \nonumber& &|\x_{1}+\x_{2}|^{\epsilon_{0}}
 \Pi_{i=1,2}\phi_{i}(\x_{i},\m_{i},\tau_{i})
\chi_{j_{i}}
 (\tau_{i}-\om(\x_{i},\m_{i}))d\x_{i}d\m_{i}
d\tau_{i}\\
\nonumber&\lesssim&\sup_{j}\|g_{j}\|_{L^{2}}
\Pi_{i=1,2}\sum_{j_{i}\geq 0}
2^{j_{i}/2}\|\phi_{i}\chi_{j_{i}}\|_{L^{2}},
 \end{eqnarray}
 where $A=A_{k}, k=1,2,3,4$ or $A=A_{5}-
\widetilde{A}_{5}(\epsilon_{0})$.
\label{mainlemma}\end{lemma}
The proof that we present below gives a little more general result
than the one stated. In particular we show that \eqref{lemma1} holds
in sets larger than $A_{3}$ and $A_{4}$, namely in
$$\tilde A_{3}=A*\cap \{|\x_1|>1,
    1<|\x_{2}|\leq 10^{-10} |\x_{1}|\},$$
and
$$\tilde A_{4}=A*\cap B\cap \{|\x_1|>1, |\x_{2}|\leq 1,
    |\x_{2}|\leq 10^{-10} |\x_{1}|\}.$$
\begin{proof}
For simplicity, for $i=1,2$,  we set
$$\phi_{i}(\x_{i},\m_{i},\tau_{i})\chi_{j_{i}}
(\tau_{i}-\om(\x_{i},\m_{i}))=\phi_{i,j_{i}}(\x_{i},\m_{i},
\tau_{i}).$$
Also, whenever we use a dyadic decomposition either with respect to
$|\x_{i}|\sim 2^{m_{i}}$ or $|\m_{i}|\sim 2^{n_{i}}$, we write
\begin{eqnarray}
\label{pmj}\phi_{i,j_{i},m_{i}}(\x_{i},\m_{i},\tau_{i})&=&
\phi_{i,j_{i}}
\chi_{\{|\x_{i}|\sim 2^{m_{i}}\}}(\x_{i},\m_{i},\tau_{i})\\
\label{pnj}\phi_{i,j_{i},n_{i}}(\x_{i},\m_{i},\tau_{i})
&=&\phi_{i,j_{i}}
\chi_{\{|\m_{i}|\sim 2^{n_{i}}\}}(\x_{i},\m_{i},\tau_{i}).
\end{eqnarray}
We prove the theorem by analyzing the integral in \eqref{lemma1} on
the different regions.

\noindent
{\bf Region $A_{1}$:} Here $|\x_{1}+\x_{2}|\lesssim 1$ and we
can simply use the Strichartz inequality \eqref{st1}.

\noindent
{\bf Region $A_{2}$:}. Here we can assume also that
$|\x_{1}+\x_{2}|>1$, otherwise we go back to the argument used in the
region $A_{1}$. We dyadically decompose with
respect to $|\x_{1}|\sim 2^{m_{1}}$ (hence $|\x_{2}|\sim 2^{m_{1}}$) and
we rewrite the left hand side of \eqref{lemma1} as
\begin{equation}
 \sum_{m_{1}\geq 0}\sum_{j,j_{1},j_{2}}2^{-j/2}\int g_{j}\chi_{j}
 |\x_{1}+\x_{2}|^{\epsilon_{0}}\Pi_{i=1,2}\phi_{i,j_{i},m_{i}}
 d\tau_{i}d\x_{i}d\m_{i}
 \label{a2step1}\end{equation}
We now consider two cases:

\noindent
{\bf Case A:} $j>2\epsilon_{0}m_{1}$. We use the Strichartz inequality
\eqref{st1} and \eqref{a2step1} can be bounded by
$$\sum_{m_{1}\geq 0}\sum_{j>2\epsilon_{0}m_{1}}\sum_{j_{1},j_{2}\geq 0}
 2^{-j/2+m_{1}\epsilon_{0}}
 \|g_{j}\chi_{j}\|_{L^{2}}\Pi_{i=1,2}2^{j_{i}/2}
 \|\phi_{i,j_{i},m_{i}}\|_{L^{2}}$$
 and \eqref{lemma1} follows in this case.

\noindent
{\bf Case B:} $0\leq j\leq 2\epsilon_{0}m_{1}$. We change variable
in $\tau_{1}$ and $\tau_{2}$ by setting
$(\tau_{i}-\om(\x_{i},\m_{i}))=\theta_{i}$ and we  write the
left hand side of
\eqref{lemma1} as
\begin{eqnarray}
 \label{a2step2}& &\sum_{j_{1},j_{2}\geq 0}\sum_{j\geq
 0}2^{-j/2}\int_{A}g_{j}(\x_{1}+\x_{2},\m_{1}+\m_{2},\theta_{1}+
 \om(\x_{1},\m_{1})+\theta_{2}+\om(\x_{2},\m_{2}))\\
 \nonumber & &\chi_{j}(\theta_{1}+
 \om(\x_{1},\m_{1})+\theta_{2}+\om(\x_{2},\m_{2})
 -\om(\x_{1}+\x_{2},\m_{1}+\m_{2}))\\
 \nonumber& &|\x_{1}+\x_{2}|^{\epsilon_{0}}
 \Pi_{i=1,2}\phi_{i}(\x_{i},\m_{i},\theta_{i}+\om(\x_{i},\m_{i}))
 \chi_{j_{i}}(\theta_{i})d\x_{i}d\m_{i}d\theta_{i}.
 \end{eqnarray}
From \eqref{den} we also have that
\begin{equation}
\label{den1}1+ \left|\theta_{1}+\theta_{2}+\frac{\x_{1}
\x_{2}}{(\x_{1}+\x_{2})}
    \left(\left(\frac{\m_{1}}{\x_{1}}-\frac{\m_{2}}
{\x_{2}}\right)^{2}
    -3(\xi_{1}+\xi_{2})^{2}\right)\right|\sim 2^{j}.
\end{equation}
{\bf Case B1:}
$\left|\frac{\m_{1}}{\x_{1}}-\frac{\m_{2}}{\x_{2}}\right|^{2}
    \leq 3/2|\xi_{1}+\xi_{2}|^{2}.$\\
Then from \eqref{den1} and \eqref{den} it follows
that $|\x_{1}||\x_{2}||\x_{1}+\x_{2}|
\lesssim
2^{\max{(j_{1},j_{2},j)}}$.
Because now $|\x_{1}+\x_{2}|>1$, if $j=\max{(j_{1},j_{2},j)}$,
then
$2^{2m_{1}}\lesssim 2^{2\epsilon_{0}m_{1}}$, a contradiction if
$\epsilon_{0}<1$ for $m_1$ large enough.
So by symmetry  we can assume that $|\x_{1}+\x_{2}|
\lesssim 2^{(j_{1}-2m_{1})}$, hence
$|\x_{1}+\x_{2}|^{\epsilon_{0}}\lesssim
2^{\epsilon_{0}(j_{1}-2m_{1})}$.
We can then continue\footnote{Here one needs to take the inverse
Fourier transform of $\phi_{2}$ and of $g$ in $L^{4}$ and that of
$\phi_{1}$ in $L^{2}$.}  the estimate of \eqref{a2step2} by using
Strichartz inequality \eqref{st1} with
\begin{equation}
    \sum_{m_{1}\geq 0}\sum_{j\leq 2\epsilon_{0}m_{1}}
    \sum_{j_{1},j_{2}\geq 0}2^{-2\epsilon_{0}m_{1}}
2^{\epsilon_{0}j_{1}}
 \|g_{j}\chi_{j}\|_{L^{2}}
 \|\phi_{1,j_{1},m_{1}}\|_{L^{2}}2^{j_{2}/2}
 \|\phi_{2,j_{2},m_{1}}\|_{L^{2}}
 \label{stg}\end{equation}
and after Cauchy-Schwarz in $m_{1}$ this gives \eqref{lemma1}
provided $\epsilon_{0}<\half$.

\noindent
{\bf Case B2:}
$\left|\frac{\m_{1}}{\x_{1}}-\frac{\m_{2}}{\x_{2}}\right|^{2}
    > 3/2|\xi_{1}+\xi_{2}|^{2}.$\\
Here we consider  the following change of variables
\begin{equation}
   \begin{array}{l}
u = \x_{1}+\x_{2}\\
v = \m_{1}+\m_{2}\\
w = \om(\x_{1},\m_{1})+\om(\x_{2},\m_{2})+\theta_{1}+\theta_{2}\\
\m_{2}=\m_{2}.
\end{array}
\label{chanm}\end{equation}
The Jacobian associated to this change of variable is
\begin{equation}
    J_{\m}=3(\x_{1}^{2}-\x_{2}^{2})-
\left(\left(\frac{|\m_{1}|}{|\x_{1}|}\right)^{2}-
\left(\frac{|\m_{2}|}{|\x_{2}|}\right)^{2}\right).
    \label{jm}\end{equation}
We observe that, for fixed $\theta_{1},
\theta_{2}, \x_{1}, \xi_{2},
\m_{1}$, the set where the free variable $\mu_{2}$ can range so that
 \eqref{den1} is verified is a
union of two symmetric intervals and the length of these intervals
is small. More precisely, if we denote with $\Delta_{\mu_{2}}$ this
length, then
\begin{equation}
    \Delta_{\m_{2}}\lesssim 2^{j-m_{1}}.
\label{deltam}    \end{equation}
To see this we introduce the function
$$f(\m)=\theta_{1}+\theta_{2}+\frac{\x_{1}\x_{2}}{(\x_{1}+\x_{2})}
    \left(\left(\frac{\m_{1}}{\x_{1}}-\frac{\m}{\x_{2}}\right)^{2}
    -3(\xi_{1}+\xi_{2})^{2}\right).$$
It's easy to see that $|f'(\m)|\gtrsim |\x_{1}|$, hence \eqref{deltam}
follows. We now consider two subcases.

\noindent
{\bf Case B2a:}
$\left|3(\x_{1}^{2}-\xi_{2}^{2})-\left(\left(\frac{\m_{1}}
{\x_{1}}\right)^{2}-
\left(\frac{\m_{2}}{\x_{2}}\right)^{2}\right)\right|> 1$. \\
We make the
change of variable \eqref{chanm} (now  $|J_{\m}|>1$).
Denote with $H(u,v,w,\m_{2},\theta_{1},\theta_{2})$ the transformation
of
$$\Pi_{i=1,2}\phi_{i, j_i, m_1 }(\x_{i},\m_{i},\theta_{i}+
\om(\x_{i},\m_{i}))
 \chi_{j_{i}}(\theta_{i})$$
under the above change of variables. Then \eqref{a2step2} becomes
\begin{eqnarray*}
 & &\sum_{m_{1},j_{1},j_{2}\geq 0}\sum_{0\leq j\leq
2\epsilon_{0}m_{1}}
 2^{-j/2}2^{\epsilon_{0}m_{1}}\int g_{j}\chi_{j}(u,v,w)\\
 & &|J_{\m}|^{-1}H(u,v,w,\m_{2},\theta_{1},\theta_{2})dudvdw
 d\m_{2}d\theta_{1}d\theta_{2}\\
 &\lesssim &\sum_{m_{1},j_{1},j_{2}\geq 0}
 \sum_{0\leq j\leq 2\epsilon_{0}m_{1}}
 2^{-j/2}2^{\epsilon_{0}m_{1}}2^{(j-m_{1})/2}\\
 &&\int g_{j} \chi_{j}(u,v,w)\left(\int_{\m_{2}}|J_{\m}|^{-2}
 H^{2}(u,v,w,\m_{2},\theta_{1},\theta_{2})d\m_{2}
\right)^{\half}dudvdw
 d\theta_{1}d\theta_{2}.
 \end{eqnarray*}
 Now we observe that by Cauchy-Schwarz and the inverse change of
 variable we have
\begin{eqnarray*}
& &\int g_{j} \chi_{j}(u,v,w)\left(\int_{\m_{2}}|J_{\m}|^{-2}
 H^{2}(u,v,w,\m_{2},\theta_{1},\theta_{2})d\m_{2}
\right)^{\half}dudvdw
 d\theta_{1}d\theta_{2}\\
&\leq& {{\| g_j \chi_j \|}_{L^2}} \left( \int \left( \int |J_\mu|^{-2}
 H^2 ( u ,v,w,\mu_2, \theta_1, \theta_2 ) du dv dw d\mu_2
 \right)^\half d\theta_1 d\theta_2 \right) \\
&\lesssim & {{\| g_j \chi_j \|}_{L^2}}  \left( \int \left( \int |J_\mu|^{-1}
 H^2 ( u ,v,w,\mu_2, \theta_1, \theta_2 ) du dv dw d\mu_2
 \right)^\half d\theta_1 d\theta_2 \right) \\
&\lesssim & {{\| g_j \chi_j \|}_{L^2}}
\left( \int \left( \int \Pi_{i=1,2} \phi_{i, j_i, m_1}^2 (\xi_i ,
 \mu_i , \theta_i + \omega(\xi_i , \mu_i ) ) \chi_{j_i} (\theta_i ) d
 \xi_1 d\xi_2 d\mu_1 d\mu_2 \right)^\half d\theta_1 d\theta_1 \right)
 \\
&\lesssim & {{\| g_j \chi_j \|}_{L^2}} 2^{j_1 /2 } 2^{j_2 /2} \left(
 \int \int \Pi_{i=1,2} \phi_{i, j_i , m_1}^2 ( \xi_i , \mu_i ,
 \theta_i + \omega( \xi_i , \mu_i )) d\xi_1 d\xi_2 d\mu_1 d\mu_2
 d\theta_1 d\theta_2 \right)^\half \\ 
&\lesssim& \|g_{j}\chi_{j}\|_{L^{2}}
 \Pi_{i=1,2}2^{j_{i}/2}\|\phi_{i,j_{i},m_{1}}\|_{L^{2}}
 \end{eqnarray*}
which inserted above, after a sum on $j$  gives
$$
\sum_{m_{1},j_{1},j_{2}\geq 0}
[1+ (2\epsilon_{0}m_{1})]
 2^{\epsilon_{0}m_{1}}2^{-m_{1}/2}
 \|g_{j}\chi_{j}\|_{L^{2}}
 \Pi_{i=1,2}2^{j_{i}/2 }\|\phi_{i,j_{1},m_{i}}\|_{L^{2}}$$
and from here we obtain again \eqref{lemma1} for
$\epsilon_{0}<\half$.

\noindent
{\bf Case B2b:}
$\left|3(\x_{1}^{2}-\xi_{2}^{2})-\left(\left(\frac{\m_{1}}
{\x_{1}}\right)^{2}-
\left(\frac{\m_{2}}{\x_{2}}\right)^{2}\right)\right|\leq 1$. \\
In this case the change of variables  above cannot be used because
the Jacobian may become zero. We consider instead the
change  of variables in which we leave $\x_{1}$ free:
\begin{equation}
    \begin{array}{l}
u = \x_{1}+\x_{2}\\
v = \m_{1}+\m_{2}\\
w = \om(\x_{1},\m_{1})+\om(\x_{2},\m_{2})+\theta_{1}+\theta_{2}\\
\x_{1}=\x_{1}.
\end{array}
\label{chanx}\end{equation}
In this case the Jacobian $J_{\x}$ is given by
\begin{equation}
    J_{\x}=
\frac{\m_{1}}{\x_{1}}-\frac{\m_{2}}{\x_{2}},
\label{jx}\end{equation}
and because we are in Case B2, it follows that
$$|J_{\x}|\gtrsim |\x_{1}+\x_{2}|>1.$$
We observe that, for fixed $\theta_{1}, \theta_{2},  \xi_{2},
\m_{1}, \m_{2}$, the set where the free variable $\x_{1}$ can range so
that we remain in Case B is a
union of two symmetric intervals and the length of these interval
is small. More precisely, if we denote with $\Delta_{\x_{1}}$ this
length, then
\begin{equation}
    \Delta_{\x_{1}}\lesssim 2^{-m_{1}}.
\label{deltax}    \end{equation}
To see this we introduce the function
$$h(\x)=3(\x^{2}-\xi_{2}^{2})-\left(\left(\frac{\m_{1}}
{\x}\right)^{2}-
\left(\frac{\m_{2}}{\x_{2}}\right)^{2}\right).$$
We compute
$$h'(\x)=6\x-2(\m_{1}/\x)(-\m_{1}/\x^{2})=6\x+
2(\m_{1}/\x)^{2}\x^{-1}$$
and we notice that $h'(\x)$ has the same sign as $\x$, hence
$|h'(\x)|\gtrsim |\x|$, and \eqref{deltax}
follows.
Again denote with $H(u,v,w,\x_{1},\theta_{1},\theta_{2})$ the
transformation of\footnote{Here we are using the notation in
\eqref{pmj}.}
$\Pi_{i=1,2}\phi_{i,j_{i}}(\x_{i},\m_{i},\theta_{i})$
under the change of variables \eqref{chanx}. Then \eqref{a2step2} becomes
\begin{eqnarray*}
 & &\sum_{m_{1},j_{1},j_{2}\geq 0}\sum_{0\leq j\leq
2\epsilon_{0}m_{1}}
 2^{-j/2}2^{\epsilon_{0}m_{1}}\int g_{j}\chi_{j}(u,v,w)\\
 & &\times |J_{\x}|^{-1}H(u,v,w,\x_{1},\theta_{1},
\theta_{2})dudvdw
 d\x_{1}d\theta_{1}d\theta_{2}\\
 &\lesssim &\sum_{m_{1},j_{1},j_{2}\geq 0}
 \sum_{0\leq j\leq 2\epsilon_{0}m_{1}}
 2^{-j/2}2^{\epsilon_{0}m_{1}}2^{-m_{1}/2}\\
 &&\int g_{j} \chi_{j}(u,v,w)\left(\int_{\x_{1}}|J_{\x}|^{-2}
 H^{2}(u,v,w,\x_{1},\theta_{1},\theta_{2})d\x_{1}
\right)^{\half}dudvdw
 d\theta_{1}d\theta_{2}\\
 &\lesssim & \sum_{m_{1},j_{1},j_{2}\geq 0}
 \sum_{0\leq j\leq 2\epsilon_{0}m_{1}}
 2^{-j/2}2^{\epsilon_{0}m_{1}}2^{-m_{1}/2}\Pi_{i=1,2}2^{j_{i}/2}
 \|\phi_{i,j_{i},m_{i}}\|_{L^{2}}
\end{eqnarray*}
and this again gives \eqref{lemma1} for $\epsilon_{0}<\half$.

\noindent
{\bf Region $A_{3}$:}  In this region, (see \eqref{a3}),
$|\x_{1}+\x_{2}|\sim|\x_{1}|$. We dyadically decompose with respect to
$|\x_{1}|\sim 2^{m_{1}}$ (hence $|\x_{1}+\x_{2}|\sim 2^{m_{1}}$). The
left hand side of \eqref{lemma1} now becomes
\begin{equation}
\sum_{m_{1}\geq 0}\sum_{j,j_{1},j_{2}\geq 0}2^{-j/2}2^{m_{1}\epsilon_{0}}
\int g_{j,m_{1}}\chi_{j}\phi_{1,j_{1},m_{1}}\phi_{2,j_{1}} d\xi_i
d\mu_i d\tau_i 
\label{a3step1}
\end{equation}
where the arguments of the functions as in \eqref{lemma1}.
We consider two subcases. For $0 < \delta \ll 1$, to be fixed, we have:\\
{\bf Case A:} $j>(2+2\delta)\epsilon_{0}m_{1}, 0<\delta<<1$.
We use Strichartz inequality \eqref{st1} and we
obtain
\begin{eqnarray*}
    \eqref{a3step1}&\lesssim&
\sum_{m_{1}\geq 0}\sum_{j_{1},j_{2}\geq 0}
\sum_{j\geq 2(1+\delta)\epsilon_{0}m_{1}}2^{-j/2}2^{m_{1}\epsilon_{0}}
\|g_{j,m_{1}}\chi_{j}\|_{L^{2}}2^{j_{1}/2}2^{j_{2}/2}
\|\phi_{1,j_{1},m_{1}}\|_{L^{2}}\|\phi_{2,j_{2}}\|_{L^{2}}\\
&\lesssim&
\sum_{m_{1}\geq 0}\sum_{j_{1},j_{2}\geq 1}2^{-m_{1}\delta\epsilon_{0}}
\sup_{j}\|g_{j}\chi_{j}\|_{L^{2}}2^{j_{1}/2}2^{j_{2}/2}
\|\phi_{1,j_{1},m_{1}}\|_{L^{2}}\|\phi_{2,j_{2}}\|_{L^{2}}.
\end{eqnarray*}
{\bf Case B:} $j\leq (2+2\delta)\epsilon_{0}m_{1}$.\\
{\bf Case B1:}
$\left|\frac{\m_{1}}{\x_{1}}-\frac{\m_{2}}{\x_{2}}\right|^{2}
    \leq 3/2|\xi_{1}+\xi_{2}|^{2}.$\\
As in region $A_{2}$,  it follows that $|\x_{1}||\x_{2}|
|\x_{1}+\x_{2}|\lesssim
2^{\max{(j_{1},j_{2},j)}}$ and since $|\x_{2}|>1$ and
$|\x_{1}|\sim |\x_{1}+\x_{2}|$, we obtain that $|\x_{1}|^{2}\lesssim
2^{\max{(j_{1},j_{2},j)}}$.
If $j=\max{(j_{1},j_{2},j)}$, then
$2^{2m_{1}}\lesssim 2^{(2+2\delta)\epsilon_{0}m_{1}}$, a contradiction if
$\epsilon_{0}<\frac{1}{1+\delta}$ for $m_1$ large enough. Assume then that $j_{1}=
\max{(j_{1},j_{2},j)}$.
It follows that
$|\x_{1}+\x_{2}|^{\epsilon_{0}}\lesssim 2^{\epsilon_{0}(j_{1}-m_{1})}$
and thanks to Strichartz\footnote{Here one needs to take the anti
Fourier transform of $g$ and of $\phi_{2}$ in $L^{4}$ and that of
$\phi_{1}$ in $L^{2}$.}  inequality \eqref{st1}
we can continue the chain of inequalities with
\begin{eqnarray*}
& &\sum_{m_{1}\geq 0}\sum_{j_{1},j_{2}\geq 0}
\sum_{j\leq 2(1+\delta)\epsilon_{0}m_{1}}2^{-j/2}
\|g_{j,m_{1}}\chi_{j}\|_{L^{2}}2^{j/2} 2^{\epsilon_{0}(j_{1}-m_{1})}
\|\phi_{1,j_{1},m_{1}}\|_{L^{2}}
 2^{j_{2}/2} \|\phi_{2,j_{2}}\|_{L^{2}}\\
    & \lesssim &\sum_{m_{1}\geq 0}\sum_{j_{1},j_{2}\geq 0}
 2^{\epsilon_{0}(j_{1}-m_{1})}[(2+2\delta) \epsilon_0 m_1 +1] 2^{j_2 /2}
 (\sup_{j}\|g_{j}\chi_{j}\|_{L^{2}})
 \|\phi_{1,j_{1},m_{1}}\|_{L^{2}}
 \|\phi_{2,j_{2}}\|_{L^{2}},
\end{eqnarray*}
and this gives \eqref{lemma1} provided $\epsilon_{0}<\half$. Clearly the
case when $j_{2}=\max{(j_{1},j_{2},j)}$ is similar.

\noindent
{\bf Case B2:} $\left|\frac{\m_{1}}{\x_{1}}-\frac{\m_{2}}
{\x_{2}}\right|^{2}
    > 3/2|\xi_{1}+\xi_{2}|^{2}.$\\
As we did for region $A_{2}$ here also we consider two subcases.

\noindent
{\bf Case B2a:}
$\left|3(\x_{1}^{2}-\xi_{2}^{2})-\left(\left(\frac{\m_{1}}
{\x_{1}}\right)^{2}-
\left(\frac{\m_{2}}{\x_{2}}\right)^{2}\right)\right|> 1$. \\
We make the
change of variable \eqref{chanm}, for which  now   $|J_{\m}|>1$
and we observe that also in this region \eqref{deltam} holds.
Then \eqref{a3step1} can be bounded by
\begin{eqnarray*}
 & &\sum_{m_{1},j_{1},j_{2}\geq 0}
 \sum_{0\leq j\leq (2+2\delta) \epsilon_{0}m_{1}}
 2^{-j/2}2^{\epsilon_{0}m_{1}}2^{(j-m_{1})/2} \|g_{j}
\chi_{j}\|_{L^{2}}
 2^{j_{1}/2}2^{j_{2}/2}\|\phi_{1,j_{1},m_{1}}\|_{L^{2}}
 \|\phi_{2,j_{2}}\|_{L^{2}}\\
&\lesssim& \sum_{m_{1},j_{1},j_{2}\geq 0}
 2^{(\epsilon_{0}-\half)m_{1}}((2+2\delta) \epsilon_{0}m_{1})
(\sup_{j} \| g_{j}\chi_{j}\|_{L^{2}})
 2^{j_{1}/2}2^{j_{2}/2}\|\phi_{1,j_{1},m_{1}}\|_{L^{2}}
 \|\phi_{2,j_{2}}\|_{L^{2}}
 \end{eqnarray*}
which again gives \eqref{lemma1} for $\epsilon_{0}<\half$.

\noindent
{\bf Case B2b:}
$\left|3(\x_{1}^{2}-\xi_{2}^{2})-\left(\left(\frac{\m_{1}}
{\x_{1}}\right)^{2}-
\left(\frac{\m_{2}}{\x_{2}}\right)^{2}\right)\right|\leq 1$.\\
We consider now the
change  of variables \eqref{chanx} and we observe that
$$
    \label{jacx}|J_{\x}|=\left|
\frac{\m_{1}}{\x_{1}}-\frac{\m_{2}}{\x_{2}}\right|
\gtrsim |\x_{1}+\x_{2}|\sim 2^{m_{1}}.
$$
We also remark that in this region
\eqref{deltax} holds too.  Repeating the argument in Case B2b
of region $A_{2}$ the left hand side
of \eqref{lemma1} can be bounded by
$$\sum_{m_{1},j_{1},j_{2}\geq 0}
 \sum_{0\leq j\leq (2+2\delta) \epsilon_{0}m_{1}}
 2^{-j/2}2^{\epsilon_{0}m_{1}}2^{-m_{1}/2}
\|g_{j}\chi_{j}\|_{L^{2}}
 2^{j_{1}/2}2^{j_{2}/2}\|\phi_{1,j_{1},m_{1}}\|_{L^{2}}
 \|\phi_{2,j_{2}}\|_{L^{2}},$$
and this concludes the estimate in $A_{3}$ for $\epsilon_{0}<\half$.

\noindent
{\bf Region $A_{4}$}. Notice that we only need to restrict the
proof to the case when
$0\leq j\leq (2+2\delta)\epsilon_{0}m_{1}$ and
$\left|\frac{\m_{1}}{\x_{1}}-\frac{\m_{2}}{\x_{2}}\right|^{2}
    \leq 3/2|\xi_{1}+\xi_{2}|^{2}$, since in the other situations (Case
    A in $A_3$ and case $B_2$ in $A_3$) and we
    didn't use the assumption $|\x_{2}|\geq 1$. Observe that by the
    above restriction we also have that
\begin{equation}
    \gib\lesssim \max(\gia, |\x_{1}|)
\label{m2simx1}\end{equation}
We consider two cases.\\
{\bf Case A:} $|\x_{2}||\x_{1}|^{1+\alpha\epsilon_{0}}\geq 1$, for
some $\alpha>0$ to be determined later. \\
Going
back to the argument presented in Case B1 in region $A_{3}$, we obtain
$|\x_{1}|^{1-\alpha\epsilon_{0}}\lesssim 2^{\max(j_{1},j_{2},j)}$. If
we let
\begin{equation}
\epsilon_{0}<\frac{1}{2+2\delta+\alpha},
\label{epsi}\end{equation}
then
$\max(j_{1},j_{2},j)=\max(j_{1},j_{2})$. Let's assume
$\max(j_{1},j_{2})=j_{1}$ and  $\theta>0$ and small.
Then $2^{\epsilon_{0}m_{1}}\lesssim
2^{j_{1}\delta_{0}}2^{-\sigma_{0}m_{1}}$, where
$$\delta_{0}=\frac{\epsilon_{0}}{(1-\alpha\epsilon_{0})(1-\theta)}, \,
\mbox{ and } \, \, \sigma_{0}=\frac{\epsilon_{0}\theta}{1-\theta}.$$
Notice that if $0<\theta<<1$, from \eqref{epsi} it follows that
$\delta_{0}\leq \half$. We then use Strichartz and we bound
the left hand side of \eqref{a2step2} with
\begin{eqnarray*}
   && \sum_{m_{1},j_{1},j_{2}\geq 0}
 \sum_{0\leq j\leq (2+2\delta) \epsilon_{0}m_{1}}
 2^{-j/2}2^{\delta_{0}j_{1}}2^{-\sigma_{0}m_{1}}
2^{j/2}(\sup_{j}\|g_{j}\chi_{j}\|_{L^{2}})
 2^{j_{2}/2}\|\phi_{1,j_{1},m_{1}}\|_{L^{2}}
 \|\phi_{2,j_{2}}\|_{L^{2}}\\
&\lesssim&\sum_{m_{1},j_{1},j_{2}\geq 0}
 2^{\delta_{0}j_{1}}2^{-\sigma_{0}m_{1}}
(2+2\delta) \epsilon_{0}m_{1}(\sup_{j}\|g_{j}\chi_{j}\|_{L^{2}})
 2^{j_{2}}/2\|\phi_{1,j_{1},m_{1}}\|_{L^{2}}
 \|\phi_{2,j_{2}}\|_{L^{2}}.
\end{eqnarray*}
The result is given by the Cauchy-Schwarz inequality in $m_{1}$.
The argument when $\max(j_{1},j_{2})=j_{2}$ is similar. \\
{\bf Case B:} $|\x_{2}||\x_{1}|^{1+\alpha\epsilon_{0}}<1$. \\ From
\eqref{m2simx1} and the definition of $A_{4}$ in \eqref{a4}, we also have
$|\m_{2}|\leq |\x_{2}||\x_{1}|\lesssim
|\x_{2}||\x_{1}|^{1+\alpha\epsilon_{0}}\lesssim 1$.
We consider two subcases.\\
{\bf Case B1:} $\gia\leq \frac{1}{100} |\x_{1}|$.\\
Here we use  smoothing
effect  and the maximal function inequalities.
Let's start by assuming that
$|\m_{2}|^{2}/|\x_{2}|\lesssim 2^{j_{2}}$. Then
$|\om(\x_{2},\m_{2})|, |\tau_{2}|, \lesssim 2^{j_{2}}, |\m_{2}|
\lesssim  1$. We now set $m(\m_{2},\tau_{2})=\chi(\m_{2},\tau_{2})$,
the characteristic function of the projection of the region of
integration onto the $\m_{2}-\tau_{2}$ plane. Then
$$\|m\|_{L^{4}_{\m_{2},\tau_{2}}}\sim
2^{j_{2}/4}.$$
We then use Plancherel, H\"older with the three spaces
$L^{4}_{x,y,t}-L^{4}_{x}L^{2}_{y}L^{2}_{t}-L^{2}_{x}L^{4}_{y}L^{4}_{t}$
and inequalities \eqref{smoothj3} and \eqref{maxfunint}
to bound \eqref{a2step2} with
\begin{eqnarray*}
   && \sum_{m_{1},j_{1},j_{2}\geq 0}
 \sum_{0\leq j\leq (2+2\delta) \epsilon_{0}m_{1}}2^{-j/2}
 2^{(\epsilon_{0}-\half)m_{1}}2^{j/2}\|g_{j}\chi_{j}\|_{L^{2}}
 2^{j_{1}/4}\|\phi_{1,j_{1},m_{1}}\|_{L^{2}}
 2^{j_{2}/4} \|\phi_{2,j_{2}}\|_{L^{2}}\\
&\lesssim&\sum_{m_{1},j_{1},j_{2}\geq 0}
2^{(\epsilon_{0}-\half)m_{1}}
[1+(2+2\delta)\epsilon_0 m_1] (\sup_{j}\|g_{j}\chi_{j}\|_{L^{2}})
 2^{j_{2}/4}2^{j_{1}/4}\|\phi_{1,j_{1},m_{1}}\|_{L^{2}}
 \|\phi_{2,j_{2}}\|_{L^{2}}
\end{eqnarray*}
and if $\epsilon_{0}<\half$ the lemma is proved also in this case.
If
$|\m_{2}|^{2}/|\x_{2}|>>2^{j_{2}}$, then $
|\om(\x_{2},\m_{2})|\sim |\m_{2}|^{2}/|\x_{2}|,
|\tau_{2}|\sim  |\m_{2}|^{2}/|\x_{2}|$ and
$|\tau_{2}-\om(\x_{2},\m_{2})|\sim  2^{j_{2}}$. From \eqref{m2simx1}
we have
\begin{equation}
    |\m_{2}|\gib\lesssim |\x_{2}||\x_{1}|^{2}\lesssim |\x_{2}|
|\x_{1}|^{1+\alpha\epsilon_{0}}|\x_{1}|^{1-\alpha\epsilon_{0}}\lesssim
1\cdot|\x_{1}|^{1-\alpha\epsilon_{0}}.
\label{m2x}\end{equation}
In this case
$$\|m\|_{L^{4}_{\m_{2},\tau_{2}}}\sim
2^{m_{1}(1-\alpha\epsilon_{0})/4)}.$$
We then use again Plancherel, H\"older with the three spaces
$L^{4}_{x,y,t}-L^{4}_{x}L^{2}_{y}L^{2}_{t}-L^{2}_{x}L^{4}_{y}L^{4}_{t}$
and inequalities \eqref{smoothj3} and \eqref{maxfunint}
to bound the left hand side of \eqref{a3step1} with
\begin{eqnarray*}
   && \sum_{m_{1},j_{1},j_{2}\geq 0}
 \sum_{0\leq j\leq (2+2\delta) \epsilon_{0}m_{1}}
 2^{(\epsilon_{0}-\half)m_{1}}\|g_{j}\chi_{j}\|_{L^{2}}
 2^{j_{1}/4}\|\phi_{1,j_{1},m_{1}}\|_{L^{2}}
 2^{m_{1}(1-\alpha\epsilon_{0})/4}
\|\phi_{2,j_{2}}\|_{L^{2}}\\
&\lesssim&\sum_{m_{1},j_{1},j_{2}\geq 0}
[1+(2+2\delta) \epsilon_0 m_1] (\sup_{j}\|g_{j}
\chi_{j}\|_{L^{2}})
 2^{(\epsilon_{0}-\half)m_{1}}2^{m_{1}
(1-\alpha\epsilon_{0})/4)}
 2^{j_{1}/4}\|\phi_{1,j_{1},m_{1}}\|_{L^{2}}
 \|\phi_{2,j_{2}}\|_{L^{2}}
\end{eqnarray*}
and if $\alpha>4$, for any $\epsilon_{0}$, the sum with
respect to $m_{1}$ can be done and
then \eqref{lemma1} is proved also in this case.\\
{\bf Case B2:} $\frac{1}{100}|\x_{1}|\leq \gia \leq 100 |\x_{1}|$. \\
Here we use \eqref{smoothj5}.
Let's start by assuming that
$|\m_{2}|^{2}/|\x_{2}|\leq 2^{j_{2}}$. Clearly now
$|\tau_{2}| \lesssim 2^{j_{2}}$ and $|\x_{2}|\lesssim  1$.
We  use Plancherel, H\"older with the three spaces
$L^{4}_{x,y,t}-L^{4}_{y}L^{2}_{x}L^{2}_{t}-L^{2}_{y}L^{4}_{x}L^{4}_{t}$
and inequalities \eqref{smoothj5} and \eqref{maxfunint1}
to bound the left hand side of \eqref{lemma1} with
$$
\sum_{m_{1},j_{1},j_{2}\geq 0}
[1+(2+2\delta) \epsilon_0 m_1 ](\sup_{j}\|g_{j}\chi_{j}\|_{L^{2}})
 2^{(\epsilon_{0}-\frac{1}{4})m_{1}}2^{j_{1}/4}\|\phi_{1,j_{1},m_{1}}\|_{L^{2}}
 2^{j_{2}/4} \|\phi_{2,j_{2}}\|_{L^{2}}
$$
and if $\epsilon_{0}<\frac{1}{4}$
then \eqref{lemma1} is proved also in this case.
Assume now that $|\m_{2}|^{2}/|\x_{2}|>>2^{j_{2}}$. Then by
\eqref{m2x}
$|\tau_{2}|\lesssim |\m_{2}|^{2}/|\x_{2}|\lesssim
|\x_{1}|^{1-\alpha\epsilon_{0}}$, so that
$\|w\|_{L^{2}_{\x_{2},\tau_{2}}}\lesssim
2^{m_{1}(1-\alpha\epsilon_{0})/4}$.
We  use Plancherel, H\"older with the three spaces
$L^{4}_{x,y,t}-L^{4}_{y}L^{2}_{x}L^{2}_{t}-L^{2}_{y}L^{4}_{x}L^{4}_{t}$
and inequalities \eqref{smoothj5} and \eqref{maxfunint1}
to bound the left hand side of \eqref{a3step1} with
$$
\sum_{m_{1},j_{1},j_{2}\geq 0}
[1+(2+2\delta) \epsilon_0 m_1] (\sup_{j}\|g_{j}\chi_{j}\|_{L^{2}})
 2^{(\epsilon_{0}-\frac{1}{4})m_{1}}2^{j_{1}/4}\|\phi_{1,j_{1},m_{1}}\|_{L^{2}}
 (2^{j_{2}/4} +2^{m_{1}(1-\alpha\epsilon_{0})/4})
 \|\phi_{2,j_{2}}\|_{L^{2}}
$$
and if $\alpha>4$, then the lemma follows in this case too.

\noindent
{\bf Region $A_{5}-\tilde{A}_{5}(\epsilon_{0})$}. Also in this case
we can assume that
$0\leq j\leq (2+2\delta)\epsilon_{0}m_{1}$, and
$\left|\frac{\m_{1}}{\x_{1}}-\frac{\m_{2}}{\x_{2}}\right|^{2}
\leq 3/2|\xi_{1}+\xi_{2}|^{2}$, since in the other situation we
didn't use the assumption $|\x_{2}|\geq 1$. Also notice that here
too \eqref{m2simx1} holds, hence $\gia\sim \gib$. Because we are
in $A_{5}-\tilde{A}_{5}(\epsilon_{0})$, to these
restrictions we have to add $|\x_{1}|^{1+\alpha\epsilon}>>\gia$, which
in turn gives
\begin{equation}
    |\m_{2}|<<|\x_{1}|^{1+\alpha\epsilon_0}|\x_{2}|.
\label{m21}\end{equation}
We consider two subcases.\\
{\bf Case A:} $|\x_{1}|^{1+\alpha\epsilon_0}|\x_{2}|\geq 1$. \\
This case is identical to Case A in region $A_{4}$. \\
{\bf Case B:} $|\x_{1}|^{1+\alpha\epsilon_0}|\x_{2}|< 1$.\\
By \eqref{m21}
we now have that $|\m_{2}|\lesssim 1$. If $|\m_{2}|^{2}/|\x_{2}|\leq
2^{j_{2}}$ we can use the same
arguments presented in the first part of Case B1 in region $A_{4}$,
by replacing
\eqref{smoothj3} with \eqref{smoothj4}. If $|\m_{2}|^{2}/|\x_{2}|>>
2^{j_{2}}$, we have that $|\tau_{2}|\sim |\m_{2}|^{2}/|\x_{2}|$ and
$|\tau_{2}-\om(\x_{2},\m_{2})|\sim 2^{j_{2}}$. We estimate
$$|\m_{2}||\m_{2}|/|\x_{2}|\lesssim |\x_{2}|(|\m_{1}|/|\x_{1}|)^{2}
\lesssim |\x_{2}||\x_{1}|^{1+\alpha\epsilon_{0}}
|\x_{1}|^{1+\alpha\epsilon_{0}}\lesssim 1\cdot
|\x_{1}|^{1+\alpha\epsilon_{0}}$$
and from here on we can proceed like in the second part of
Case B1 in region $A_{4}$, again by
replacing \eqref{smoothj3} with \eqref{smoothj4}. 
We then obtain: 
$$\sum_{m_1, j_1, j_2 \geq 0} \sum_{0 \leq j \leq
    (2+2\delta) \epsilon_0 m_1} ( \sup_j {{\| g_j \chi_j
        \|}_{L^2}} ) 2^{(\epsilon_0 - \half) m_1} 2^{m_1 (1 + \alpha
      \epsilon_0) /4} 2^{j_1 /4} {{\| \phi_{1, j_1 , m_1 } \|}_{L^2}}
    {{\| \phi_{2, j_2} \|}_{L^2}},$$
    and if $(\epsilon_0 - \half) +
    (1+\alpha \epsilon_0 )/4 <0 $, or $\epsilon_0 (1+ \alpha /4 ) <
    \frac{1}{4},$ we can sum in $m_1$. This is always possible for
    $\epsilon_0 < \frac{1}{8}$, with some $\alpha > 4$.
The analysis of this
case concludes the proof of the lemma.
\end{proof}
We are now ready to prove Theorem \ref{xbil}.
To make the presentation more clear we summarize below
the main cases  considered in our
analysis \footnote{Recall that here $\x=\x_{1}+\x_{2}$ and
$\m=\m_{1}+\m_{2}$.}.
 \begin{itemize}
  \item   Region $A_{1}$
  \begin{itemize}
      \item Case A: $|\x|\geq \half|\m|/|\x|$
      \item Case B: $|\x|< \half|\m|/|\x|$
      \begin{itemize}
          \item Case B1: $|\m_{1}|\leq |\m_{2}|$
          \item Case B2: $|\m_{1}|> |\m_{2}|$
          \end{itemize}
  \end{itemize}
  \item Region $A_{2}$
  \begin{itemize}
      \item Case A: $|\x|\geq \half|\m|/|\x|$
      \item Case B: $|\x|< \half|\m|/|\x|$
      \begin{itemize}
          \item Case B1: $|\m_{1}|\leq |\m_{2}|$
          \item Case B2: $|\m_{1}|> |\m_{2}|$
      \end{itemize}
  \end{itemize}
  \item Region $A_{3}\cup A_{4}\cup (A_{5}-\tilde{A}_{5}(\epsilon_{0}))$
  \begin{itemize}
      \item Case A: $|\x|\geq \half |\m|/|\x|$
      \begin{itemize}
          \item Case A1: $|\x_{1}|\geq \half\gia$
          \item Case A2: $|\x_{1}|\leq \half \gia$
          \end{itemize}
      \item Case B: $|\x|< 1/2|\m|/|\x|$
      \begin{itemize}
          \item Case B1: $|\m_{1}|\leq |\m_{2}|$
          \item Case B2: $|\m_{1}|> |\m_{2}|$
      \end{itemize}
  \end{itemize}
  \item Region $A_{5}\cap \tilde{A}_{5}(\epsilon_{0})$
  \begin{itemize}
      \item Case A: $|\x|\geq \half |\m|/|\x|$
      \item Case B: $|\x|< \half |\m|/|\x|$
 \end{itemize}
 \item Region $A_{6}$
  \begin{itemize}
      \item Case A: $|\x_{1}|\geq 10^{2}\gia$
      \item Case B: $|\x_{1}|\leq  10^{-2}\gia$
      \item Case C: $10^{-2}\gia \leq |\x_{1}|\leq  10^{2}\gia.$
      \end{itemize}
\end{itemize}
\begin{proof}[{Proof of Theorem \ref{xbil}}]
 We reexpress the left hand side of \eqref{Xbilinear}
in Theorem \ref{xbil} using duality
and we obtain \eqref{duality-1} and \eqref{duality0}. We analyze
these expressions on the  regions described in
\eqref{a1}-\eqref{a6}. For the estimates in the regions $A_{1}$
through $A_{5}$
we find it convenient to normalize the functions $u$ and $v$ so that
the expression in the right hand side of the bilinear inequality
involves only   $L^{2}$ norms. So define
\begin{eqnarray}
 \label{normalizeu} \phi_{1,j_{1}}(\x,\m,\tau)&=&
 \max(1,|\x|,|\m|/|\x|)^{1-\epsilon_{0}}|\hat{u}|
 \chi_{j_{1}}(\x,\m,\tau),\\
\label{normalizev}
\phi_{2,j_{2}}(\x,\m,\tau)&=&\max(1,|\x|,|\m|/|\x|)^{1-\epsilon_{0}}
|\hat{v}|
 \chi_{j_{2}}(\x,\m,\tau).
\end{eqnarray}
If we use the identities
$\m=\m_{1}+\m_{2}$ and $\x=\x_{1}+\x_{2}, \tau = \tau_1 + \tau_2$, \eqref{normalizeu} and
\eqref{normalizev}, we can rewrite the left hand side of
\eqref{duality-1} as
\begin{eqnarray}
\label{duality1}&&\sum_{j_{1},j_{2}\geq 0}
\sum_{j\geq 0}\sum_{m\geq 0}2^{-j/2}
\int_{A^*} g_j ( \xi , \mu , \tau) \chi_1 ( \xi_1 , \mu) \chi_j ( \tau -
\omega (\xi , \mu))
 \theta_{m}(\x) |\x_{1}+\x_{2}|\\
\nonumber&&\max(1,
 |\x_{1}+\x_{2}|)^{1-\epsilon_{0}}\frac{\phi_{1,j_{1}}(\x_{1},\m_{1},\tau_{1})}
{\max(1, |\x_{1}|,\gia)^{1-\epsilon_{0}}}
 \frac{\phi_{2,j_{2}}(\x_{2},\m_{2},\tau_{2})}
 {\max(1, |\x_{2}|,\gib)^{1-\epsilon_{0}}} d\xi_1 d\xi_2 d\mu_1 d\mu_2
 d\tau_1 d\tau_2
\end{eqnarray}
and the left hand side of \eqref{duality0} as
\begin{eqnarray}
\label{duality2}&&\sum_{j_{1},j_{2}\geq 0}
\sum_{j\geq 0}\sum_{n\geq 0}2^{-j/2}
\int_{A^*} g_j ( \xi , \mu , \tau) \chi_1 ( \xi_1 , \mu) \chi_j ( \tau -
\omega (\xi , \mu)) 
 \theta_{n}(\m) |\x_{1}+\x_{2}|\\
 \nonumber&&\max(1,
 |\m_{1}+\m_{2}|/|\x_{1}+\x_{2}|)^{1-\epsilon_{0}}
\frac{\phi_{1,j_{1}}(\x_{1},\m_{1},\tau_{1})}
{\max(1, |\x_{1}|,\gia)^{1-\epsilon_{0}}}
 \frac{\phi_{2,j_{2}}(\x_{2},\m_{2},\tau_{2})}
 {\max(1, |\x_{2}|,\gib)^{1-\epsilon_{0}}} d\xi_1 d\xi_2 d\mu_1 d\mu_2 .
\end{eqnarray}
Below, case A will alway correspond to the estimate \eqref{duality1},
while case B will correspond to the estimate for
\eqref{duality2}. Note that in what follows, $\chi_1$ and $\chi_2$
will always denote the functions introduced in Definition 1, and
$\chi_j$ the ones defined in Definition 2.
\noindent \newline
{\bf Region $A_{1}$}. \\
{\bf Case A:} $|\x|\geq \half|\m|/|\x|$.\\
Note that here $|\x_{1}+\x_{2}|\leq 2$ so the sum in $m$ is finite
and we can simply use the Strichartz inequality \eqref{st1} and the
fact that $l^{1}\subset l^{2}$, to obtain
$$\eqref{duality1}\lesssim \sup_{m,j}\|g_{j}\chi_{1}\chi_{j}
 \theta_{m}\|_{L^{2}}\sum_{j_{1},j_{2}\geq 0}\Pi_{i=1,2}2^{j_{i}/2}
 \|\phi_{i,j_{i}}\|_{L^{2}}$$
and the theorem follows in this case.

\noindent
{\bf Case B:} $|\x|< \half |\m|/|\x|$. \\
If $\gi\leq 1$,
it follows that $|\x|\leq \half$ and hence $|\m|\leq \half$, that is
the sum on $n$ in \eqref{duality2} reduces
to a finite sum and we proceed as above using the Strichartz
inequality \eqref{st1}.
So we can assume that  $\gi>1$.
In this case \eqref{duality2} reduces to
\begin{eqnarray}
\label{a1mstep1}
&&\sum_{j_{1},j_{2}\geq 0}\sum_{j\geq 0}\sum_{n\geq 0}2^{-j/2}
\int_{A_{1}}g_{j}\chi_{2}\chi_{j}(\x,\m,\tau)
 \theta_{n}(\m) |\x_{1}+\x_{2}|^{\epsilon_{0}}
 |\m_{1}+\m_{2}|^{1-\epsilon_{0}}\\
\nonumber
&&\frac{\phi_{1,j_{1}}(\x_{1},\m_{1},\tau_{1})}
{\max(1, |\x_{1}|,\gia)^{1-\epsilon_{0}}}
 \frac{\phi_{2,j_{2}}(\x_{2},\m_{2},\tau_{2})}
 {\max(1, |\x_{2}|,\gib)^{1-\epsilon_{0}}}.
\end{eqnarray}
We then consider two subcases.

\noindent
{\bf Case B1:} $|\m_{1}|\leq |\m_{2}|$. \\
If $\half\gib\leq |\x_{2}|$ then
$|\m_{2}|\leq 2, |\m_{1}+\m_{2}|\leq 4$. Again the sum on $n$
reduces to a finite sum and because still $|\x_{1}+\x_{2}|\leq 1$ we
go back to the previous case. If $\half \gib> |\x_{2}|$ we introduce a
dyadic decomposition with respect to $\m_{2}$ and  we set
$|\m_{2}|\sim 2^{n_{2}}$. Then
$|\m_{1}+\m_{2}|\leq C|\m_{2}|
\leq C2^{n_{2}}$ and we can write $1+ |\m_{1}+\m_{2}|\sim
2^{n_{2}+1-r}, \, 0\leq r\leq n_{2}$. We can bound  \eqref{a1mstep1} with
$$\sum_{j_{1},j_{2}\geq 0}\sum_{j\geq 0}\sum_{n_{2}\geq 0}\sum_{0\leq r
\leq n_{2}}
2^{-j/2}\int_{A_{1}}g_{j}\chi_{2}\chi_{j}(\x,\m,\tau)
 \theta_{n_{2}+1-r}2^{(n_{2}+1-r)(1-\epsilon_{0})}\phi_{1,j_{1}}(\x_
{1},\m_{1},\tau_{1})
 \frac{\chi_{2}\phi_{2,j_{2},n_{2}}(\x_{2},\m_{2},\tau_{2})}
 {(\gib)^{1-\epsilon_{0}}}.
$$
We now use the fact that $|\x_{2}|\leq 2$ and again the
Strichartz inequality  \eqref{st1} to continue  with
$$
\sum_{j,j_{1},j_{2}\geq 0}\sum_{n_{2}\geq 0}\sum_{0\leq r\leq n_{2}}
2^{-j/2}\|g_{j}\chi_{2}\chi_{j}
 \theta_{n_{2}+1-r}\|_{L^{2}}2^{-r(1-\epsilon_{0})}2^{j_{1}/2}
 \|\phi_{1,j_{1}}\|_{L^{2}}
 2^{j_{2}/2}\|\chi_{2}\phi_{2,j_{2},n_{2}}\|_{L^{2}}
$$
and this is enough to prove the theorem in this case. \\
{\bf Case B2:} $|\m_{1}|\geq  |\m_{2}|$. \\
This case can be treated like Case B1 by replacing the role of $(\x_{2},\m_{2})$
by $(\x_{1},\m_{1})$.

\noindent
{\bf Region $A_{2}$}.

\noindent
{\bf Case A:} $|\x|\geq \half|\m|/|\x|$. \\
Here \eqref{duality1} becomes
\begin{eqnarray}
\label{a2x}
&&\sum_{j_{1},j_{2}\geq 0}\sum_{j\geq 0}\sum_{m\geq 0}2^{-j/2}
\int_{A_{2}}g_{j}\chi_{1}\chi_{j}(\x,\m,\tau)
 \theta_{m}(\x) |\x_{1}+\x_{2}|^{2-\epsilon_{0}}\\
\nonumber
&&\frac{\phi_{1,j_{1}}(\x_{1},\m_{1},\tau_{1})}
{\max(1, |\x_{1}|,\gia)^{1-\epsilon_{0}}}
 \frac{\phi_{2,j_{2}}(\x_{2},\m_{2},\tau_{2})}
 {\max(1, |\x_{2}|,\gib)^{1-\epsilon_{0}}}.
\end{eqnarray}
We dyadically decompose $|\x_{1}|\sim|\x_{2}|\sim 2^{m_{1}}$. Then
$1+ |\x|\sim 2^{m}$, with $m=m_{1}+1-r, 0\leq r\leq m_{1}$ and
\begin{eqnarray*}
\eqref{a2x}&\lesssim &\sum_{j_{1},j_{2}\geq 0}
\sum_{j\geq 0}\sum_{m_1 \geq 0}\sum_{0\leq r \leq  m_{1}}
2^{-j/2}\int_{A_{2}}g_{j}\chi_{1}\chi_{j}(\x,\m,\tau)
 \theta_{m_{1}+1-r}(\x) |\x_{1}+\x_{2}|^{2-\epsilon_{0}}\\
\nonumber
&&\frac{\phi_{1,j_{1},m_{1}}(\x_{1},\m_{1},\tau_{1})}
{2^{m_{1}(1-\epsilon_{0})}}
 \frac{\phi_{2,j_{2},m_{1}}(\x_{2},\m_{2},\tau_{2})}
 {2^{m_{1}(1-\epsilon_{0})}}\\
 \nonumber &\lesssim&\sum_{j_{1},j_{2}\geq 0}
 \sum_{j\geq 0}\sum_{m_1 \geq 0}
 \sum_{0\leq r \leq m_{1}}
2^{-j/2}\int_{A_{2}}g_{j}\chi_{1}\chi_{j}(\x,\m,\tau)
 \theta_{m_{1}+1-r}(\x) 2^{-r(2-2\epsilon_{0})}\\
 & &|\x_{1}+\x_{2}|^{\epsilon_{0}}\Pi_{i=1,2}\phi_{i,j_{1},m_{i}}
 (\x_{i},\m_{i},\tau_{i}).
 \end{eqnarray*}
We  apply  Lemma \ref{mainlemma} relative to the region $A_{2}$ and
we continue with
\begin{eqnarray*}
&\lesssim&\sup_{j,m}\|g_{j}\chi_{1}\chi_{j}\theta_{m}\|_{L^{2}}
\sum_{m_1 \geq 0}\sum_{0\leq r \leq m_{1}}2^{-r(2-2\epsilon_{0})}
\Pi_{i=1,2}\sum_{j_{i}\geq 0}2^{j_{i}/2}\|\phi_{i,j_{i},m_{1}}\|_{L^{2}}\\
&\lesssim&\sum_{j_{1},j_{2}\geq 0}2^{j_{1}/2}2^{j_{2}/2}
\left(\sum_{m_{1}\geq 0}\|\phi_{1,j_{1},m_{1}}\|_{L^{2}}
\|\phi_{2,j_{2},m_{1}}\|_{L^{2}}\right) \\
&\lesssim& \sum_{j_1 , j_2 \geq 0} 2^{j_1 /2 } 2^{j_2 /2} \left(
  \sum_{m_1 \geq 0} {{\| \phi_{1, j_1 , m_1} \|}^2_{L^2}}
\right)^\half
\left(
  \sum_{m_1 \geq 0} {{\| \phi_{2, j_2 , m_1} \|}^2_{L^2}}
\right)^\half \\
&\lesssim &
\left(
  \sum_{j_1 \geq 0} 2^{j_1 /2} {{\| \phi_{1, j_1 } \|}^2_{L^2}}
\right)^\half
\left(
  \sum_{j_2 \geq 0} 2^{j_2 /2} {{\| \phi_{2, j_2 } \|}^2_{L^2}}
\right)^\half,
\end{eqnarray*}
which concludes the argument.
%
%% then the inclusion  $l^{1}\subset l^{2}$ and Cauchy-Schwarz in $m_{1}$
%% conclude the argument.

\noindent
{\bf Case B:} $|\x|< \half|\m|/|\x|$.\\
If $\gi\leq 1$
it follows that
$|\x|\leq \half$ and hence $|\m|\leq \half$ and we go back to the same
estimates presented for region $A_{1}$. So we can assume
$\gi\geq 1$. We have to estimate \eqref{a1mstep1} where now the integral
takes place in $A_{2}$. Again we consider two subcases.\\
{\bf Case B1:} $|\m_{1}|\leq |\m_{2}|$. \\
If $\half\gib\leq |\x_{2}|$ then
$|\m_{2}|\leq 2|\x_{2}|^{2}, |\m_{1}+\m_{2}|\leq 4|\x_{2}|^{2}$.
We dyadically decompose with respect to
$|\x_{1}|\sim |\x_{2}|\sim 2^{m_{2}}$, so that
$1+ |\m_{1}+\m_{2}|\sim 2^{-r+2m_{2}}, 0\leq r\leq 2m_{2}.$ Then in this
case we bound  \eqref{duality2} with
\begin{eqnarray*}
&&\sum_{j_{1},j_{2}\geq 0}
\sum_{j\geq 0}\sum_{m_2 \geq 0}
 \sum_{0\leq r \leq  2m_{2}}
2^{-j/2}\int_{A_{2}}g_{j}\chi_{2}\chi_{j}(\x,\m,\tau)
 \theta_{2m_{2}+1-r}(\m) 2^{-r(1-\epsilon_{0})}\\
& & |\x_{1}+\x_{2}|^{\epsilon_{0}}\Pi_{i=1,2}\phi_{i,j_{i},m_{2}}
 (\x_{i},\m_{i},\tau_{i}).
 \end{eqnarray*}
and one can use again Lemma \ref{mainlemma} in region $A_{2}$ as above.
If $\half\gib> |\x_{2}|$ we introduce a
dyadic decomposition with respect to $\m_{2}$ and  we set
$|\m_{2}|\sim 2^{n_{2}}$. Then
$|\m_{1}+\m_{2}|\leq C|\m_{2}|
\leq C2^{n_{2}}$ and we can write $1+|\m_{1}+\m_{2}|\sim 2^{n_{2}+1-r},
\, \,
0\leq r\leq n_{2}$. We then  bound  \eqref{duality2} with
\begin{eqnarray*}
 & &  \sum_{j_{1},j_{2}\geq 0}
 \sum_{j\geq 0}\sum_{n_{2}\geq 0}\sum_{0\leq r\leq n_{2}}
2^{-j/2}\int_{A_{2}}g_{j}\chi_{2}\chi_{j}(\x,\m,\tau)
\theta_{n_{2}+1-r}2^{(n_{2}+1-r)(1-\epsilon_{0})}\\
& &\times|\x_{1}+\x_{2}|^{\epsilon_{0}}
 \frac{\phi_{1,j_{1}}(\x_{1},\m_{1},\tau_{1})}{|\x_{1}|^{1-\epsilon_{0}}}
 \frac{\chi_{2}\phi_{2,j_{2},n_{2}}(\x_{2},\m_{2},\tau_{2})}
 {(\gib)^{1-\epsilon_{0}}}\\
 &\lesssim&\sum_{j_{1},j_{2}\geq 0}
 \sum_{n_{2}\geq 0}\sum_{0\leq r\leq n_{2}}\sum_{j\geq 0}
 2^{-j/2}\int_{A_{2}}g_{j}\chi_{2}\chi_{j}(\x,\m,\tau)
\theta_{n_{2}+1-r}2^{-r(1-\epsilon_{0})}\\
& &\times
|\x_{1}+\x_{2}|^{\epsilon_{0}}
 \phi_{1,j_{1}}(\x_{1},\m_{1},\tau_{1})
 \chi_{2}\phi_{2,j_{2},n_{2}}(\x_{2},\m_{2},\tau_{2})
\end{eqnarray*}
We use again Lemma \ref{mainlemma} and we continue with
$$
\sup_{j,n}\|g_{j}\chi_{2}\chi_{j}
 \theta_{n}\|_{L^{2}}\sum_{n_{2}\geq 0}\sum_{0\leq r\leq n_{2}}
2^{-r(1-\epsilon_{0})}\left(\sum_{j_{1}\geq 0}2^{j_{1}/2}
\|\phi_{1,j_{1}}\|_{L^{2}}\right)
 \left(\sum_{j_{2}\geq 0}2^{j_{2}/2}
 \|\chi_{2}\phi_{2,n_{2},j_{2}}\|_{L^{2}}\right)
$$
and this is enough to prove the theorem in this case. \\
{\bf Case B2:} $|\m_{1}|\geq |\m_{2}|$. One can use the same argument
presented for Case B1 inverting the role of $(\x_{1},\m_{1})$ and
$(\x_{2},\m_{2})$.

\noindent
{\bf Region $A_{3}\cup A_{4} \cup (A_{5}-\tilde{A}_{5}(\epsilon_{0}))$}.\\
{\bf Case A:} $ \half\gi\leq|\x|$. \\ We consider two subcases.\\
{\bf Case A1:} $|\xi_{1}|\geq \half\gia$.
We dyadically decompose with respect to $|\x|\sim|\x_{1}|\sim 2^{m_{1}}$
and we bound \eqref{duality1}, now integrated over the region
$A=A_{3}\cup A_{4} \cup (A_{5}-\tilde{A}_{5} (\epsilon_0))$,  with
\begin{eqnarray*}
&&\sum_{j_{1},j_{2}\geq 0}\sum_{j\geq 0}\sum_{m_{1}\geq 0}2^{-j/2}
\int_{A}g_{j}\chi_{1}\chi_{j}(\x,\m,\tau)
 \theta_{m_{1}}(\x)
 |\x_{1}+\x_{2}|^{2-\epsilon_{0}}\\
&&\frac{\chi_{2}\phi_{1,j_{1},m_{1}}}{|\x_{1}|^{1-\epsilon_{0}}}
(\x_{1},\m_{1},\tau_{1})
\frac{\phi_{2,j_{2}}}{(\gib)^{1-\epsilon_{0}}}(\x_{2},\m_{2},\tau_{2}).
\end{eqnarray*}
But in this region $\gib\gtrsim |\xi_{1}|\sim |\x_{1}+\x_{2}|$, hence
we can continue with
\begin{equation}
\sum_{j_{1},j_{2}\geq 0}\sum_{j\geq 0}\sum_{m_{1}\geq 0}2^{-j/2}
\int_{A}g_{j}\chi_{1}\chi_{j}(\x,\m,\tau)
 \theta_{m_{1}}(\x)
 |\x_{1}+\x_{2}|^{\epsilon_{0}}\chi_{1}\phi_{1,j_{1},m_{1}}
\phi_{2,j_{2}}
\label{helpy1}\end{equation}
We then apply Lemma \ref{mainlemma} and we obtain the desired
result.\\
{\bf Case A2:} $|\xi_{1}|\leq \half\gia$. \\
We do a dyadic decomposition of $\phi_{2}$ in $(\m_{2}/\x_{2})$, so
that we write $\phi_{2}=\sum_{r_{2} \geq 0}\tilde\phi_{2,r_{2}}= \sum_{r_2 \geq 0} \sum_{j_2 \geq 0} 
\widetilde{\phi}_{2, r_2, j_2} $, where
$1+ |\m_{2}|/|\x_{2}|\sim 2^{r_{2}}$ and $1 + | \tau_2 - \omega( \xi_2 , \mu_2 ) | \thicksim 2^{j_2}$. 
Note that
$|\x_{1}|\leq C\gia \leq C\gib$, and that if $|\x_{1}|\sim 2^{m_{1}}$,
then $2^{r_{2}}\sim 2^{m_{1}+r}, \, r\leq -C$. We bound
 \eqref{duality1} with
\begin{eqnarray}
\nonumber&&\sum_{j_{1},j_{2}\geq 0}\sum_{j\geq 0}\sum_{m_{1}\geq 0}
\sum_{r\geq -C}2^{-j/2}
\int_{A}g_{j}\chi_{1}\chi_{j}(\x,\m,\tau)
 \theta_{m_{1}}(\x)\\
\nonumber & &|\x_{1}+\x_{2}|^{2-\epsilon_{0}}
\frac{\chi_{2}\phi_{1,j_{1},m_{1}}}{|\x_{1}|^{1-\epsilon_{0}}}
(\x_{1},\m_{1},\tau_{1})
\frac{\tilde\phi_{2,m_{1}+r,j_{2}}}
{(\gib)^{1-\epsilon_{0}}}(\x_{2},\m_{2},\tau_{2})\\
\label{helpy2}&\lesssim&\sum_{j_{1},j_{2}\geq 0}
\sum_{j\geq 0}\sum_{m_{1}\geq 0}\sum_{-C\leq
r}2^{-j/2}2^{-r(1-\epsilon_{0})}
\int_{A}g_{j}\chi_{1}\chi_{j}(\x,\m,\tau)
 \theta_{m_{1}}(\x) \\
\nonumber & & |\x_{1}+\x_{2}|^{\epsilon_{0}}\chi_{2}\phi_{1,j_{1},m_{1}}
\tilde\phi_{2,m_{1}+r,j_{2}}.
\end{eqnarray}
If we use again Lemma \ref{mainlemma} we can continue the chain of
inequalities with
\begin{equation}
\sup_{j,m_1 }\|g_{j}\chi_{1}\chi_{j}\theta_{m_{1}}\|_{L^{2}}
\sum_{m_{1}\geq 0}\sum_{-C \leq r}2^{-r(1-\epsilon_{0})}
\sum_{j_{1}\geq 0}2^{j_{1}/2}\|\chi_{2}\phi_{1,j_{1},m_{1}}\|_{L^{2}}
\sum_{j_{2}\geq 0}2^{j_{2}/2}\|\tilde\phi_{2,m_{1}+r,j_{2}}\|_{L^{2}}.
\label{helpy3}\end{equation}
Now Cauchy-Schwarz in $m_{1}$ is enough to obtain Theorem
\ref{xbil} in this case.\\
{\bf Case B:} $ \half\gi>|\x|$. As we observed in the analysis
of region $A_{i}, i=1,2$, without loss of generality we can assume that
$\gi\geq 1$. We consider two subcases.\\
{\bf Case B1:}  $|\m_{1}|\leq |\m_{2}|$. \\
We recall that in this
region we also have $\gib>>|\x_{2}|$. We repeat the argument presented
in the second part of Case B1 of region $A_{2}$.\\
{\bf Case B2:} $|\m_{1}|\geq |\m_{2}|$.\\
If $\half\gia\leq|\x_{1}|$ then
we use an  argument similar to the one  in the first part of
Case B1 in region
$A_{2}$, where the role of $(\x_{2},\m_{2})$ is now played by
$(\x_{1},\m_{1})$. 
% In particular, if $|\xi_{1}|\sim 2^{m_{1}}$, we can
% bound \eqref{a1mstep1}  with
% %
% $$
% \sum_{j_{1},j_{2}\geq 0}\sum_{j\geq 0}\sum_{m_1 \geq 0}
%  \sum_{0\leq r \leq  2m_{1}}
% 2^{-j/2}\int_{A}g_{j}\chi_{2}\chi_{j}(\x,\m,\tau)
%  \theta_{2m_{1}+1-r}(\m) 2^{-r(1-\epsilon_{0})}
% |\x_{1}+\x_{2}|^{\epsilon_{0}}\chi_{1}\phi_{1,j_{1},m_{1}}\phi_{2,j_{2}}.
% $$
% %
In particular, if $|\xi_1 | \thicksim 2^{m_1}$, since $|\mu_2 | / |\xi_2 | \gtrsim |\xi_1 |$ and $|\xi_1 + \xi_3| \thicksim |\xi_1|$, we can bound  \eqref{a1mstep1} with
\begin{eqnarray*}
& &  \sum_{j_1, j_2 \geq 0} \sum_{j \geq 0} \sum_{m_1 \geq 0} 2^{-j/2} \int_A g_j \chi_s \chi_j (\xi, \mu, \tau) |\xi_1 + \xi_2 |^{\epsilon_0} \\
& \times & \frac{|\mu_1 + \mu_2 |^{1 - \epsilon_0}}{2^{2(1-\epsilon_0) m_1}} \theta_{m_1} ( \xi_1)
\phi_{1,j_1} (\xi_1, \mu_1, \tau_1) \phi_{2,j_2} (\xi_2, \mu_2, \tau_2). 
\end{eqnarray*}
But, $|\mu_1 + \mu_2 | \leq 2 |\mu_1 | \lesssim \frac{|\mu_1|}{|\xi_1|} |\xi_1 | \lesssim |\xi_1|^2$,
so that 
\begin{equation*}
  (1 + |\mu_1 + \mu_2 |) \thicksim 2^{2 m_1 -r + 1}, 0 \leq r \leq 2m_1 ,
\end{equation*}
so we can continue with
\begin{eqnarray*}
& &  \sum_{j_1, j_2 \geq 0} \sum_{m_1 \geq 0} \sum_{0 \leq r \leq 2 m_1}  2^{-j/2}
\int_A g_j \chi_2 \chi_j (\xi, \mu, \tau) \theta_{2m_1 + 1 -r } (\mu ) \\
&  \times& 2^{-r (1-\epsilon_0 )} |\xi_1 + \xi_2 |^{\epsilon_0} \chi_1 \phi_{1,j_1, m_1} \phi_{2, j_2}.
\end{eqnarray*}
We then use Lemma \ref{mainlemma}.
If $\half\gia \geq |\x_{1}|$, we dyadically decompose so that
$|\m_{1}|\sim 2^{n_{1}}$. Then $|\m_{1}+\m_{2}|\sim 2^{n_{1}+1-r},
0\leq r\leq n_{1}$. We bound \eqref{a1mstep1} with
\begin{eqnarray*}
 & &   \sum_{j\geq 1}\sum_{n_{1}\geq 0}\sum_{0\leq r\leq n_{1}}
2^{-j/2}\int_{A}g_{j}\chi_{2}\chi_{j}(\x,\m,\tau)
\theta_{n_{1}+1-r}2^{(n_{1}+1-r)(1-\epsilon_{0})}\\
& &\times|\x_{1}+\x_{2}|^{\epsilon_{0}}
 \frac{\chi_{2}\phi_{1,n_{1}}(\x_{1},\m_{1},\tau_{1})}
{(\gia)^{1-\epsilon_{0}}}
 \frac{\phi_{2}(\x_{2},\m_{2},\tau_{2})}
 {(\gib)^{1-\epsilon_{0}}}\\
 &\lesssim&\sum_{n_{1}\geq 0}\sum_{0\leq r\leq n_{1}}\sum_{j\geq 1}
 2^{-j/2}\int_{A}g_{j}\chi_{2}\chi_{j}(\x,\m,\tau)
\theta_{n_{1}+1-r}2^{-r(1-\epsilon_{0})}\\
& &\times
|\x_{1}+\x_{2}|^{\epsilon_{0}}2^{n_{1}(1-\epsilon_{0})}
\frac{\chi_{2}\phi_{1,n_{1}}(\x_{1},\m_{1},\tau_{1})}
{(\gia)^{1-\epsilon_{0}}}
 \frac{\phi_{2}(\x_{2},\m_{2},\tau_{2})}
 {|\x_{1}|^{1-\epsilon_{0}}}\\
 &\lesssim&\sum_{n_{1}\geq 0}\sum_{0\leq r\leq n_{1}}\sum_{j\geq 1}
 2^{-j/2}\int_{A}g_{j}\chi_{2}\chi_{j}(\x,\m,\tau)
\theta_{n_{1}+1-r}2^{-r(1-\epsilon_{0})}
|\x_{1}+\x_{2}|^{\epsilon_{0}}\chi_{2}\phi_{1,n_{1}}
\phi_{2}
\end{eqnarray*}
Now again one uses Lemma \ref{mainlemma} to conclude the argument.

\noindent
{\bf Region $A_{5}\cap\tilde{A}_{5}(\epsilon_{0})$}.
We summarize the restrictions that
occur in this region: for $\alpha>4$
\begin{equation}
 \begin{array}{l}
     |\x_{2}|\leq 1, \, \, \gia\lesssim \gib, \, \,
     |\m_{1}/\x_{1}-\m_{2}/\x_{2}|\lesssim |\x_{1}+\x_{2}|\\
     |\x_{1}|^{1+\alpha\epsilon_{0}}\leq \gia,\, \, \,
     \gia\geq 100|\x_{1}|.
     \end{array}
\label{conditions}    \end{equation}
In this case, $A = A^* \cap A_5 \cap {\widetilde{A_5}} ( \epsilon_0 ).$

\noindent
{\bf Case A:} $|\x|\geq \half\gi$.\\
We dyadically decompose with respect to $\x$, so that
$|\x|\sim|\xi_{1}|\sim 2^{m}$. We write
\begin{eqnarray*}
& &\sum_{j_{1},j_{2}\geq 0}\sum_{j\geq 0}\sum_{m\geq 0}2^{-j/2}
\int_{A}g_{j}\chi_{1}\chi_{j}(\x,\m,\tau)
 \theta_{m}(\x)
 |\x_{1}+\x_{2}|^{2-\epsilon_{0}}\\
 & &\times
\frac{\phi_{1,j_{1},m}}{(\gia)^{1-\epsilon_{0}}}
(\x_{1},\m_{1},\tau_{1})
\frac{\phi_{2,j_{2}}}
{(\gib)^{1-\epsilon_{0}}}(\x_{2},\m_{2},\tau_{2}).
\end{eqnarray*}
Now let's consider the multiplier in the above integral. Using
\eqref{conditions} (with $\alpha>4$) we can write
\begin{equation}
    \frac{|\x_{1}|^{2-\epsilon_{0}}}{(\gia)^{2-2\epsilon_{0}}}
\lesssim \frac{|\x_{1}|^{1-2\epsilon_{0}}}{(\gia)^{1-2\epsilon_{0}}}
\lesssim |\x_{1}|^{-\alpha\epsilon_{0}(1-2\epsilon_{0})}.
\label{newmult}\end{equation}
Using \eqref{newmult} and Strichartz  inequality \eqref{st1}
we can continue the chain of
inequalities with
\begin{eqnarray*}
&\lesssim&\sum_{j\geq 0}2^{-j/2}\sum_{m\geq 0}
\|g_{j}\theta_{m}\chi_{1}\chi_{j}\|_{L^{2}}
 2^{-\alpha\epsilon_{0}(1-2\epsilon_{0})m}
 \sum_{j_{1},j_{2}\geq 0}2^{j_{1}/2}\|\phi_{1,m,j_{1}}\|_{L^{2}}
 2^{j_{2}/2} \|\phi_{2,j_{2}}\|_{L^{2}}
\end{eqnarray*}
and Cauchy-Schwarz in $m$ is enough to prove the theorem in this case.

\noindent
{\bf Case B:} $|\x|<\half\gi$.\\
From \eqref{conditions} we have that
\begin{equation}
    \gia\sim\gib,
\label{m1x1m2x2}\end{equation}
hence
$|\m_{2}|\sim |\x_{2}|\gia<<|\m_{1}|$, hence
$|\m_{1}+\m_{2}|\sim |\m_{1}|$. We dyadically decompose with respect to
$|\m_{1}|\sim 2^{n_{1}}$. We have to estimate
\begin{eqnarray*}
& &\sum_{j_{1},j_{2}\geq 0}\sum_{j\geq 0}\sum_{n_{1}\geq 0}2^{-j/2}
\int_{A}g_{j}\chi_{2}\chi_{j}(\x,\m,\tau)
 \theta_{n_{1}}(\m)
 |\x_{1}+\x_{2}|^{\epsilon_{0}}|\m_{1}+\m_{2}|^{1-\epsilon_{0}}\\
 & &\times
\frac{\chi_{2}\phi_{1,n_{1},j_{1}}}{(\gia)^{1-\epsilon_{0}}}
(\x_{1},\m_{1},\tau_{1})
\frac{\phi_{2,j_{2}}}
{(\gib)^{1-\epsilon_{0}}}(\x_{2},\m_{2},\tau_{2}).
\end{eqnarray*}
Now let's consider the multiplier in this  integral. Using
\eqref{conditions} we can write
$$
    \frac{|\x_{1}|(\gia)^{1-\epsilon_{0}}}{(\gia)^{2-2\epsilon_{0}}}
\lesssim \frac{(\gia)^{1/(1+\alpha\epsilon_{0})}}{(\gia)^{1-\epsilon_{0}}}
\lesssim 1,
$$
provided $\epsilon_{0}<(\alpha-1)/\alpha$. Then using Strichartz we can
continue the chain of inequality with
\begin{eqnarray*}
&\lesssim&\sum_{j\geq 0}2^{-j/2}\sum_{n_{1}\geq 0}
\|g_{j}\theta_{n_{1}}\chi_{2}\chi_{j}\|_{L^{2}}
 \sum_{j_{1},j_{2}\geq 0}2^{j_{1}/2}\|\chi_{2}
 \phi_{1,n_{1},j_{1}}\|_{L^{2}}
 2^{j_{2}/2} \|\phi_{2,j_{2}}\|_{L^{2}},
\end{eqnarray*}
and the theorem is proved also in this case.

\noindent
{\bf Region $A_{6}$}. Of the whole theorem this is the region in which
the estimates are the most delicate.
We summarize the restrictions on this region:
\begin{equation}
 \begin{array}{l}
     |\x_{1}|\geq 1, \, \, \, |\x_{2}|\leq  10^{-10}|\x_{1}|\\
     \gib\leq 10^{-10}\max(|\x_{1}|, \gia).
     \end{array}
\label{conditionsa6}    \end{equation}
We observe that the multipliers appearing in
\eqref{duality1} and \eqref{duality2} can be bounded in the following
way:
\begin{equation}
 \frac{|\x|\max(1,|\x|,\gi)^{1-\epsilon_{0}}}
 {\max(1,|\x_{1}|,\gia)^{1-\epsilon_{0}}
 \max(1,|\x_{2}|,\gib)^{1-\epsilon_{0}}}\lesssim
 \frac{|\x_{1}|}{\max(1,|\x_{2}|,\gib)^{1-\epsilon_{0}}}
    \label{a6mult}\end{equation}
This is obvious when $|\x|\geq \half\gi$ or when
$|\x|\leq  \half\gi$ and $\gi\leq 1$. In the remaining region we estimate
the numerator
\begin{eqnarray*}
|\x|\max(1,|\x|,\gi)^{1-\epsilon_{0}}&\sim& |\x_{1}|^{\epsilon_{0}}
|\m_{1}+\m_{2}|^{1-\epsilon_{0}}\\
&\lesssim&|\x_{1}|^{\epsilon_{0}}
|\m_{1}|^{1-\epsilon_{0}}+
|\x_{1}|^{\epsilon_{0}}|\m_{2}|^{1-\epsilon_{0}}\\
&\lesssim &|\x_{1}|^{\epsilon_{0}}|\x_{1}|^{1-\epsilon_{0}}
(\gia)^{1-\epsilon_{0}}+|\x_{1}|^{\epsilon_{0}}|\x_{2}|^{1-\epsilon_{0}}
(\gib)^{1-\epsilon_{0}}\\
&\lesssim&|\x_{1}|((\gia)^{1-\epsilon_{0}}+
\max(|\x_{1}|, \gia)^{1-\epsilon_{0}}).
\end{eqnarray*}
{\bf Case A:} $|\x_{1}|\geq 10^{2}\gia$.\\
Let's show that in this case we also have $|\x|\geq \gi$.
\begin{eqnarray*}
\gi&\leq &(1-10^{-10})^{-1}\frac{|\m_{1}+\m_{2}|}{|\x_{1}|}
\leq (1-10^{-10})^{-1}\left(\gia+\gib \frac{|\x_{2}|}{|\x_{1}|}\right)
\\
&\leq&(1-10^{-10})^{-1}\left(\gia+10^{-20}|\x_{1}|\right)\leq
(1-10^{-10})^{-1}(10^{-2}+10^{-20})|\xi_{1}|\\
&\leq&
(1-10^{-10})^{-2}(10^{-2}+10^{-20})|\x_{1}+\x_{2}|\leq |\x|.
\end{eqnarray*}
We have to estimate \eqref{duality1} and we use again
the functions $\phi_{i,j_{i}}$. We change variables
in $\tau_{1}$ and $\tau_{2}$ as in
\eqref{a2step2} and we use \eqref{a6mult} to bound
\eqref{duality1} with
\begin{eqnarray}
 \label{a6step1}&&
\sum_{j\geq 0} \sum_{m\geq 0} \sum_{j_1 , j_2 \geq 0}
2^{-j/2}
 \int \chi_{1}\theta_{m_{1}}g_{j}(\x_{1}+\x_{2},\m_{1}+\m_{2},\theta_{1}+
 \om(\x_{1},\m_{1})+\theta_{2}+\om(\x_{2},\m_{2}))\\
 \nonumber & &\chi_{j}(\theta_{1}+
 \om(\x_{1},\m_{1})+\theta_{2}+\om(\x_{2},\m_{2})
 -\om(\x_{1}+\x_{2},\m_{1}+\m_{2}))\chi_{1}(\xi_{1}, \mu_1 )\\
 \nonumber& &\frac{|\x_{1}|}{\max(1,|\x_{2}|, \gib)^{1-\epsilon_{0}}}
 \Pi_{i=1,2}\phi_{i}(\x_{i},\m_{i},\theta_{i}+\om(\x_{i},\m_{i}))
 \chi_{j_{i}}(\theta_{i})d\x_{i}d\m_{i}d\tau_{i}.
 \end{eqnarray}
We change variables again and this time we use \eqref{chanm}. From
\eqref{jm} we deduce that $|J_{\m}|\geq C|\x_{1}|^{2}$. We also
perform a dyadic decompositions by setting $|\x_{i}|\sim 2^{m_{i}}$
(hence $|\x|\sim 2^{m_{1}}$)
and  $|\mu_{2}|\sim 2^{n_{2}}$. Let $m^{*}_{2}=\max(n_{2}-m_{2},m_{2})$,
(here $m_{2}, n_{2}\in \dbZ$).
\\
{\bf Case A1:} $j\geq \max(0, m^{*}_{2}).$\\
Denote with $H(u,v,w,\m_{2},\theta_{1},\theta_{2})$ the transformation
of
$$\chi_{1}(\xi_{1}, \mu_1)\phi_{1}(\x_{1},\m_{1},\theta_{1}+\om(\x_{1},\m_{1}))
 \chi_{j_{1}}(\theta_{1})\frac{\phi_{2}
 (\x_{2},\m_{2},\theta_{2}+\om(\x_{2},\m_{2}))
 \chi_{j_{2}}(\theta_{2})}{\max(1,|\x_{2}|,\gib)^{1-\epsilon_{0}}}$$
under the  change of variables  in \eqref{chanm}. Here we use
the fact that $\Delta_{\m_{2}}\sim 2^{n_{2}}$.
Then we can rewrite \eqref{a6step1} as
\begin{eqnarray*}
 & &\sum_{m_{1}\geq 0, m_{2}}\sum_{n_{2}}
 \sum_{j_{1},j_{2}\geq 0}\sum_{j\geq \max(0, m^{*}_{2})}
 2^{-j/2}2^{m_{1}}\int g_{j}\theta_{m_{1}}\chi_{1}\chi_{j}(u,v,w)\\
 & &|J_{\m}|^{-1}H(u,v,w,\m_{2},\theta_{1},\theta_{2})dudvdw
 d\m_{2}d\theta_{1}d\theta_{2}\\
 &\lesssim &\sum_{m_{1}\geq 0, m_{2}}\sum_{n_{2}}
 \sum_{j_{1},j_{2}\geq 0}\sum_{j\geq \max(0, m^{*}_{2})}
 2^{-j/2}2^{m_{1}}2^{n_{2}/2}\\
 &&\int g_{j}\theta_{m_{1}}\chi_{1} \chi_{j}(u,v,w)
 \left(\int_{\m_{2}}|J_{\m}|^{-2}
 H^{2}(u,v,w,\m_{2},\theta_{1},\theta_{2})d\m_{2}\right)^{\half}dudvdw
 d\theta_{1}d\theta_{2}.
 \end{eqnarray*}
 Now we observe that by Cauchy-Schwarz and the inverse change of
 variable we have
\begin{eqnarray*}
& &\int g_{j}\theta_{m_{1}}\chi_{1} \chi_{j}(u,v,w)
\left(\int_{\m_{2}}|J_{\m}|^{-2}
 H^{2}(u,v,w,\m_{2},\theta_{1},\theta_{2})d\m_{2}
\right)^{\half}dudvdw
 d\theta_{1}d\theta_{2}\\
&\lesssim& 2^{-m_{1}}\|g_{j}\chi_{j} \theta_{m_{1}}
\chi_{1} \|_{L^{2}}
 2^{j_{1}/2}\|\chi_{1}\phi_{1,j_{1},m_{1}}\|_{L^{2}}
2^{j_{2}/2}\frac{\|\phi_{2,j_2, m_{2},n_{2}}\|_{L^{2}}}
{\max(1,2^{m^{*}_{2}})^{1-\epsilon_{0}}},
 \end{eqnarray*}
where $\phi_{2, j_2, m_2, n_2} = \phi_{2, j_2 , m_2 } \theta_{n_2},$
and inserting this above we can continue with
\begin{eqnarray}
\label{a6step2}&\lesssim &\sum_{m_{1}
\geq 0, m_{2}}\sum_{n_{2}}
 \sum_{j_{1},j_{2}\geq 0}\sum_{j\geq \max(0, m^{*}_{2})}
 2^{-j/2}2^{n_{2}/2}\\
 \nonumber& &\|g_{j}\chi_{j} \theta_{m_{1}}\chi_{1} \|_{L^{2}}
 2^{j_{1}/2}
 \|\chi_{1}\phi_{1,j_{1},m_{1}}\|_{L^{2}}
2^{j_{2}/2}\frac{\|\phi_{2,j_2 m_{2},n_{2}}\|_{L^{2}}}
{\max(1,2^{m^{*}_{2}})^{1-\epsilon_{0}}}.
\end{eqnarray}
{\bf Case A1a:} $m_{2}\leq 0, \, n_{2}-m_{2}>0$.\\
Now  $m^{*}_{2}=n_{2}-m_{2}$ and \eqref{a6step2} can be bounded by
\begin{eqnarray*}
& &\sum_{m_{1}
\geq 0, m_{2}\leq 0}\sum_{n_{2}\geq m_{2}}
 \sum_{j_{1},j_{2}\geq 0}\sum_{j\geq \max(0, m^{*}_{2})}
 2^{-j/2}2^{n_{2}/2}2^{(-n_{2}+m_{2})(1-\epsilon_{0})}\\
 & &\|g_{j}\chi_{j} \theta_{m_{1}}\chi_{1} \|_{L^{2}}
 2^{j_{1}/2}
 \|\chi_{1}\phi_{1,j_{1},m_{1}}\|_{L^{2}}
2^{j_{2}/2}\|\phi_{2,j_2 , m_{2},n_{2}}\|_{L^{2}},
\end{eqnarray*}
and if we sum for $j\geq n_{2}-m_{2}$ we can continue our chain of
inequalities with
$$\lesssim \sum_{m_{1}
\geq 0}\sum_{n_{2}\geq m_{2}, m_{2}\leq 0}
 \sum_{j_{1},j_{2}\geq 0}
2^{m_{2}/2}2^{(-n_{2}+m_{2})(1-\epsilon_{0})}
 2^{j_{1}/2}
 \|\chi_{1}\phi_{1,j_{1},m_{1}}\|_{L^{2}}
2^{j_{2}/2}\|\phi_{2,j_2 , m_{2},n_{2}}\|_{L^{2}},
$$
and by Cauchy-Schwarz on $n_{2}$ and $m_{2}\leq 0$ we obtain
\begin{equation}\eqref{duality1}\lesssim \sum_{m_{1}\geq 0}
 \sum_{j_{1},j_{2}\geq 0}2^{j_{1}/2}
 \|\chi_{1}\phi_{1,j_{1},m_{1}}\|_{L^{2}}
2^{j_{2}/2}\|\phi_{2,j_{2}}\|_{L^{2}}
\label{a6step3}\end{equation}
and this proves the theorem in this case.\\
%%%%%%%%%%%%%%%%%%%%%%
{\bf Case A1b:} $m_{2}\leq 0, \, n_{2}-m_{2}\leq 0$.\\
In this case $\max(1,2^{m^{*}_{2}})=1$. We repeat the argument above
and we bound \eqref{duality1} with
\begin{eqnarray}
& &\sum_{m_{1}
\geq 0}\sum_{n_{2}\leq m_{2}\leq 0}
 \sum_{j_{1},j_{2}\geq 0}\sum_{0\leq j}
 2^{-j/2}2^{n_{2}/2}\\
 \nonumber& &\|g_{j}\chi_{j}\theta_{m_{1}}\chi_{1} \|_{L^{2}}
 2^{j_{1}/2}
 \|\chi_{1}\phi_{1,j_{1},m_{1}}\|_{L^{2}}
2^{j_{2}/2}\|\phi_{2,j_2, m_{2},n_{2}}\|_{L^{2}}
\end{eqnarray}
and one obtains again \eqref{a6step3} after summing in $j$ and
using Cauchy-Schwarz first in $n_{2}$ and then in $m_{2}$.\\
%%%%%%%%%%%%%%%%
{\bf Case A1c:} $m_{2}>0, \, n_{2}-m_{2}\leq 0$.\\
In this case $m_{2}^{*}=m_{2}>0$. We obtain
\begin{eqnarray}
\eqref{duality1}&\lesssim &\sum_{m_{1}
\geq 0}\sum_{n_{2}\leq m_{2}}\sum_{0\leq m_{2}}
 \sum_{j_{1},j_{2}\geq 0}\sum_{j\geq m_{2}}
 2^{-j/2}2^{n_{2}/2}2^{-m_{2}(1-\epsilon_{0})}\\
 \nonumber& &\|g_{j}\chi_{j}\theta_{m_{1}}\chi_{1} \|_{L^{2}}
 2^{j_{1}/2}
 \|\chi_{1}\phi_{1,j_{1},m_{1}}\|_{L^{2}}
2^{j_{2}/2}\|\phi_{2,j_2 , m_{2},n_{2}}\|_{L^{2}}
\end{eqnarray}
and summing over $j$
$$
\lesssim \sum_{m_{1}
\geq 0}\sum_{n_{2}\leq m_{2}}\sum_{0\leq m_{2}}
 \sum_{j_{1},j_{2}\geq 0}
 2^{-m_{2}/2}2^{n_{2}/2}2^{-m_{2}(1-\epsilon_{0})}
 2^{j_{1}/2}
 \|\chi_{1}\phi_{1,j_{1},m_{1}}\|_{L^{2}}
2^{j_{2}/2}\|\phi_{2, j_2 , m_2 , n_2 }\|_{L^{2}}
$$
and by Cauchy-Schwarz in $n_{2}$
$$
\lesssim \sum_{m_{1}
\geq 0}\sum_{0\leq m_{2}}
 \sum_{j_{1},j_{2}\geq 0}
 2^{-m_{2}(1-\epsilon_{0})}
 2^{j_{1}/2}
 \|\chi_{1}\phi_{1,j_{1},m_{1}}\|_{L^{2}}
2^{j_{2}/2}\|\phi_{2,j_{2},m_{2}}\|_{L^{2}}
$$
and a final Cauchy-Schwarz in $m_{2}$ concludes the argument.\\
%%%%%%%%%%%%%%%%%%%%%%%%%%%%%%%%%%
{\bf Case A1d:} $m_{2}>0, \, n_{2}-m_{2}\geq 0$.\\
In this case it is easy to see that
\begin{equation}
2^{n_{2}/2}2^{-m^{*}_{2}}\leq 1.
\label{restr}\end{equation}
Using the same type of estimates presented above, after summing on
$j\geq m^{*}_{2}\geq 0$ we obtain
$$
\eqref{duality1}\lesssim \sum_{m_{1}
\geq 0}\sum_{n_{2}, m_{2}}
 \sum_{j_{1},j_{2}\geq 0}
 2^{-m^{*}_{2}(3/2-\epsilon_{0})}2^{n_{2}/2}
 2^{j_{1}/2}
 \|\chi_{1}\phi_{1,j_{1},m_{1}}\|_{L^{2}}
2^{j_{2}/2}\|\phi_{2, j_2 , m_2 , n_2 }\|_{L^{2}},
$$
and after using \eqref{restr} we can continue with
$$
\lesssim \sum_{m_{1}
\geq 0}\sum_{n_{2}, m_{2}}
 \sum_{j_{1},j_{2}\geq 0}
 2^{-m^{*}_{2}(\half-\epsilon_{0})}
 2^{j_{1}/2}
 \|\chi_{1}\phi_{1,j_{1},m_{1}}\|_{L^{2}}
2^{j_{2}/2}\|\phi_{2, j_2 , m_2 , n_2 }\|_{L^{2}}.
$$
If $n_{2}-m_{2}\geq m_{2}$ then
$$
\lesssim \sum_{m_{1}
\geq 0}\sum_{m_{2}>0}\sum_{n_{2}\geq 2m_{2}}
 \sum_{j_{1},j_{2}\geq 0}
 2^{(-n_{2}+m_{2})(\half-\epsilon_{0})}
 2^{j_{1}/2}
 \|\chi_{1}\phi_{1,j_{1},m_{1}}\|_{L^{2}}
2^{j_{2}/2}\|\phi_{2, j_2 , m_2 , n_2 }\|_{L^{2}}.
$$
and we conclude by Cauchy-Schwarz first with respect to $n_{2}$ and then
$m_{2}$, (here we assume $\epsilon_{0}<\half$).
If $0\leq n_{2}-m_{2}< m_{2}$, we have $m_{2}\leq n_{2}<2 m_{2}$,
$$
\lesssim \sum_{m_{1}
\geq 0}\sum_{m_{2}\leq n_{2}<2 m_{2}}\sum_{m_{2}\geq 0}
 \sum_{j_{1},j_{2}\geq 0}
 2^{-m_{2}(\half-\epsilon_{0})}
 2^{j_{1}/2}
 \|\chi_{1}\phi_{1,j_{1},m_{1}}\|_{L^{2}}
2^{j_{2}/2}\|\phi_{2, j_2 , m_2 , n_2 }\|_{L^{2}}.
$$
and we proceed as in the previous case. This concludes the analysis of
Case A1.\\
{\bf Case A2:} $0\leq j\leq m_{2}^{*}$.\\
{\bf Case A2a:} $m_{2}\geq 0$ or $m_{2}<0$ and $n_{2}>0$.\\
We claim that in this case
$j\leq\max{(j_{1},j_{2}})$. Recall the fundamental identity
\eqref{den}. Then combining this with the restrictions of Case A, we
conclude that $2^{2m_{1}+m_{2}}\lesssim 2^{\max(j,j_{1},j_{2})}$.
If $j\geq \max(j_{1},j_{2})$, then
$2^{2m_{1}+m_{2}}\lesssim 2^{j}$ and if $m_{2}\geq 0$ this implies that
$2^{2m_{1}}\leq \max(2^{m_{2}}, 2^{n_{2}-m_{2}})$, a contradiction
if one compares this with \eqref{conditionsa6}. When $m_{2}<0$ and $
n_{2}\geq 0$, then  $n_{2}-m_{2}\geq 0$ and we obtain
$2^{2m_{1}}\lesssim 2^{n_{2}-2m_{2}}\lesssim 2^{2(n_{2}-m_{2})}$, again
a contradiction if one compares with \eqref{conditionsa6}.
If we assume that $j_{1}=\max(j_{1},j_{2})$, then
$2^{2m_{1}+m_{2}}\lesssim 2^{j_{1}}$, that is
$$
    |\x_{1}|\sim 2^{m_{1}}\lesssim 2^{j_{1}/2-m_{2}/2}.
$$
We can then bound \eqref{duality1} with
\begin{equation}
\lesssim \sum_{m_{1}
\geq 0}\sum_{n_{2}, m_{2}}\sum_{j\leq  m_{2}^{*}}
 \sum_{j_{1},j_{2}\geq 0}
 2^{-j/2}2^{j_{1}/2-m_{2}/2}\int \chi_{j}g_{j}\chi_{1}\theta_{m_{1}}
 \chi_{1}\phi_{1,j_{1},m_{1}}
2^{-m_{2}^{*}(1-\epsilon_{0})}\phi_{2, j_2 , m_2 , n_2 }
\label{a6step5}\end{equation}
We use the Strichartz inequality \eqref{st1}
for $\chi_j g_j \chi_1 \theta_{m_1}$ and for $\phi_{2, j_2 , m_2 , n_2}$ and Plancherel for
$\phi_{1, j_1, m_1}$ and H\"older's inequality to continue the chain
of inequalities with
$$
\lesssim \sum_{m_{1}\geq 0}\sum_{n_{2}, m_{2}}
\sum_{j\leq   m_{2}^{*}}\sum_{j_{1},j_{2}\geq 0}
 2^{j_{1}/2-m_{2}/2}2^{-m_{2}^{*}(1-\epsilon_{0})}
\|\chi_{1}\phi_{1,j_{1},m_{1}}\|_{L^{2}}
2^{j_{2}/2}\|\phi_{2, j_2 , m_2 , n_2 }\|_{L^{2}}.
$$
We now sum over $j$ to get

\begin{equation}
\lesssim \sum_{m_{1}\geq 0}\sum_{n_{2}, m_{2}}
\sum_{j_{1},j_{2}\geq 0}
(1+ m_{2}^{*})2^{-m_{2}/2}2^{-m_{2}^{*}(1-\epsilon_{0})}
2^{j_{1}/2}\|\chi_{1}\phi_{1,j_{1},m_{1}}\|_{L^{2}}
2^{j_{2}/2}\|\phi_{2, j_2 , m_2 , n_2 }\|_{L^{2}}.
\label{a6step4}\end{equation}
Now assume that $m_{2}\geq 0$. We split the $m_{2}$
sum in \eqref{a6step4} into $n_{2}-m_{2}> m_{2}$ and
$n_{2}-m_{2}\leq  m_{2}$. In the first case \eqref{a6step4} becomes
$$
\lesssim \sum_{m_{1}\geq 0}\sum_{n_{2}\geq 2m_{2}}
\sum_{j_{1},j_{2}\geq 0}
[1+(n_{2}-m_{2})]2^{-m_{2}/2}2^{-(n_{2}-m_{2})(1-\epsilon_{0})}
2^{j_{1}/2}\|\chi_{1}\phi_{1,j_{1},m_{1}}\|_{L^{2}}
2^{j_{2}/2}\|\phi_{2, j_2 , m_2 , n_2 }\|_{L^{2}}.
$$
We first use Cauchy-Schwarz on $n_{2}$ and then on $m_{2}$ to
finish. In the second case $n_{2}\leq 2m_{2}$ and $m_{2}^{*}=m_{2}$.
In this case we go back to \eqref{a6step5} and we sum with respect to
$n_{2}$. Then we use Strichartz inequality \eqref{st1}  in the order
$L^{4}L^{2}L^{4}$ (as above) to get
$$
\lesssim \sum_{m_{1}\geq 0}\sum_{m_{2}\geq 0}\sum_{0\leq j\leq m_{2}}
\sum_{j_{1},j_{2}\geq 0}
2^{-m_{2}/2}2^{-m_{2}(1-\epsilon_{0})}
2^{j_{1}/2}\|\chi_{1}\phi_{1,j_{1},m_{1}}\|_{L^{2}}
2^{j_{2}/2}\|\phi_{2,j_{2},m_{2}}\|_{L^{2}}.
$$
Summing in $j$ and then using  Cauchy-Schwarz in $m_2$
will prove the theorem also in this case
\footnote{Observe that if $j_{2}=\max(j_{1},j_{2})$ one
does the same analysis by applying Strichartz inequality \eqref{st1}
in the order
$L^{4}L^{4}L^{2}$.}. Assume now that $m_{2}<0$ and $n_{2}>0$, hence
$m_{2}^{*}=n_{2}-m_{2}\geq 0$.
Then \eqref{a6step4} becomes
$$
 \sum_{m_{1}\geq 0}\sum_{n_{2}\geq 0}\sum_{m_{2}\leq 0}
\sum_{j_{1},j_{2}\geq 0}
[1+ (n_{2}-m_{2})]2^{-m_{2}/2}2^{-(n_{2}-m_{2})(1-\epsilon_{0})}
2^{j_{1}/2}\|\chi_{1}\phi_{1,j_{1},m_{1}}\|_{L^{2}}
2^{j_{2}/2}\|\phi_{2, j_2 , m_2 , n_2 }\|_{L^{2}}.
$$
Observe that
$$2^{-m_{2}/2}(n_{2}-m_{2})2^{-(n_{2}-m_{2})(1-\epsilon_{0})}\lesssim
2^{-(1-\sigma)(1-\epsilon_{0})n_{2}}2^{m_{2}
(-\half+(1-\sigma)(1-\epsilon_{0}))}
$$
for some $0 < \sigma<<1$ and $\epsilon_{0}<\half$. This is  enough for
Cauchy-Schwarz with respect to $n_{2}$ and $m_{2}$.\\
{\bf Case A2b:} $m_{2}, \, n_{2}\leq 0$ \\
In this case $|\x_{2}|, |\m_{2}|\leq 1$. We bound \eqref{duality1}
with
\begin{eqnarray}
 \label{a6step7}&&\sum_{m_{1}\geq 0}
 \sum_{0\leq j}
\sum_{j_{1},j_{2}\geq 0}2^{-j/2}\int_{|\x_{2}|,|\m_{2}|\leq 1}
(|\x_{1}+\x_{2}|^{\half}g_{j}\theta_{m_{1}}
\chi_{1}\chi_{j})\\
\nonumber & &(|\x_{1}|^{\half}\chi_{1}\phi_{1,j_{1},m_{1}})
\frac{\phi_{2,j_{2}}}
{\max{(1, |\x_{2}|,|\m_{2}|/|\x_{2}|)^{1-\epsilon_{0}}}}
    \end{eqnarray}
If $|\m_{2}|^{2}/|\x_{2}|\lesssim 2^{j_{2}}$ then
$|\om(\x_{2},\m_{2})|\lesssim 2^{j_{2}}$ and because
$|\tau_{2}-\om(\x_{2},\m_{2})|\sim 2^{j_{2}}$ we also have
$|\tau_{2}|\lesssim 2^{j_{2}}$. We then apply H\"older's inequalities in
the order $L^{4}_{x}L^{2}_{y}L^{2}_{t}-L^{4}_{x}L^{2}_{y}L^{2}_{t}
-L^{2}_{x}L^{\infty}_{y}L^{\infty}_{t}$ combined with \eqref{smoothj3}
and \eqref{maxfun} to continue the chain of inequalities in
\eqref{a6step7} with
$$
\lesssim \sum_{m_{1}\geq 0}\sum_{0\leq j}
\sum_{j_{1},j_{2}\geq 0}2^{-j/2}2^{j/4}\|g_{j} \theta_{m_1}
\chi_{1}\chi_{j}\|_{L^{2}}
2^{j_{1}/4}\|\phi_{1,j_{1},m_{1}}\chi_{1}\|_{L^{2}}
2^{j_{2}/2}\|\phi_{2,j_{2}}\|_{L^{2}}
$$
and in this case we are done. Assume now that
$|\m_{2}|^{2}/|\x_{2}|>>  2^{j_{2}}$. Then $|\tau_2|\sim
|\m_{2}|^{2}/|\x_{2}|$ and we can rewrite \eqref{a6step7}
as follows
\begin{eqnarray*}
& &\sum_{m_{1}\geq 0}\sum_{0\leq j}
\sum_{j_{1},j_{2}\geq 0}2^{-j/2}\int(|\x_{1}+\x_{2}|^{\half}g_{j}
\theta_{m_{1}}
\chi_{1}\chi_{j})(|\x_{1}|^{\half}\phi_{1,j_{1},m_{1}}\chi_{1})
\frac{\phi_{2,j_{2}}|\m_{2}|^{1-\epsilon_{0}}}
{(|\m_{2}|^{2}/|\xi_{2}|)^{1-\epsilon_{0}}}\\
&\sim& \sum_{m_{1}\geq 0}\sum_{0\leq j}
\sum_{j_{1},j_{2}\geq 0}2^{-j/2}\int(|\x_{1}+\x_{2}|^{\half}g_{j}
\theta_{m_{1}}
\chi_{1}\chi_{j})(|\x_{1}|^{\half}\phi_{1,j_{1},m_{1}}\chi_{1})
\frac{\phi_{2,j_{2}}|\m_{2}|^{1-\epsilon_{0}}}
{|\tau_{2}|^{1-\epsilon_{0}}}\\
&\sim& \sum_{m_{1}\geq 0}\sum_{0\leq j}
\sum_{j_{1},j_{2}\geq 0}2^{-j/2+j/4}\|g_{j}\theta_{m_{1}}
\chi_{1}\|_{L^{2}}2^{j_{1}/4}\|\phi_{1,j_{1},m_{1}}\chi_{1}\|_{L^{2}}
\|\phi_{2,j_{2}}\|_{L^{2}}
\end{eqnarray*}
where in the last step we used again the fact that $|\m_{2}|\leq 1$,
$\epsilon_{0}<\half$,
H\"older's inequality in the order
$L^{4}_{x}L^{2}_{y}L^{2}_{t}-L^{4}_{x}L^{2}_{y}L^{2}_{t}
-L^{2}_{x}L^{\infty}_{y}L^{\infty}_{t}$ and \eqref{smoothj3}
and \eqref{maxfun}. This concludes the analysis of
Case A. \\
{\bf Case B:} $|\x_{1}|\leq 10^{-2}\gia$.\\
Using \eqref{conditionsa6}
one can prove that $|\m_{1}+\m_{2}|/|\x_{1}+\x_{2}|\geq
2|\x_{1}+\x_{2}|$, that
$|\m_{1}+\m_{2}|/|\x_{1}+\x_{2}|\sim |\m_{1}|/|\x_{1}|$, and that
$|\m_{1}+\m_{2}|\sim |\m_{1}|$. We dyadically decompose
$|\m_{1}+\m_{2}|\sim |\m_{1}|\sim 2^{n_{1}}$ and we write
\eqref{duality2} as follows
\begin{equation}
\label{a6step9}\sum_{0\leq j}\sum_{j_{1},j_{2}\geq 0}
\sum_{n_{1}\geq 0}
 2^{-j/2}\int|\x_{1}|g_{j}\theta_{n_{1}}(\m_{1}+\m_{2})
\chi_{2}\chi_{j}\phi_{1,j_{1},n_{1}}\chi_{2}
\frac{\phi_{2,j_{2}}}
{\max{(1, |\x_{2}|,|\m_{2}|/|\x_{2}|)^{1-\epsilon_{0}}}}.
    \end{equation}
We dyadically decompose also $|\x_{2}|\sim 2^{m_{2}}$
and $|\m_{2}|\sim 2^{n_{2}}$. As in Case A,  we define $m^{*}_{2}=
\max(n_{2}-m_{2}, m_{2})$ and we analyze two subcases. \\
{\bf Case B1:} $j\geq \max(0, m^{*}_{2}).$\\
We  use the change of variable
\eqref{chanm}. Now $|J_{\m}|\gtrsim (\gia)^{2}
\gtrsim |\x_{1}|^{2}$. Then we proceed like in Case A1 above, where the
sum in $m_{1}$ is replaced by a sum in $n_{1}$.
\\
{\bf Case B2:} $0\leq j\leq  m^{*}_{2}.$\\
We again consider the two  subcases $m_{2}\geq 0$ or
$m_{2}<0$ and $n_{2}\geq 0$,  and   $m_{2}< 0$ and $ n_{2}<0$.\\
{\bf Case B2a:} $m_{2}\geq 0$ or $m_{2}<0$ and $n_{2}>0$.\\
The fundamental identity \eqref{den} now gives
$$(\gia)^{2}|\x_{2}|\leq (\half-10^{-10})^{-1}2^{\max(j_{1},j_{2},j)}.$$
This again forces $j\leq \max(j_{1},j_{2})$. In fact after setting
$C=(\half-10^{-10})^{-1}$, if $j>
\max(j_{1},j_{2})$ and
$m_{2}\geq 0$, then $(\gia)^{2}\leq C
2^{\max(n_{2}-m_{2},m_{2})}$ and so $(\gia)^{2}\leq C
\max(\gib, |\x_{2}|)$ which is a contradiction in this region.
If $m_{2}< 0$ and $ n_{2}\geq 0$, then $(\gia)^{2}|\x_{2}|
\leq C2^{n_{2}-m_{2}}$ or $(\gia)^{2}\leq C2^{n_{2}-2m_{2}}
\leq C(\gib)^{2}$, which is again a contradiction.
Thus (if for example $j_{2}\leq j_{1}$)
$$
    |\x_{1}|\leq C\gia\leq 10^{-2}2^{j_{1}/2}2^{-m_{2}/2},
$$
and \eqref{a6step9} can be bounded by
\begin{equation}
\sum_{j\leq m^{*}_{2}}\sum_{j_{1},j_{2}\geq 0}
\sum_{n_{1}\geq 0}\sum_{n_{2},m_{2}}
 2^{-j/2}2^{j_{1}/2-m_{2}/2}\int
 \chi_{2}\chi_{j}g_{j,n_{1}}\theta_{n_{1}}
 \phi_{1,j_{1},n_{1}}\chi_{2}
2^{-m_{2}^{*}(1-\epsilon_{0})}\phi_{2,j_{2},m_{2},n_{2}}.
\label{a6step10}\end{equation}
At this point we argue like in Case A2a by replacing
\eqref{a6step5} with \eqref{a6step10}.\\
{\bf Case B2b:} $n_{2},m_{2}\leq 0$.\\
This case can be treated like Case A2b by replacing \eqref{a6step7}
with
\begin{eqnarray}
\label{a6ste10}&&\sum_{n_{1}\geq 0}
 \sum_{0\leq j}
\sum_{j_{1},j_{2}\geq 0}2^{-j/2}\int_{|\x_{2}|,|\m_{1}|\leq 1}
(|\x_{1}+\x_{2}|^{\half}g_{j}\theta_{n_{1}}
\chi_{2}\chi_{j})\\
\nonumber & &(|\x_{1}|^{\half}\chi_{2}\phi_{1,j_{1},n_{1}})
\frac{\phi_{2,j_{2}}}
{\max{(1, |\x_{2}|,|\m_{2}|/|\x_{2}|)^{1-\epsilon_{0}}}}.
    \end{eqnarray}
{\bf Case C:} $10^{-2}\gia\leq |\x_{1}|\leq 10^{2}\gia$.\\
It is easy to show that in this case
$|\m_{1}+\m_{2}|/|\x_{1}+\x_{2}|\lesssim |\x_{1}+\x_{2}|$.
We dyadically decompose with respect to
$|\x_{1}|\sim |\x_{1}+\x_{2}|\sim 2^{m_{1}}$. We go back to
\eqref{a6step1} and we consider two subcases: when $j>m_{1}$
and when $j\leq m_{1}$. \\
{\bf Case C1:} $j>m_{1}$.\\
In this case we  use the change of variable \eqref{chanx}, where now
the free variable is $\x_{2}$ instead of $\x_{1}$. It is easy
to check that also in this case \eqref{jx} holds true and in
particular $|J_{\x}|\gtrsim \gia \gtrsim |\x_{1}|$. We perform a dyadic
decomposition in $\x_{2}$, but only for large frequencies, that is
for $|\x_{2}|\geq 1$. Then \eqref{duality1} becomes
\begin{eqnarray*}
 & &\sum_{m_{1}, m_{2}\geq 0}\sum_{j_{1},j_{2}\geq 0}
 \sum_{j> m_{1}}
 2^{-j/2}2^{m_{1}}\int \chi_{1}g_{j}\theta_{m_{1}}\chi_{j}(u,v,w)\\
 & &\times |J_{\x}|^{-1}H(u,v,w,\x_{2},\theta_{1},\theta_{2})
 dudvdw
 d\x_{2}d\theta_{1}d\theta_{2},
\end{eqnarray*}
where $H(u,v,w,\x_{2},\theta_{1},\theta_{2})$
is the transformation of
$$\chi_{1}\phi_{1,m_{1}}
(\x_{1},\m_{1},\theta_{1}+\om(\x_{1},\m_{1}))
 \chi_{j_{1}}\frac{\phi_{2, m_2 }
 (\x_{2},\m_{2},\theta_{2}+\om(\x_{2},\m_{2}))
 \chi_{2}\chi_{j_{2}}}
 {\max(1,|\x_{2}|,\gib)^{1-\epsilon_{0}}}$$
under the above change of variables. % Then we can continue our
% estimates with
% %
% \begin{eqnarray*}
%  &\lesssim &\sum_{m_{2},m_{1}>0}\sum_{j_{1},j_{2}\geq 0}
%  \sum_{j> m_{1}}
%  2^{-j/2}2^{m_{1}}\|g_{j,m_{1}}\chi_{1}\|_{L^{2}}\\
%  &&\int \left(\int
% |J_{\x}|^{-2}H^{2}_{m_{2}}(u,v,w,\x_{2},\theta_{1},\theta_{2})dudvdw
% \right)^{\half}d\x_{2}d\theta_{1}d\theta_{2}\\
% &\lesssim&\sum_{m_{2},m_{1}>0}\sum_{j_{1},j_{2}\geq 0}
%  \sum_{j> m_{1}}2^{-j/2}2^{m_{1}/2}2^{j_{1}/2}2^{j_{2}/2}2^{m_{2}/2}
% \|g_{j,m_{1}}\chi_{1}\|_{L^{2}}
% \|\phi_{1,j_{1},m_{1}}\chi_{1}\|_{L^{2}}2^{-m_{2}(1-\epsilon_{0})}
% \|\phi_{2,j_{2},m_{2}}\|_{L^{2}}.
% \end{eqnarray*}
% % REPLACED WITH WHAT FOLLOWS ON 10 FEB 2003
We first use Cauchy-Schwarz in $(u,v,w)$ to obtain
\begin{eqnarray*}
  &\lesssim & \sum_{m_1, m_2 \geq 0} \sum_{j_1, j_2 \geq 0} \sum_{j> m_1 }
2^{-j/2} 2^{m_1} {{\| g_j \theta_{m_1} \chi_1 \|}_{L^2}} \\
&\times & \left( \int \left( |J_\xi |^{-2} H^2 ( u, v, w , \xi_2 , \theta_1 , \theta_2 )
du dv dw \right)^\half d\xi_2 d\theta_1 d\theta_2 \right)
\end{eqnarray*}
We now use the lower bound $|J_\xi| \gtrsim 2^{m_1}$, and Cauchy-Schwarz in $\xi_2 , \theta_1 
, \theta_2,$ to obtain
\begin{eqnarray*}
 & \lesssim& \sum_{m_1, m_2 \geq 0} \sum_{j_1 , j_2 \geq 0} \sum_{j > m_1 }
2^{-j/2} 2^{m_1 /2} 2^{m_2 /2 } 2^{j_1 /2 } 2^{j_2 /2} \\
&\times & {{\| g_j \theta_{m_1} \chi_1 \|}_{L^2}} \left( \int H^2 ( u , v, w , \xi_2 , \theta_1 ,
\theta_2 ) |J_\xi |^{-1} du dv dw d\xi_2 d\theta_1 d\theta_2 \right)^\half,
\end{eqnarray*}
which, upon undoing the change of variables yields (since $\max ( 1 , |\xi_2 | , |\mu_2 |/ |\xi_2 |) \geq 2^{m_2}$)
\begin{eqnarray*}
  &\lesssim & \sum_{m_1, m_2 \geq 0} \sum_{j_1 , j_2 \geq 0} \sum_{j > m_1 }
2^{-j/2} 2^{m_1 /2} 2^{m_2 /2 } 2^{j_1 /2 } 2^{j_2 /2} \\
&\times & {{\| g_j \theta_{m_1} \chi_1 \|}_{L^2}} {{\| \phi_{1, j_1 , m_1 } \chi_1 \|}_{L^2}}
{{\| \phi_{2, j_2 , m_2 } \|}_{L^2}} 2^{-m_2 (1 - \epsilon_0 )}.
\end{eqnarray*}
Then a sum in $j$ and  Cauchy-Schwarz in $m_{2}$ gives the result. \\
{\bf Case C2:} $j\leq m_{1}$.\\
Let's introduce the region
\begin{equation}
    R=\{10^{-2}\gia\leq |\x_{1}|\leq 10^{2}\gia\}\cap A_{6}
\end{equation}
This is the region where we need to introduce the space
$Y_{1-\epsilon_{0},\epsilon_{0},\half}$. We
go back to \eqref{duality-1} and this time we keep $|\hat v|$ and we
only normalize $|\hat u|$. Then \eqref{duality-1} becomes
\begin{eqnarray}
 \label{a6step11}& &\sum_{m_{1}\geq 0}\sum_{0\leq j\leq
 m_{1}}\sum_{j_{1},j_{2}\geq 0}2^{-j/2}
 \int_{R} \chi_{1}\theta_{m_{1}}g_{j}(\x_{1}+\x_{2},\m_{1}+\m_{2},
 \theta_{1}+
 \om(\x_{1},\m_{1})+\theta_{2}+\om(\x_{2},\m_{2}))\\
 \nonumber & &\chi_{1}\chi_{j}(\theta_{1}+
 \om(\x_{1},\m_{1})+\theta_{2}+\om(\x_{2},\m_{2})
 -\om(\x_{1}+\x_{2},\m_{1}+\m_{2}))
 |\x_{1}|\\
 \nonumber &&\chi_{1}\chi_{j_{1}}\phi_{1,m_{1}}(\x_{1},\m_{1},\theta_{1}+
 \om(\x_{1},\m_{1}))\chi_{j_{2}}|\hat{v}|(\x_{2},\m_{2},\theta_{2}+
 \om(\x_{2},\m_{2})).
 \end{eqnarray}
Define $D_{m_{1},m_{2}}$ to be the dyadic block such that
$|\x_{i}|\sim 2^{m_{i}}, i=1,2$. We observe that for fixed
$(\x_{1},\xi_{2}, \m_{1},
\theta_{1}, \theta_{2})$, the set of $\m_{2}$ such that
$(\x_{1},\xi_{2}, \m_{1}, \m_{2},\theta_{1}, \theta_{2})\in R\cap
D_{m_{1},m_{2}}$ and such that \eqref{den1} is true, is
a union of two symmetric intervals  with length satisfying
\eqref{deltam} .
Similarly, for fixed
$(\xi_{2}, \m_{1}, \m_{2},\theta_{1}, \theta_{2})$ the set of $\x_{1}$
such that
$(\x_{1},\xi_{2}, \m_{1}, \m_{2},\theta_{1}, \theta_{2})\in R\cap
D_{m_{1},m_{2}}$ and such that \eqref{den1} is true, is
a union of two symmetric intervals  with length satisfying
\begin{equation}
    |\Delta_{\x_{1}}|\lesssim 2^{j-m_{1}-m_{2}}.
\label{deltam1}\end{equation}
To prove this it's enough to use the mean value theorem, and estimate
from below $|g'(\xi_1)|$, where
$$g(\x_1 )=\theta_{1}+\theta_{2}+\frac{\x_1 \x_{2}}{(\x_1 +\x_{2})}
    \left(\left(\frac{\m_{1}}{\x_1 }-\frac{\m_{2}}{\x_{2}}\right)^{2}
    -3(\xi_1 +\xi_{2})^{2}\right).$$
After a short calculation one has
$$g'(\x_1 )=\frac{\x_{2}^{2}}{(\x_1 +\x_{2})^{2}}
\left(\left(\frac{\m_{1}}{\x_1 }-\frac{\m_{2}}{\x_{2}}\right)^{2}
    -3(\xi_1 +\xi_{2})^{2}
\right) - \xi_{2}\left(\frac{2\m_{1}}{\xi_1 (\x_1 +\x_{2})}
\left(\frac{\mu_{1}}{\x_1 }-\frac{\mu_{2}}{\x_{2}}\right)+6\x_1 \right).$$
Note that in $R$
$$\left|\frac{\x_{2}^{2}}{(\x_{1}+\x_{2})^{2}}
\left(\left(\frac{\m_{1}}{\x_{1}}-\frac{\m_{2}}{\x_{2}}\right)^{2}
    -3(\xi_{1}+\xi_{2})^{2}
\right)\right| \leq C\frac{|\x_{2}|^{2}}{|\x_{1}|^{2}}
\left(\frac{|\m_{1}|}{|\x_{1}|}\right)^{2}\leq C10^{-10}
|\x_{2}|^{2}.$$
On the other hand it's easy to check that
$\x_1$ and $\frac{2\m_{1}}{\xi_1 (\x_1 +\x_{2})}
\left(\frac{\mu_{1}}{\x_1 }-\frac{\mu_{2}}{\x_{2}}\right)$ have the same
sign, hence
$$\left|\xi_{2}\left(\frac{2\m_{1}}{\xi_1 (\x_1 +\x_{2})}
\left(\frac{\mu_{1}}{\x_1 }-\frac{\mu_{2}}{\x_{2}}\right)+6\x_1 \right)
\right|\geq 6|\x_{2}|\x_1 |$$
and the claim follows.\\
{\bf Case C2a:} $m_{2}\geq 0$\\
In this case we use the change of variables \eqref{chanx}, and
by \eqref{jx}, $|J_{\x}|\gtrsim 2^{m_{1}}$. We can write \eqref{a6step11}
as
$$
\sum_{m_{1}\geq 0}\sum_{0\leq j\leq m_{1}}\sum_{m_{2},j_{1},j_{2}\geq 0}
2^{-j/2}2^{m_{1}}\int \chi_{1}\chi_{j}g_{j,m_{1}}(u,v,w)|J_{\x}|^{-1}
H(u,v,w, \xi_{1},\theta_{1},\theta_{2})dudvdwd\xi_{1}d\theta_{1}
d\theta_{2},
$$
where $H(u,v,w,\xi_{1},\theta_{1},\theta_{2})$ is the transformation
of
$$\chi_{j}\chi_{1}\chi_{j_{1}}\phi_{1,m_{1}}(\x_{1},\m_{1},\theta_{1}+
 \om(\x_{1},\m_{1}))\chi_{j_{2}}|\hat{v}|(\x_{2},\m_{2},
\theta_{2}+
 \om(\x_{2},\m_{2})) \theta_{m_2} (\xi_2 )$$
 under the above change of variables. 
We will define 
$$\widehat{v}_{j_2 , m_2 } (\xi_2 , \mu_2 , \tau_2 ) = \chi_{j_2} (\tau_2 -
\omega ( \xi_2 , \mu_2 )) \widehat{v} ( \xi_2 , \mu_2 , \tau_2 ) \theta_{m_2 } (\xi_2 ).$$
We can then continue the
 estimate with
\begin{eqnarray*}
&\lesssim&\sum_{m_{1}\geq 0}\sum_{0\leq j\leq m_{1}}
\sum_{m_{2},j_{1},j_{2}\geq 0}
2^{-j/2}2^{m_{1}/2}2^{(j-m_{1}-m_{2})/2}\int
\chi_{j}\chi_{1}g_{j,m_{1}}(u,v,w)\\
&&\left(\int H^{2}(u,v,w, \xi_{1},\theta_{1},\theta_{2})
|J_{\x}|^{-1}d\x_{1}\right)^{\half}dudvdwd\theta_{1}
d\theta_{2}\\
&\lesssim&\sum_{0\leq j\leq m_{1}}\sum_{m_{2}, m_{1},j_{1},j_{2}\geq 0}
2^{-j/2}2^{m_{1}/2}2^{(j-m_{1}-m_{2})/2}
\|g_{j}\theta_{m_{1}}\chi_{j}\chi_{1}\|_{L^{2}}2^{j_{1}/2}2^{j_{2}/2}\\
& &\left(\int
|\chi_{1}\phi_{m_{1},j_{1}}
(\xi_{1},\m_{1},\theta_{1}+\om(\x_{1},\m_{1}))|^{2}
|\hat{v}_{m_{2},j_{2}}
(\xi_{2},\m_{2},\theta_{2}+\om(\x_{2},\m_{2}))|^{2}\right.\\
&&\left.\chi_{j}(\theta_{1}+
 \om(\x_{1},\m_{1})+\theta_{2}+\om(\x_{2},\m_{2})
 -\om(\x_{1}+\x_{2},\m_{1}+\m_{2}))\right)^{\half}.
\end{eqnarray*}
We now use  H\"older inequality in $\m_{2}$ and Lemma \ref{sobolev}.
More precisely  set $w(\x,\m)=(1+|\x|+\gi)$, then we write, for
$\theta=(p-2)/2p, \, p=2r, \, r>1, 1 = \frac{1}{r} + \frac{1}{r'}$,
\begin{equation}
\int |\hat{v}_{m_{2},j_{2}}|^{2}\chi_{j}\leq
|\Delta_{\m_{2}}|^{1/r'}\|\hat{v}_{m_{2},j_{2}}\|_{L^{2r}}^{2}
\lesssim 2^{(j-m_{1})/r'}\|w^{\epsilon_{0}}\hat{v}_{m_{2},j_{2}}
\|_{L^{2}}^{2(1-\theta)}
\|w^{-\epsilon_{0}}\partial_{\m_{2}}\hat{v}_{m_{2},j_{2}}
\|_{L^{2}}^{2\theta}.
    \label{rr'}\end{equation}
We insert this in the chain of inequalities above and
we continue with
\begin{eqnarray*}
&\lesssim&\sum_{m_{1}\geq 0}\sum_{0\leq j\leq m_{1}}\sum_{m_{2}\geq 0}
\sum_{j_{1},j_{2}\geq 0}
2^{-j/2}2^{m_{1}/2}2^{(j-m_{1}-m_{2})/2}2^{j_{1}/2}2^{j_{2}/2}
2^{(j-m_{1})/2r'}\|g_{j}\theta_{m_{1}}\chi_{1}\chi_{j}\|_{L^{2}}
\|\chi_{1}\phi_{m_{1},j_{1}}\|_{L^{2}}\\
&&\left(\int \|w^{\epsilon_{0}}\hat{v}_{m_{2},j_{2}}
(\xi_{2},\m_{2},\theta_{2}+\om(\x_{2},\m_{2}))
\|_{L^{2}_{\m_{2}}}^{2(1-\theta)}
\|w^{-\epsilon_{0}}\partial_{\m_{2}}
\hat{v}_{m_{2},j_{2}}
(\xi_{2},\m_{2},\theta_{2}+\om(\x_{2},\m_{2}))
\|_{L^{2}_{\m_{2}}}^{2\theta}
d\xi_{2}d\theta_{2}\right)^{\half},
\end{eqnarray*}
and after applying H\"older's inequality with respect to $\xi_{2}, \theta_{2}$ with
exponents $\theta^{-1}$ and $(1-\theta)^{-1}$,
\begin{eqnarray*}
&\lesssim&\sum_{m_{1}\geq 0}\sum_{0\leq j\leq m_{1}}\sum_{m_{2}\geq 0}
\sum_{j_{1},j_{2}\geq 0}
2^{-j/2}2^{m_{1}/2}2^{(j-m_{1}-m_{2})/2}2^{j_{1}/2}2^{j_{2}/2}
2^{(j-m_{1})/2r'}\|g_{j}\theta_{m_{1}}\chi_{1}\chi_{j}\|_{L^{2}}\\
&&\|\chi_{1}\phi_{m_{1},j_{1}}\|_{L^{2}}
\|w^{\epsilon_{0}}\hat{v}_{m_{2},j_{2}}\|_{L^{2}}^{(1-\theta)}
\|w^{-\epsilon_{0}}\partial_{\m_{2}}
\hat{v}_{m_{2},j_{2}}\|_{L^{2}}^{\theta}.
\end{eqnarray*}
We first  sum on $0\leq j\leq m_{1}$, then  on $m_{1}$
and $j_{1}$ so that the norm $\|u\|_{X_{1-\epsilon_{0}},\half}$ appears.
After Cauchy-Schwarz in $m_{2}$ we are left with the following term to
estimate
$$\sum_{j_{2}\geq 0}2^{j_{2}/2}\left(\sum_{m_{2}\geq 0}
\|w^{\epsilon_{0}}\hat{v}_{m_{2},j_{2}}\|_{L^{2}}^{2(1-\theta)}
\|w^{-\epsilon_{0}}\partial_{\m_{2}}
\hat{v}_{m_{2},j_{2}}\|_{L^{2}}^{2\theta}\right)^{\half}.
$$
We use a H\"older inequality with respect to the sum on $m_{2}$
to obtain
\begin{eqnarray*}&\leq &\sum_{j_{2}\geq 0}2^{j_{2}/2}
\left(\sum_{m_{2}\geq 0}
\|w^{\epsilon_{0}}\hat{v}_{m_{2},j_{2}}\|_{L^{2}}^{2}
\right)^{(1-\theta)/2}
\left(\sum_{m_{2}\geq 0}\|w^{-\epsilon_{0}}\partial_{\m_{2}}
\hat{v}_{m_{2},j_{2}}\|_{L^{2}}^{2}\right)^{\theta/2}\\
&=&\sum_{j_{2}\geq 0}2^{j_{2}/2}\|w^{\epsilon_{0}}\hat{v}_{j_{2}}
\|_{L^{2}}^{(1-\theta)}
\|w^{-\epsilon_{0}}\partial_{\m_{2}}\hat{v}_{j_{2}}\|_{L^{2}}^{\theta},
\end{eqnarray*}
where 
\begin{equation*}
  \widehat{v}_{j_2} = \sum_{m_2 \geq 0}  \widehat{v}_{j_2 , m_2 } = \chi_{j_2}
( \tau - \omega ( \xi_2 , \mu_2 )) \widehat{v} ( \xi_2 , \mu_2 , \tau_2 ).
\end{equation*}
 then  a H\"older inequality in $j_{2}$ to finish with
\begin{equation}
\left(\sum_{j_{2}\geq 0}2^{j_{2}/2}\|w^{\epsilon_{0}}
\hat{v}_{j_{2}}\|_{L^{2}}
\right)^{(1-\theta)}
\left(\sum_{j_{2}\geq 0}2^{j_{2}/2}\|w^{-\epsilon_{0}}
\partial_{\m_{2}}
\hat{v}_{j_{2}}\|_{L^{2}}\right)^{\theta}.
\label{a6step12}\end{equation}
Clearly the first coefficient of \eqref{a6step12} is controlled by
$\|v\|_{X_{\epsilon_{0},\half}}^{1-\theta}$. For the second one we write
\begin{eqnarray*}
&&\partial_{\m_{2}}(\hat{v}_{j_{2}}(\x_{2},\m_{2},
\theta_{2}+\om(\x_{2},\m_{2})))=\partial_{\m_{2}}
(\chi_{j_{2}}(\theta_{2})\hat{v}(\x_{2},\m_{2},
\theta_{2}+\om(\x_{2},\m_{2})))\\
&=&\chi_{j_{2}}\partial_{\m_{2}}
\hat{v}(\x_{2},\m_{2},
\theta_{2}+\om(\x_{2},\m_{2})+
2\m_{2}/\x_{2}\chi_{j_{2}}\partial_{\tau_{2}}
\hat{v}(\x_{2},\m_{2},
\theta_{2}+\om(\x_{2},\m_{2})
    \end{eqnarray*}
    which shows that the second term is controlled by
    $\|v\|_{Y_{1-\epsilon_{0},-\epsilon_{0},\half}}^{\theta}$, with
    $\theta=(p-2)/2p, \, p=2r, \, r>1$.\\
{\bf Case C2b:} $m_{2}< 0$\\
If $j-m_{1}\leq 2m_{2}$ we proceed like in Case C2a and we obtain
\begin{eqnarray*}
&\lesssim&\sum_{m_{1}\geq 0}\sum_{m_{2}<0}\sum_{0\leq j\leq m_{1}+2m_{2}}
\sum_{j_{1},j_{2}\geq 0}
2^{-j/2}2^{m_{1}/2}2^{(j-m_{1}-m_{2})/2}2^{j_{1}/2}2^{j_{2}/2}
2^{(j-m_{1})/2r'}\|g_{j}\theta_{m_{1}}\chi_{1}\chi_{j}\|_{L^{2}}
\|\chi_{1}\phi_{m_{1},j_{1}}\|_{L^{2}}\\
&&\|w^{\epsilon_{0}}\hat{v}_{m_{2},j_{2}}\|_{L^{2}}^{(1-\theta)}
\|w^{-\epsilon_{0}}\partial_{\m_{2}}
\hat{v}_{m_{2},j_{2}}\|_{L^{2}}^{\theta}.
\end{eqnarray*}
We sum on $j$, and then choose $r'<2$ so that we can use the fact that
$m_{2}<0$ and Cauchy-Schwarz in
$m_{2}$, to finish like in Case C2a.
If $j-m_{1}> 2m_{2}$ we use the change of variables \eqref{chanx},
where we leave the variable $\x_{2}$ free. It is easy to check
that  we have $|J_{\x}|\gtrsim |\x_{1}|$. Arguing as
in Case C2a we are led to
\begin{eqnarray*}
&\lesssim&\sum_{m_{1}\geq 0}\sum_{0\leq j\leq m_{1}}
\sum_{m_{2}<(j-m_{1})/2}\sum_{j_{1},j_{2}\geq 0}
2^{(m_{1}-j)/2}2^{j_{1}/2}2^{j_{2}/2}
2^{(j-m_{1})/2r'}2^{m_{2}/2}\|g_{j,m_{1}}\chi_{1}\chi_{j}\|_{L^{2}}
\|\chi_{1}\phi_{m_{1},j_{1}}\|_{L^{2}}\\
&&\|w^{\epsilon_{0}}\hat{v}_{m_{2},j_{2}}\|_{L^{2}}^{(1-\theta)}
\|w^{-\epsilon_{0}}\partial_{\m_{2}}
\hat{v}_{m_{2},j_{2}}\|_{L^{2}}^{\theta}.
\end{eqnarray*}
We first apply Cauchy-Schwarz with respect to $m_{2}$ and we obtain an
extra factor $2^{(-m_{1}+j)/4}$. If we now have that
$\frac{1}{4}+ \frac{1}{2r'}> \half$, we can then sum in $j$, and 
repeat the argument
 in Case
C2a. This is again the restriction  $r'<2 $, which when we go to
 \eqref{fp} gives $p=2r>4$, and hence
$\theta>\frac{1}{4}$.

This concludes the proof of Theorem \ref{xbil}.

\end{proof}

%%%%%%%%%%%%%%%%%%%%%%%%%%%Proof of Theorem \ref{ybil}}%%%%%%%%%%%%
                                                                  %

\begin{proof}[{Proof of Theorem \ref{ybil}}]
 We first check the part of the statement involving the term
 $\partial_{\tau}(\widehat{\partial_{x}(uv)})$. We proceed by
 writing
 \begin{eqnarray*}
& &\partial_{\tau}\left(\x\int
\hat{v}(\x-\x_{1},\m-\m_{1},\tau-\tau_{1})\hat{u}
(\x_{1},\m_{1},\tau_{1})d\x_{1}d\m_{1}d\tau_{1}\right)\\
&=&\partial_{\tau}\left(\x\int_{|\x-\x_{1}|\leq |\x_{1}|}
\hat{v}(\x-\x_{1},\m-\m_{1},\tau-\tau_{1})\hat{u}
(\x_{1},\m_{1},\tau_{1})d\x_{1}d\m_{1}d\tau_{1}\right)\\
&+&\partial_{\tau}\left(\x\int_{|\x-\x_{1}|\geq |\x_{1}|}
\hat{v}(\x-\x_{1},\m-\m_{1},\tau-\tau_{1})\hat{u}
(\x_{1},\m_{1},\tau_{1})d\x_{1}d\m_{1}d\tau_{1}\right)\\
&=&\x\int_{|\x-\x_{1}|\leq |\x_{1}|}
\hat{v}(\x-\x_{1},\m-\m_{1},\tau-\tau_{1})\partial_{\tau_{1}}\hat{u}
(\x_{1},\m_{1},\tau_{1})d\x_{1}d\m_{1}d\tau_{1}\\
&+&\x\int_{|\x-\x_{1}|\geq |\x_{1}|}
\partial_{\tau}\hat{v}(\x-\x_{1},\m-\m_{1},\tau-\tau_{1})\hat{u}
(\x_{1},\m_{1},\tau_{1})d\x_{1}d\m_{1}d\tau_{1}.
     \end{eqnarray*}
Then the estimates for each one of these two terms follows from the
proof of
Theorem \ref{xbil}. Next we check the part involving $\partial_{\m}
\widehat{\partial_{x}(uv)}$. If we proceed as above, we only need to
check
$$I(\x,\m,\tau)=\frac{|\x|}{\max(1,|\x|,|\m|/|\x|)^{\epsilon_{0}}}
\int_{|\x-\x_{1}|\leq  |\x_{1}|}
|\hat{v}|(\x-\x_{1},\m-\m_{1},\tau-\tau_{1})|\partial_{\m_{1}}\hat{u}|
(\x_{1},\m_{1},\tau_{1})d\x_{1}d\m_{1}d\tau_{1}.
$$
We introduce the two functions
\begin{eqnarray}
 \label{normalizu1} \phi_{1, j_1}(\x,\m,\tau)&=&
 \max(1,|\x|,|\m|/|\x|)^{-\epsilon_{0}}|\partial_{\m}\hat{u}|
 \chi_{j_1 }(\x,\m,\tau),\\
\label{normalizv1}
\phi_{2, j_2}(\x,\m,\tau)&=&\max(1,|\x|,|\m|/|\x|)^{1-\epsilon_{0}}
|\hat{v}|\chi_{j_2 }(\x,\m,\tau).
\end{eqnarray}
corresponding to \eqref{normalizeu} and \eqref{normalizev}  in Theorem
\ref{xbil}. We now observe that in the region $|\x-\x_{1}|\leq
|\x_{1}|$ we also have $\gia\leq |\m-\m_{1}|/|\x-\x_{1}| +
2\gi$, so
\begin{equation}I(\x,\m,\tau)\lesssim J^{1}+J^{2},
\label{ybilstep1}\end{equation}
where %{\bf [GAFA: labels encroach display. help?]}
\begin{eqnarray}
\label{j1}J^{1}&=&\frac{|\x|}{\max(1,|\x|,\gi)^{\epsilon_{0}}}
\int_{*}\phi_{1,j_{1}}(\x_{1},\m_{1},\tau_{1})
\frac{\phi_{2,j_{2}}}
{\max(1,|\x-\x_{1}|,|\m-\m_{1}|/|\x-\x_{1}|)^{1-2\epsilon_{0}}}\\
\label{j2}J^{2}&=&|\x|
\int_{*}\phi_{1,j_{1}}(\x_{1},\m_{1},\tau_{1})
\frac{\phi_{2,j_{2}}}
{\max(1,|\x-\x_{1}|,|\m-\m_{1}|/|\x-\x_{1}|)^{1-\epsilon_{0}}}.
\end{eqnarray}
where $\int_{*}$ is the integral over the region given by
 $|\x-\x_{1}|\leq |\x_{1}|$. Also notice that if
\begin{equation}\tilde J^{2}=|\x|
\int_{*}\phi_{1,j_{1}}(\x_{1},\m_{1},\tau_{1})
\frac{\phi_{2,j_{2}}}
{\max(1,|\x-\x_{1}|,|\m-\m_{1}|/|\x-\x_{1}|)^{1-2\epsilon_{0}}},
\label{tildej2}\end{equation}
then
\begin{equation}
    J^{1}+J^{2}\leq \tilde J^{2}.
\label{j1j2}\end{equation}
 By duality we need to estimate
\begin{eqnarray}
\label{ybildualx}& &\sum_{i=1,2}\sum_{j\geq 0}\sum_{m\geq 0}
2^{-j/2} \int g_{j}(\x,\m,\tau)J^{i}(\x,\m,\tau)\theta_{m}(\x)
    \chi_{1}(\x,\m)\chi_{j}(\tau-\om(\x,\m))d\x d\m d\tau\\
 \label{ybildualm}& &\sum_{i=1,2} \sum_{j\geq 0}
 \sum_{n\geq 0}2^{-j/2} \int g_{j}(\x,\m,\tau)J^{i}(\x,\m,\tau) \theta_{n}(\m)
    \chi_{2}(\x,\m)\chi_{j}(\tau-\om(\x,\m)) d\x d\m d\tau.
\end{eqnarray}
We  begin a case by case analysis for \eqref{ybildualx} and
\eqref{ybildualm}. We introduce
again the notation $\m=\m_{1}+\m_{2}$ and $\x=\x_{1}+\x_{2}$. \\
{\bf Region $A_{1}$}. In this region we use \eqref{j1j2}
and we  bound \eqref{ybildualx} and \eqref{ybildualm} by replacing
$J^{i}, i=1,2$ with $\tilde J^{2}$ in \eqref{tildej2}.\\
{\bf Case A:}  $|\x|\geq \half\gi$.\\
We only have to estimate \eqref{ybildualx}.  We observe that
$|\x|\leq 2$ and hence, from \eqref{tildej2}, 
$$\tilde J^{2}\lesssim \int_{*}\phi_{1,j_{1}}(\x_{1},\m_{1},\tau_{1})
\phi_{2,j_{2}}(\x-\x_{1},\m-\m_{1},\tau-\tau_{1}).$$
Then if we dyadically decompose for $1+ |\x|\sim 2^{m}$, the sum in $m$ is
a finite sum and we can use Strichartz inequality \eqref{st1}, applied to $\phi_{1, j_1 }$ 
and $\phi_{2, j_2}$, and Plancherel and H\"older, 
to prove the theorem in this case.\\
{\bf Case B:} $|\x|< \half \gi$.\\
We  have to estimate \eqref{ybildualm}. If
$\gi\leq 1$, it follows that $|\m_{1}+\m_{2}|\leq 2$, and we proceed
like in Case A. If  $\gi\geq  1$. We consider two subcases.\\
{\bf Case B1:} $|\m_{1}|\lesssim |\m_{2}|$.\\
If $\half \gib\leq |\x_{2}|$, because $|\x_{2}|\leq 1$, it follows that
$|\m_{1}|, |\m_{2}|\leq 4$ and we go back to the previous case.
If $\half \gib> |\x_{2}|$, we dyadically decompose with
respect to $\m_{2}$. Then
$$|\x_{1}+\x_{2}|\leq |\x_{1}+\x_{2}|^{2\epsilon_{0}}
|\m_{1}+\m_{2}|^{1-2\epsilon_{0}}\leq 2|\m_{1}+
\m_{2}|^{1-2\epsilon_{0}}$$
this corresponds to Case B1 of region $A_{1}$
in the proof of Theorem \ref{xbil}. In fact we
have the following bound for \eqref{ybildualm}:
\begin{eqnarray*}
& &\sum_{j_{1},j_{2}\geq 0}\sum_{j\geq 0}\sum_{n_{2}\geq 0}
\sum_{0\leq r\leq n_{2}}
2^{-j/2}\int_{A_{1}}g_{j}\chi_{2}\chi_{j}(\x,\m,\tau)
 \theta_{n_2 + 1 -r}2^{(n_{2}+1-r)(1-2\epsilon_{0})}\\
 &&
 \phi_{1,j_{1}}(\x_{1},\m_{1},\tau_{1})
 \frac{\chi_{2}\phi_{2,j_{2},n_{2}}(\x_{2},\m_{2},\tau_{2})}
 {(\gib)^{1-2\epsilon_{0}}},
\end{eqnarray*}
where $\phi_{2, j_2 , n_2} = \theta_{n_2} ( \mu_2 ) \phi_{2, j_2}$. (Note that 
$\phi_{1, j_1} , \phi_{2, j_2}$ are in this case as defined in \eqref{normalizu1}, \eqref{normalizv1},
and one concludes like in that case for $\epsilon_{0}<\frac{1}{4}$.\\
{\bf Case B2:} $|\m_{1}|>>|\m_{2}|$.\\
Then $|\m_{1}|\sim |\m_{1}+\m_{2}|$. If
$\half \gia\leq |\x_{1}|$, it follows that
$|\m_{1}+\m_{2}|\sim |\m_{1}|\leq 4$. So if we set $1+ |\m_{1}+\m_{2}|\sim
2^{n}$, then the sum on $n$ is finite and we use the Strichartz
inequality \eqref{st1}. If $\half \gia> |\x_{1}|$, we
 dyadically
decompose  $|\m_{1}+\m_{2}|\sim |\m_{1}|\sim 2^{n_{1}}$.
Then we use Strichartz inequality \eqref{st1} again.

\noindent
{\bf Region $A_{2}$}. \\
{\bf Case A:} $|\x|\geq \half \gi$.\\
This case is similar to the corresponding case in Theorem \ref{xbil}.
If $|\x_{1}|\sim |\x_{2}|\sim 2^{m_{1}}$ and
$|\x_{1}+\x_{2}|\sim 2^{m_{1}+1-r}, 0\leq r<m_{1}$, (and
$|\x_{1}+\x_{2}|\leq 1$ when $r=m_{1}$), then
one can bound \eqref{ybildualx} with
$$
\sup_{j,m}\|g_{j}\chi_{1}\chi_{j}\theta_{m}\|_{L^{2}}
\sum_{m_1 \geq 0}\sum_{0\leq r \leq  m_{1}}2^{-r(1-2\epsilon_{0})}
\Pi_{i=1,2}\sum_{j_{i}\geq 0}2^{j_{i}/2}\|\phi_{i,j_{i},m_{1}}\|_{L^{2}},
$$
and this proves the estimate for $\epsilon_{0}<\frac{1}{4}$. \\
{\bf Case B:} $|\x|< \half \gi$.\\
If $\gi\leq 1$, then $|\x|, |\m|\leq 4$ and we go back to the estimate
in region $A_{1}$.
So we assume that $\gi> 1$. We consider two subcases.\\
{\bf Case B1:} $|\m_{1}| \lesssim |\m_{2}|$\\
If $\half \gib\leq |\x_{2}|$ we use the fact that
\begin{equation}
|\x_{1}+\x_{2}|\leq
1/\sqrt{2}|\m_{1}+\m_{2}|^{\half}\leq C|\m_{2}|^{\half}\leq C|\x_{2}|,
\label{estimate}\end{equation}
and we dyadically decompose with respect to $\x_{2}$. Then if
$|\x_{2}|\sim 2^{m_{2}}$, $1+ |\m_{1}+\m_{2}|\sim 2^{-r}2^{2m_{2}}, \, \,
0\leq r\leq 2m_{2}$. We proceed now like in the corresponding case for
Theorem \ref{xbil} where we replace $\epsilon_{0}$ by $2\epsilon_{0}$
and $r$ by $r/2$.
If $\half \gib> |\x_{2}|$  we
dyadically decompose with respect to $\m_{2}$ by setting
$|\m_{2}|\sim 2^{n_{2}}$. Then
$|\m_{1}+\m_{2}|\sim 2^{n_{2}-r+1}, \, \, 0\leq r\leq n_{2}$ and from
\eqref{estimate} also $|\x_{1}+\x_{2}|\lesssim  2^{n_{2}/2-r/2+1/2}$.
We reduce our estimate to
\begin{eqnarray*}
 & &   \sum_{j_{1},j_{2}\geq 0}
 \sum_{j\geq 0}\sum_{n_{2}\geq 0}\sum_{0\leq r\leq n_{2}}
2^{-j/2}\int_{A_{2}}g_{j}\chi_{2}\chi_{j}(\x,\m,\tau)
\theta_{n_2 + 1 -r}2^{(n_{2}+1-r)(1-2\epsilon_{0})/2}\\
& &\times|\x_{1}+\x_{2}|^{2\epsilon_{0}}
 \phi_{1,j_{1}}(\x_{1},\m_{1},\tau_{1})
 \frac{\chi_{2}\phi_{2,j_{2},n_{2}}(\x_{2},\m_{2},\tau_{2})}
 {(\gib)^{1-2\epsilon_{0}}}\\
 &\lesssim&\sum_{j_{1},j_{2}\geq 0}
 \sum_{n_{2}\geq 0}\sum_{0\leq r\leq n_{2}}\sum_{j\geq 0}
 2^{-j/2}\int_{A_{2}}g_{j}\chi_{2}\chi_{j}(\x,\m,\tau)
\theta_{n_2 + 1 -r}2^{-r(1-2\epsilon_{0})/2}\\
& &\times
|\x_{1}+\x_{2}|^{2\epsilon_{0}}
 \phi_{1,j_{1}}(\x_{1},\m_{1},\tau_{1})
 \frac{|\x_{2}|^{1-2\epsilon_{0}}}{|\m_{2}|^{(1-2\epsilon_{0})/2}}
 \chi_{2}\phi_{2,j_{2},n_{2}}(\x_{2},\m_{2},\tau_{2})
\end{eqnarray*}
We use again Lemma \ref{mainlemma} and, in view of the fact that
here $|\x_{2}|^{2}/|\m_{2}|\leq \half$, we continue with
$$
\sup_{j,n}\|g_{j}\chi_{2}\chi_{j}
 \theta_{n}\|_{L^{2}}\sum_{n_{2}\geq 0}\sum_{0\leq r\leq n_{2}}
2^{-r(1-2\epsilon_{0})/2}\left(\sum_{j_{1}\geq 0}2^{j_{1}/2}
\|\phi_{1,j_{1}}\|_{L^{2}}\right)
 \left(\sum_{j_{2}\geq 0}2^{j_{2}/2}
 \|\chi_{2}\phi_{2,n_{2},j_{2}}\|_{L^{2}}\right),
$$
and we  sum in $r$.
This is enough to prove the estimate  as long as $\epsilon_{0}<\frac{1}{4}$. \\
{\bf Case B2:} $|\m_{1}|>>|\m_{2}|$. \\
 If $\half \gia\leq |\x_{1}|$,
then
$$
|\x_{1}+\x_{2}|\leq
1/\sqrt{2}|\m_{1}+\m_{2}|^{\half }\leq C|\m_{1}|^{\half}\leq C|\x_{1}|,
$$
replaces \eqref{estimate} and we can repeat the argument
given in the first part of Case B1. If $\half \gia> |\x_{1}|$, since
$|\m_{1}+\m_{2}|\sim |\m_{1}|$, we dyadically decompose with respect to
$|\m_{1}|\sim 2^{n_{1}}$.
Because $|\x_{1}+\x_{2}|\lesssim |\x_{2}|$, we reduce our estimate to
\begin{equation}
\sum_{j_{1},j_{2}\geq 0}\sum_{j\geq 0}\sum_{n_{1}\geq 0}
2^{-j/2}\int_{A_{2}}g_{j}\chi_{2}\chi_{j}(\x,\m,\tau)
\theta_{n_{1}}
|\x_{1}+\x_{2}|^{2\epsilon_{0}}
 \chi_{2}\phi_{1,n_{1},j_{1}}(\x_{1},\m_{1},\tau_{1})
 \phi_{2,j_{2}}(\x_{2},\m_{2},\tau_{2}),
\label{estimate1}\end{equation}
and at this point we can proceed by using Lemma \ref{mainlemma}.

\noindent
{\bf Region $A_{3}\cup A_{4} \cup (A_{5}-\tilde{A}_{5}( 2 \epsilon_{0}))$}.\\
{\bf Case A:} $ \half \gi\leq|\x|$. \\ We consider two subcases.\\
{\bf Case A1:} $|\xi_{1}|\geq  \half \gia$. \\
We dyadically decompose with respect
to $|\x|\sim|\x_{1}|\sim 2^{m_{1}}$.
Because in this region $\gib\gtrsim |\xi_{1}|\sim |\x_{1}+\x_{2}|$,
we reduce our estimate to \eqref{helpy1} of the corresponding case in
the proof of Theorem \ref{xbil}, where $\epsilon_{0}$ is replaced
by $2\epsilon_{0}$.\\
{\bf Case A2:} $|\xi_{1}|\leq \half \gia$. \\
We use again a dyadic
decomposition with respect to $\x_{1}$. Observe that
$(|\x_{2}|<<)|\x_{1}|\leq \half \gia \leq C\gib$ hence
$\gib \sim 2^{m_{1}+r}, r\geq
-C_{2}$. We then reduce our estimate to \eqref{helpy2}, again replacing
$\epsilon_{0}$ with $2\epsilon_{0}$.\\
{\bf Case B:} $ \half \gi>|\x|$. As we observed in the analysis
of region $A_{i}, i=1,2$, without loss of generality we can assume that
$\gi\geq 1$. We consider two subcases.\\
{\bf Case B1:}  $|\m_{1}|\lesssim |\m_{2}|$. \\
We recall that in this
region we also have $\gib>>|\x_{2}|$ and we repeat the argument presented
in the second part of Case B1 of region $A_{2}$.\\
{\bf Case B2:} $|\m_{1}|>> |\m_{2}|$.\\
If $\half \gia\leq|\x_{1}|$ then
we use an  argument similar to the one used in the first part of
Case B1 in region
$A_{2}$, where the role of $(\x_{2},\m_{2})$ is now played by
$(\x_{1},\m_{1})$.
If $\half \gia \geq |\x_{1}|$,
we dyadically decompose so that
$|\m_{1}+\m_{2}|\sim |\m_{1}|\sim 2^{n_{1}}$. Then because
$|\x_{1}+\x_{2}|\sim |\x_{1}|\lesssim \gib$, we obtain the estimate
\eqref{estimate1} and also this case is done.

\noindent
{\bf Region $A_{5}\cap\tilde{A}_{5}(2\epsilon_{0})$}.
In this case we cannot use Lemma \ref{mainlemma}. \\
{\bf Case A:} $|\x|\geq \half \gi$.\\
We observe that
if \footnote{Based on the proof of Lemma \ref{mainlemma},
 we assume that $\alpha>4$ and $2\epsilon_{0}< \frac{1}{8}$.}
 $\epsilon_{0}<(\alpha -1)/2\alpha$, then
$(1+\alpha 2\epsilon_{0})(1-2\epsilon_{0})>1$ and
\begin{equation}(\gia)^{1-2\epsilon_{0}}\geq |\x_{1}|^{1+\delta}\sim
|\x_{1}+\x_{2}|^{1+\delta},
\label{estimate3}\end{equation}
for some $\delta$ such that $0<\delta<<1$.
We dyadically decompose so that $|\x_{1}|\sim
|\x_{1}+\x_{2}|\sim 2^{m_{1}}$  and using \eqref{estimate3}
we reduce the estimate to
$$
\sum_{j_{1},j_{2}\geq 0}\sum_{j\geq 0}\sum_{m_{1}\geq 0}
2^{-j/2}\int g_{j}\chi_{1}\chi_{j}(\x,\m,\tau)
\theta_{m_{1}}2^{-m_{1}\delta}
 \chi_{1}\phi_{1,j_{1},m_{1}}(\x_{1},\m_{1},\tau_{1})
 \phi_{2,j_{2}}(\x_{2},\m_{2},\tau_{2}),
 $$
and Strichartz inequality \eqref{st1} and Cauchy-Schwarz can be used to
finish the proof also in this case.\\
{\bf Case B:} $|\x|\leq \half\gi$.\\
Notice that from \eqref{m1x1m2x2} $|\m_{2}|\sim |\x_{2}|\gia<<|\m_{1}|$,
hence
$|\m_{1}+\m_{2}|\sim |\m_{1}|$. We dyadically decompose with
respect to
$|\m_{1}|\sim 2^{n_{1}}$. Using \eqref{estimate3} we obtain the
estimate
\begin{eqnarray*}
& &\sum_{j_{1},j_{2}\geq 0}\sum_{j\geq 0}\sum_{n_{1}\geq 0}2^{-j/2}
\int g_{j}\chi_{2}\chi_{j}(\x,\m,\tau)
 \theta_{n_{1}}(\m)
 |\x_{1}|^{-\delta}(\gia)^{1-2\epsilon_{0}}\\
 & &\times
\chi_{2}\phi_{1,n_{1},j_{1}}(\x_{1},\m_{1},\tau_{1})
\frac{\phi_{2,j_{2}}}
{(\gib)^{1-2\epsilon_{0}}}(\x_{2},\m_{2},\tau_{2}).
\end{eqnarray*}
Then Strichartz inequality \eqref{st1} is enough to conclude the
proof also in this case.

\noindent
{\bf Region $A_{6}$}. As we did for the proof of Theorem \ref{xbil},
also in this case we consider three subcases: {\bf Case A:}
$|\x_{1}|\geq 10^{2}\gia$; {\bf Case B:}
$|\x_{1}|< 10^{-2}\gia$ and {\bf Case C:}
$10^{-2}\gia<|\x_{1}|< 10^{2}\gia$. Since in the proof of Theorem \ref{xbil}
we used \eqref{a6mult} in cases A and B, we can treat these cases in the
same way, with the understanding that
now $\epsilon_{0}$ is replaced by $2\epsilon_{0}$. For Case C we go
back to \eqref{j1} and \eqref{j2} and we show that in this region
\begin{equation}
    J^{1}\lesssim J^{2}.
\label{j1<j2}\end{equation}
If one assumes this for a moment, then it is easy to see that we can
repeat exactly the argument we gave for Case C in the proof of
Theorem \ref{xbil}. To prove \eqref{j1<j2} it's enough to show that
$$
    \max(|\x|, \gi)\gtrsim \max(1,|\x_{2}|, \gib).
$$
To simplify the notation we set $\max(|\x|, \gi)=M$ and
$\max(1,|\x_{2}|, \gib)=M_{2}$.  Assume that $M=|\x|$. If
$M_{2}=1$, then
$M_{2}=1\lesssim |\x|=M$.
If  $M_{2}=|\x_{2}|$, then $M\sim
|\x_{1}|>>|\x_{2}|=M_{2}$ and also this case is done.
If  $M_{2}=\gib$, then $M=|\x|\sim |\x_{1}|\sim \gia>>
\gib=M_{2}$. Now assume that  $M=\gi$. If
$M_{2}=1$, then
$M_{2}=1\lesssim |\x|\leq \gi=M$.
If  $M_{2}=|\x_{2}|$, then $M_{2}<<
|\x_{1}|\sim |\x|\leq \gi=M_{2}$ and also this case is done.
Finally If  $M_{2}=\gib$, then $M_{2}<< \gia\sim |\x_{1}|\sim |\x|
\leq \gi=M$.

The proof of Theorem \ref{ybil} is now complete.

\end{proof}

We conclude this section with a counterexample that shows that
if $\epsilon<\frac{1}{4}$ in Theorem \ref{xbil}, then the
theorem does not hold. This counterexample is important because,
as we will discuss below in Remark \ref{rescaling}, if we
could have taken $\epsilon<\frac{1}{4}$, then we could have removed the
smallness assumption in the initial data and at the same time we would
have obtained a global result in the modified energy space
$E\cap P$.

\begin{proposition}\label{caunter}
The bilinear estimate
\begin{eqnarray}
\label{bilinearcounter}\|\partial_{x}(uv)\|_{X_{1,-\half}}&\leq& C
\|u\|_{X_{1,\half}}(\|v\|_{X_{1,\half}}+
\|v\|_{X_{1,\half}}^{1-\epsilon }
\|v\|_{Y_{1,0, \half}}^{\epsilon })\\
\nonumber&+&C
\|v\|_{X_{1,\half }}(\|u\|_{X_{1, \half}}+
\|u\|_{X_{1, \half }}^{1-\epsilon }
\|u\|_{Y_{1,0, \half}}^{\epsilon })
\end{eqnarray}
fails for $\epsilon < \frac{1}{4}$.
\end{proposition}
\begin{proof}
The proof of the proposition is based on the example of Molinet,
Saut and Tzvetkov from \cite{MST, MST1}. We introduce the sets:
\begin{equation*}
  E_1 = \{ (\xi , \mu , \tau): \xi \in [\half\alpha, \alpha],~ \mu \in [
-6 \alpha^2 , 6 \alpha^2 ], ~|\tau - \omega ( \xi , \mu ) |
\lesssim 1
\}
\end{equation*}
\begin{equation*}
E_2 = \{ (\xi , \mu , \tau): \xi \in [N , N +\alpha],~ \mu \in [
\sqrt{3} N^2  , \sqrt{3} N^2 + \alpha^2 ],
~|\tau - \omega ( \xi , \mu  ) | \lesssim 1
\}.
\end{equation*}
where again $\omega (\xi , \mu)= \xi^3 + \frac{\mu^2}{\xi}$.
We observe that $\partial_\xi \omega (\xi, \mu)
= 3 \xi^2 - \frac{\mu^{2}}{\xi^2}$, and $\partial_\mu \omega (\xi , \mu) =
\frac{2\mu}{\xi}.$ The reason for the $\sqrt{3}$ in the definition of
$E_2$
may be seen by calculating $\partial_\xi \omega $ inside
$E_2$. The $O(N^2)$ terms cancel and the next biggest term is
$O(N \alpha )$.
In the region $E_2$, the dispersive
surface has slope $O(N \alpha)$ along the $\xi$ direction and
slope $O(N)$
along $\mu$.
An $\alpha \times \alpha^2$ piece of the tangent plane to the
dispersive
surface in $E_2$ stays within $N \alpha^2 $ of the surface.
Therefore,
we select $\alpha \thicksim N^{-\half}$  so that $E_2$ is
a $O(1)$ vertical thickening of this piece of tangent plane.
These calculations then ``explain''
the choice of the $\alpha$,  $\alpha^2$
scaling in the $\xi$ and $\mu$ directions. Note that $E_1$ is
essentially
an $\alpha \times \alpha^2 \times 1$ box.

Let $\widehat{u}= \alpha^{-\frac{3}{2}} \chi_{E_1}, ~ \widehat{v}
= N^{-1} \alpha^{-\frac{3}{2}} \chi_{E_2},$ where
the functions $\chi_{E_i}$
are smoothed out characteristic functions.

We calculate
\begin{equation*}
  | [\partial_x (u v)]\sphat(\xi, \mu, \tau) | \thicksim |\xi| |\widehat{u} *
\widehat{v} | ( \xi , \mu , \tau ) \thicksim |\xi | \alpha^{-3} N^{-1}
~\chi_{E_1}* \chi_{E_2} (\xi, \mu , \tau ).
\end{equation*}
We have
\begin{equation*}
  \chi_{E_1}* \chi_{E_2} (\xi, \mu , \tau ) \thicksim \sup_{{trans}}
|{{trans}} (E_1) \cap E_2 | ~\chi_{E_1 + E_2 } ( \xi , \mu , \tau ),
\end{equation*}
where ${{trans}}$ denotes an arbitrary translation in the
$(\xi, \mu , \tau)$ space. Geometric considerations similar to those
discussed above show that
$E_1 +E_2 $ contains an $\alpha \times \alpha^2 \times
N \alpha^2 $
box containing the point
$(N+\alpha, \sqrt{3} N^2 + \alpha^2 , 4 N^3 )$.
This point is at vertical distance $O(1)$ from the dispersive surface.
Since $E_2$ has slope
$N \alpha$
along $\xi$ and slope $N$ along $\mu$, we observe that the sup above
is bounded by $(N \alpha)^{-1} \times N^{-1} \times 1$. Combining these
remarks gives,
\begin{equation*}
   | {\widehat{\partial_x u v}} | \thicksim
N ~ \alpha^{-3} N^{-1} ~ N^{-2} \alpha^{-1} \chi_{E_1 + E_2 } \thicksim
\alpha^{-4} N^{-2} \chi_{E_1 + E_2}.
\end{equation*}
Therefore, we have
$$
  \|\partial_{x}(u v)\|_{X_{1,-\half}} \gtrsim  \alpha^{-4} N^{-2}~ N
(N \alpha^2 )^{-\half}~ | E_1 + E_2 |^{\half} \gtrsim  \alpha^{-5}
N^{-1 - \half} ~(\alpha^5 N )^{\half}.
$$
The choice of $\alpha = N^{-\half}$ yields
$$
   \|\partial_{x}(u v)\|_{X_{1,-\half}} \gtrsim  N^{\frac{1}{4}}.
$$
as the size of the left-side of \eqref{bilinearcounter}.

We now consider the right-side of \eqref{bilinearcounter}.
The functions $u,~
v$ are normalized to have size $O(1)$ in the various $X_{s,b}$-norms. The
$Y_{1,0,\half}$-norm has two pieces. The term arising from $\partial_\tau
u$ essentially reproduces $u$ since $E_1$ is of size $O(1)$
along the $\tau$ direction. The other term involves $\partial_\mu u$.
Since $E_1$ has size $\alpha^2$ along $\mu$, we have that
$|\partial_\mu
\chi_{E_2}| \thicksim \alpha^{-2} \chi_{E_2}$ so this part of the
$Y_{1,0,\half}$-norm is of size $O(\alpha^{-2}) = O( N )$. Upon taking
this to the power $\epsilon$ and comparing with the size of the
left-side,
 $N^{\frac{1}{4}}$, we see the failure of \eqref{bilinearcounter}
 when $\epsilon <\frac{1}{4}.$

\end{proof}

\section{Proof of Theorem \ref{energytheorem} and Theorem \ref{main0}  }
We start with the proof of Theorem \ref{main0}, because Theorem
\ref{energytheorem}
is a corollary of Theorem \ref{main0}. The proof uses a classical
fixed point theorem (see for example \cite{KPV3}). We will first carry
out the proof when $T=\half$
\begin{proof}[Proof of Theorem \ref{main0}]
We start by transforming \eqref{ivp} into the
integral equation
\begin{equation}
\label{d} u=\psi(t)S(t)u_{0}-\psi(t)
\int_{0}^{t}
S(t-t')\partial_{x}(u^{2})(t') dt'.
\end{equation}
(Here we have fixed $\beta = 1$.)
Then it is clear that a solution for \eqref{d}
is a fixed point for the operator
\begin{equation}
    \label{oper}L(v)=\psi(t)S(t)u_{0}
    -\psi(t)\int_{0}^{t}
S(t-t')\partial_{x}(v^{2})(t')dt'
    \end{equation}
To simplify the notation we set
$Z_{1-\epsilon_{0}}=X_{1-\epsilon_{0},\half}\cap
Y_{1-\epsilon_{0},- \epsilon_{0}, \half}$. Then, for fixed $\sigma>0$,
we assume that
$\|u_{0}\|_{{B^{1,2}_{1-\epsilon_{0}}}\cap
P^{1,2}_{-\epsilon_{0}}}\leq \sigma$
and we set $a=4C\|u_{0}\|_{B^{1,2}_{1-\epsilon_{0}}\cap
P^{1,2}_{-\epsilon_{0}}}$, where $C$ is the constant in
\eqref{xsgroup} and \eqref{ysgroup}. We show that,
if $B_{a}$ is the ball centered at the origin and radius $a$ in
$Z_{1-\epsilon_{0}}$, then
\begin{eqnarray}
    \label{ball-ball}& & L:B_{a}\longrightarrow B_{a},\\
    \label{contr} & & \|L(u-v)\|_{Z_{1-\epsilon_{0}}}
    \leq \half\|u-v\|_{Z_{1-\epsilon_{0}}}
    \end{eqnarray}
and this is enough to finish the proof of the theorem.
To prove \eqref{ball-ball} we use \eqref{xsgroup}, \eqref{ysgroup},
\eqref{xsnonhom} and \eqref{ysnonhom} to show that
$$
\|L(v)\|_{Z_{1-\epsilon_{0}}}\lesssim \frac{a}{4} +C_{1}
(\|\partial_{x}(v^{2})\|_{X_{1-\epsilon_{0},-\half}}+
\|\partial_{x}(v^{2})\|_{Y_{1-\epsilon_{0},-\epsilon_{0},-\half}})
$$
and if we continue with Theorems \ref{xbil} and \ref{ybil} we obtain
\begin{eqnarray}
\label{smallassumption}
& &\|L(v)\|_{Z_{1-\epsilon_{0}}}\lesssim \frac{a}{4}\\
\nonumber&&+C_{1}(\|v\|_{X_{1-\epsilon_{0},\half}}^{2}+
\|v\|_{X_{1-\epsilon_{0},\half}}
\|v\|_{Y_{1-\epsilon_{0},-\epsilon_{0},\half}}+
\|v\|_{X_{1-\epsilon_{0},\half}}^{1-\epsilon}
\|v\|_{Y_{1-\epsilon_{0},-\epsilon_{0},\half}}^{1+\epsilon}),
\end{eqnarray}
hence \eqref{ball-ball} follows with our choice of $a$, for small
$\sigma$.
\end{proof}

We finish this section with a remark that should convince the reader
that in a sense the fixed point method used above is performed in a
critical regime. This criticality appears in an
unusual way. This remark also shows how to obtain the case of general
T in Theorem \ref{main0}, from the case $T=\half$. One simply chooses
$T=\lambda^{3}/2$ below, $\lambda$ large, and then the norm small
depending also on $\lambda$.

\begin{remark}\label{rescaling}

If $u(x,y,t)$ is a solution of the IVP
\eqref{ivp}, then
\begin{equation}
u_{\lambda}(x,y,t)=\lambda^{2}u(\lambda x,
\lambda^{2}y,\lambda^{3}t)
\label{rescsol}\end{equation}
is a solution for the IVP \eqref{ivp} with initial data
$u_{\lambda,0}(x,y)=\lambda^{2}u_{0}(\lambda x,
\lambda^{2}y)$.
Following the directions in the literature we define the
{\it{critical Sobolev indices for the KP
equation}} as the couple of real numbers  $(s_{c}^{1}, s_{c}^{2})$ such
that the homogeneous Sobolev norm
$$\|u_{\lambda,0}(x,y)\|_{\dot{H}^{s_{c}^{1}, s_{c}^{2}}_{x,y}}\sim C,$$
where $C$ is independent of $\lambda$. With a simple calculation we
obtain
\begin{equation}
    s_{c}^{1}+2s_{c}^{2}=-\half.
\label{critind}\end{equation}
While for the KP-II equation one can get a well-posedness theory for
Sobolev spaces with indices satisfying a relationship pretty close to
\eqref{critind} (see \cite{TT}), for the KP-I, due to the observations
made in the introduction, we do not really expect to be able to reach
near the
critical indices in \eqref{critind}.  The type of criticality that
occurs in our case can be summarized in the following lemma:
\begin{lemma}
In $u_{\lambda}$ is defined as in \eqref{rescsol}, $0<\lambda\leq 1$,
then for any $s\geq 1$,
\begin{equation}
    \|u_{0,\lambda}\|_{B^{2,1}_{s}}^{1-\epsilon}
    \|u_{0,\lambda}\|_{P^{2,1}_{s-1}}^{\epsilon}\lesssim
    \lambda^{1-4\epsilon}.
    \label{badrescal}\end{equation}
    \end{lemma}
The proof of the lemma follows from simple changes of variables.

Observe that we can consider $\epsilon=\frac{1}{4}$ to be a
critical exponent in this case. In
fact, if we could take $\epsilon<\frac{1}{4}$ in Theorems \ref{xbil} and
\ref{ybil}, then we could use \eqref{badrescal} to remove the smallness
assumption needed in the contraction argument presented for the proof
of Theorem \ref{main0}. But,
like the counterexample in Proposition \ref{caunter}
shows, this is not possible.
    \end{remark}
We are now ready to sketch the  proof of  Theorem \ref{energytheorem}.
\begin{proof}[Proof of Theorem \ref{energytheorem}]
% Let's fix an interval of time $[-T,T]$. If $u_{0}\in
% E\cap P$ and $\|u_{0}\|_{E\cap P}\leq \delta_{1}$, then by Remark
% \ref{embedding} we also have that
% $u_{0}\in B^{2,1}_{1-\epsilon}\cap P^{2,1}_{-\epsilon}$ and
% $\|u_{0}\|_{ B^{2,1}_{1-\epsilon}\cap P^{2,1}_{-\epsilon}}\leq \delta,$
% for any positive
% small $\epsilon$. Then by Theorem \ref{main0} there exists a unique
% solution $u\in Z_{1-\epsilon}=X_{1-\epsilon,\half}\cap
% Y_{1-\epsilon,\epsilon, \half}$. Now  by continuity with respect
% to the initial data, if $u^{k}_{0}$ is a smooth sequence that
% approximate $u_{0}$ in $B_{\delta_{1}}\subset E\cap P$, then the
% associated sequence of solutions $u_{k}$ is smooth and approximates
% $u$ in $C([-T,T], B_{1-\epsilon}\cap P_{-\epsilon})$. On the
% other hand, by Remark \ref{natep}
% %
% $$\|u_{k}\|_{L^{\infty}_{[-T,T]}(E\cap P)}\leq C(T,\|u_{0}\|_{E\cap P})$$
% %
% uniformly with respect to $k$. Then $u_{k}$ has a weak limit in
% $L^{\infty}_{[-T,T]}(E\cap P)$ that must coincide with $u$ almost
% everywhere. Then
% %
% $$\|u\|_{L^{\infty}_{[-T,T]}(E\cap P)}\leq C(T,\|u_{0}\|_{E\cap P}),$$
% %
% and this concludes the proof.
% ABOVE WAS REPLACED ON 10 FEB 2003 WITH WHAT FOLLOWS
Let's fix an interval of time $[-T, T]$ and a small $\epsilon > 0 ( \epsilon < \frac{1}{16})$. 
In light of Remark \ref{embedding}, we have that, if $u_0 \in E \cap P$, and $\| u_0 \|_{E\cap P} \leq
\delta_1 $, then $u_0 \in B^{2,1}_{1-\epsilon} \cap P^{2,1}_{-\epsilon}$, and we can choose
$\delta_1 = \delta_1 ( \epsilon )$ so that ${{\| u_0 \|}_{B^{2,1}_{1-\epsilon} \cap P^{2,1}_{- \epsilon}}} \leq \delta$ where $\delta$ is chosen as in Theorem \ref{main0}, depending upon $\epsilon, \delta_1$.
We can thus apply Theorem \ref{main0} to obtain a unique solution $u \in Z_{1-\epsilon} = X_{1-\epsilon , \half} \cap Y_{1-\epsilon, -\epsilon, \half}.$ Now, by continuity with respect to the initial
data, if $u_0^k$ is a smooth sequence that approximates $u_0$ in $B_{\delta_1} \subset E\cap P$, which
also approximates it in $B^{2,1}_{1-\epsilon} \cap P^{2,1}_{-\epsilon}$, then the associated sequence
of solutions $u_k$ is smooth, and it approximates $u$ in $C([-T,T]; B^{2,1}_{1-\epsilon} \cap 
P^{2,1}_{-\epsilon})$. On the other hand, by Remark \ref{natep},
\begin{equation*}
  {{\| u_k \|}_{L^\infty_{[-T,T]} (E \cap P )}} \leq C ( T , {{\| u_0 \|}_{E\cap P}}), 
\end{equation*}
uniformly with respect to $k$. Thus, $u_k$ has a weak limit in $L^\infty_{[-T,T]} (E \cap P)$
that must coincide with $u$ for almost every $t$. Then,
\begin{equation*}
   {{\| u \|}_{L^\infty_{[-T,T]} (E \cap P )}} \leq C ( T , {{\| u_0 \|}_{E\cap P}}),
\end{equation*}
\and this concludes the proof.
\end{proof}


\begin{thebibliography}{10}

\bibitem{BS}
M. Ben-Artzi, J.-C. Saut,
\emph{Uniform decay estimates for a class
of oscillatory integrals and applications},
Diff. Int. Eq. \textbf{12} (1999), 137--145.
%

\bibitem{BL}
J. Bergh, J. L\"ofstr\"om,
\emph{Interpolation Spaces},
Springer-Verlag,  (1976).

\bibitem{B}
J. Bourgain, \emph{On the Cauchy problem for the
Kadomtsev-Petviashvili equation},
Geometric and Funct. Anal. \textbf{3} (1993), 315--341.



\bibitem{C}
J. Colliander,
\emph{Globalizing estimates for the periodic KPI
 equation}, Illinois J. Math.
 \textbf{40} (1996), no. 4, 692--698.


\bibitem{CKS}
J. Colliander, C. Kenig, G. Staffilani,
\emph{Small solutions for the Kadomtsev-Petviashvili I equation}, 
Mosc. Math. J. \textbf{1} (2001), no. 4, 491--520.

\bibitem{F}
A. V. Faminskii, \emph{ Cauchy problem for the generalized Kadomtsev-
Petviashvili equation},
Sibirsk. Mat. Zh.  \textbf{33} (1992), 160--172.

\bibitem{FS}
A. S. Fokas, L. Y. Sung,
\emph{On the solvability of the N-wave,
Davey-Stewartson and Kadomtsev-Petviashvili equations},
Inverse Problems \textbf{8} (1992), 673-708.

\bibitem{GTV:ZS}
J. Ginibre, Y. Tsutsumi, G. Velo,
\emph{On the Cauchy problem for the Zakharov system},
J. Funct. Anal. \textbf{151} (1997), no. 2, 384--436.


\bibitem{IN}
R. J. I\'orio, W. V. L. Nunes,
\emph{On equations of KP-type},
Proc. Royal Soc. Edin.  \textbf{128A} (1998), 725--743.

\bibitem{IMS}
P. J. Isaza, J. L. Mej\'ia, V. Stallbohm,
\emph{Local solutions for  the Kadomtsev-
Petviashvili equation in $\dbR^{2}$},
J. Math. Anal. Appl.   \textbf{196} (1995), 566--587.

\bibitem{IMS1}
P. J. Isaza, J. L. Mej\'ia, V. Stallbohm,
\emph{Regularizing effects for the linearized
   Kadomtsev-Petviashvili (KP) equation}, Rev.
   Colombiana Mat.    \textbf{31} (1997), 37--61.


\bibitem{KP}
B. B. Kadomtsev, V. I. Petviashvili,
\emph{On the stability of solitary waves in weakly dispersive media},
Soviet Phys. Dokl. \textbf{15} (1970), 539--541.


%
\bibitem{KPV}
C. Kenig, G. Ponce,  L. Vega ,
\emph{The Cauchy problem for the
Korteweg-de Vries equation in Sobolev spaces of negative indices},
Duke Math. J.
 \textbf{71}, No.1, (1993), 1--21.
%


\bibitem{KPV1}
C.E. Kenig, G. Ponce, L. Vega,
\emph{Oscillatory integrals and regularity of dispersive equations},
Indiana Univ. Math. J. \textbf{40} (1991), 33-69.



\bibitem{KPV2}
C.E. Kenig, G. Ponce, L. Vega,
\emph{Well-posedness and scattering results for generalized
Korteweg-de Vries equation via the contraction principle},
Comm. Pure Appl. Math. \textbf{46} (1993), 527-620.

\bibitem{KPV3}
C.E. Kenig, G. Ponce, L. Vega,
\emph{The Cauchy problem for the Korteweg-de Vries equation in Sobolev
spaces of negative indices},
Duke  Math. J. \textbf{71} (1993), no. 1, 1-21.


\bibitem{KPV4}
C.E. Kenig, G. Ponce, L. Vega,
\emph{Small solutions to nonlinear Schr\"odinger equations},
Ann. Inst. H. PoincarÈ Anal. Non LinÈaire,
 \textbf{10} (1993), 255--288 .

\bibitem{KPV5}
C.E. Kenig, G. Ponce, L. Vega,
\emph{A bilinear estimate with applications to KdV equations},
J. AMS, \textbf{9} (1996), 573-603 .

\bibitem{MRKP}
L. Molinet, F. Ribaud,
\emph{The global Cauchy problem in Bourgain's
type spaces for a dispersive dissipative semilinear equation},
 SIAM J. Math. Anal.  \textbf{33}  (2002),  no. 6, 1269--1296.

\bibitem{MRKdV}
L. Molinet, F. Ribaud,
\emph{On the low regularity of the Korteweg-de
Vries-Burgers equation}, 
 Int. Math. Res. Not.  \textbf{2002},  no. 37, 1979--2005.


\bibitem{MST}
L. Molinet, J.-C. Saut, N. Tzvetkov,
\emph{Well-posedness and ill-posedness results for the
 Kadomtsev-Petviashvili-I equation},
Duke. Math. J. \textbf{115} (2002), no. 2, 353--384.

\bibitem{MST1}
L. Molinet, J.-C. Saut, N. Tzvetkov,
\emph{Global well-posedness for the KP-I equation},
 Math. Ann.  \textbf{324}  (2002),  no. 2, 255--275.


 \bibitem{S}
J.-C. Saut,
\emph{Remarks on generalized
Kadomtsev-Petviashvili equations},
Indiana Math. J. \textbf{42}, No. 3 (1993), 1011--1026.
%

\bibitem{ST:higherKP}
J.-C. Saut, N. Tzvetkov, 
\emph{The Cauchy problem for higher-order KP equations}, 
J. Differential Equations, \textbf{153} (1999), no. 1, 196--222.



\bibitem{Sw}
M. Schwarz,
\emph{Periodic solutions of   Kadomtsev-
Petviashvili equations},
Adv. Math.  \textbf{66} (1987), 217--233.



\bibitem{Ta1}
H. Takaoka,
\emph{Time local well-posedness for the
Kadomtsev-Petviashvili II equation},
Harmonic Analysis and nonlinear PDE \textbf{1102} (1999), 1--8.
%


\bibitem{TT}
H. Takaoka and N. Tzvetkov,
\emph{On the local regularity of
Kadomtsev-Petviashvili II equation},
Internat. Math. Res. Notices \textbf{2} (2001), 77--114..

\bibitem{T}
M. M. Tom,
\emph{On the  generalized    Kadomtsev-
Petviashvili equation},
Contemp. Math. AMS   \textbf{200} (1996), 193--210.

\bibitem{U}
S. Ukai,
\emph{Local solutions of the  Kadomtsev-
Petviashvili equation},
J. Fac. Sci. Univ. Tokyo Sect. 1A Math.  \textbf{36} (1989), 193-209.


\bibitem{Z}
X. Zhou,
\emph{Inverse scattering transform for the time dependent
Schr\"odinger equations with application to the KP-I equations},
Comm. Math. Phys. \textbf{128} (1990), 551--564.




\end{thebibliography}
\end{document}